\documentclass[notitlepage,amsmath,amssymb,nobibnotes,nofootinbib,longbibliography,showkeys,aps,prfluids,superscriptaddress]{revtex4-2}

\usepackage[T1]{fontenc}
\usepackage{vmargin}
\usepackage{float}
\usepackage{enumitem}
\usepackage[english]{babel}
\usepackage{amssymb}
\usepackage{amsmath,mathtools}
\usepackage{xcolor}
\usepackage{graphicx}
\usepackage{multirow}
\usepackage{overpic}
\usepackage{psfrag}
\usepackage{bm}
\usepackage{epsfig}
\usepackage{marvosym}
\usepackage{natbib}
\usepackage{verbatim}
\usepackage{ulem}
\usepackage{lmodern}
\usepackage{booktabs} 
\usepackage{moresize} 
\usepackage{tikz}
\usepackage{pgfplots}
\pgfplotsset{compat=1.16}
\usetikzlibrary{shapes, arrows, arrows.meta, calc, positioning, decorations,
decorations.pathmorphing, decorations.text}
\usepackage[bookmarks=false]{hyperref}
\usepackage{array}
\usepackage{booktabs}
\usepackage{makecell}
\usepackage{threeparttable}
\usepackage{kantlipsum}
\usepackage{breakcites}

\definecolor{myblue1}	{RGB}{0,177,234}				
\definecolor{myblue2}	{RGB}{76,200,239}				
\definecolor{myblue3}	{RGB}{127,215,244}				
\definecolor{myblue4}	{RGB}{178,231,248}				
\definecolor{myblue5}	{RGB}{198,251,255}				
\definecolor{mybluegray1}{RGB}{0,127,167}				
\definecolor{mybluegray2}{RGB}{76,165,193}				
\definecolor{mybluegray3}{RGB}{127,191,211}				
\definecolor{mybluegray4}{RGB}{178,216,228}				
\definecolor{mygray1}	{RGB}{76,84,93}				
\definecolor{mygray2}	{RGB}{129,135,141}				
\definecolor{mygray3}	{RGB}{165,169,174}				
\definecolor{mygray4}	{RGB}{201,203,206}				
\definecolor{myorange1}	{RGB}{255,126,46}				
\definecolor{myorange2}	{RGB}{255,164,108}				
\definecolor{myorange3}	{RGB}{255,190,150}				
\definecolor{myorange4}	{RGB}{255,216,192}				

\newcommand\red[1]{\textcolor{black}{#1}}

\newcommand\gray[1]{\textcolor{gray}{#1}}

\newcommand\white[1]{\textcolor{white}{#1}}

\newcommand{\bsymb}[1]{\boldsymbol{#1}}
\renewcommand\vec[1]{\bsymb{#1}}
\newcommand\op[1]{{\textbf{#1}}}
\newcommand\myrefeq[1]{$($\ref{#1}$)$}

\def\bal#1\eal{\begin{align}#1\end{align}}
\def\bov#1\eov{\begin{overpic}#1\end{overpic}}

\def\nn{\vec{n}}

\def\uu{\vec{u}}
\def\ff{\vec{f}}

\def\bsigma{\bsymb{\sigma}}

\def\00{\mathbf{0}}
\def\rey{\textrm{Re}}

\DeclareMathOperator{\sgn}{sgn}

\begin{document}
\selectlanguage{english}


\title{Single-step deep reinforcement learning for\break open-loop control of laminar and turbulent flows}



\author{H. Ghraieb}
\affiliation{MINES ParisTech, PSL Research University, Centre de mise en forme des matériaux (CEMEF), CNRS UMR 7635, 06904 Sophia Antipolis Cedex, France}

\author{J. Viquerat}
\affiliation{MINES ParisTech, PSL Research University, Centre de mise en forme des matériaux (CEMEF), CNRS UMR 7635, 06904 Sophia Antipolis Cedex, France}

\author{A. Larcher}
\affiliation{MINES ParisTech, PSL Research University, Centre de mise en forme des matériaux (CEMEF), CNRS UMR 7635, 06904 Sophia Antipolis Cedex, France}

\author{P. Meliga}
\affiliation{MINES ParisTech, PSL Research University, Centre de mise en forme des matériaux (CEMEF), CNRS UMR 7635, 06904 Sophia Antipolis Cedex, France}

\author{E. Hachem}
\affiliation{MINES ParisTech, PSL Research University, Centre de mise en forme des matériaux (CEMEF), CNRS UMR 7635, 06904 Sophia Antipolis Cedex, France}


\date{\today}

\begin{abstract} 
This research gauges the ability of deep reinforcement learning (DRL) techniques to assist the optimization and control of fluid mechanical systems. 
It relies on single-step PPO, a novel, ``degenerate'' version of the proximal policy optimization (PPO) algorithm, intended for situations where the optimal policy to be learnt by a neural network does not depend on state, as is notably the case in open-loop control problems. 
The numerical reward fed to the neural network is computed with an in-house stabilized finite elements environment implementing the variational multiscale (VMS) method.
Several prototypical separated flows in two dimensions are used as testbed. 
\red{The method is applied first to two relatively simple optimization test cases (maximizing the mean lift of a NACA 0012 airfoil and the fluctuating lift of two side-by-side circular cylinders, both in laminar regimes) to assess convergence and accuracy by comparing to in-house DNS data.
The potential of single-step PPO for reliable black-box optimization of computational fluid dynamics (CFD) systems is then showcased by tackling several problems of open-loop control with parameter spaces large enough to dismiss DNS. The approach proves relevant to map the best positions for placement of a small control cylinder in the attempt to reduce drag in laminar and turbulent cylinder flows. 
All results are consistent with in-house data obtained by the adjoint method, and the drag of a square cylinder at Reynolds numbers in the range of a few thousands is reduced by $30\%$, which matches well reference experimental data available from literature. 
The method also successfully reduces the drag of the fluidic pinball, an equilateral triangle arrangement of rotating cylinders immersed in a turbulent stream. 
Consistently with reference machine learning results from the literature, drag is reduced by almost  $60\%$ using a so-called boat tailing actuation made up of a slowly rotating front cylinder and two downstream cylinders rotating in opposite directions so as to reduce the gap flow in between them.}

\end{abstract}

\keywords{Deep Reinforcement Learning; \red{Proximal Policy Optimization}; Neural Networks; Computational fluid dynamics; \red{Open-loop} flow control; Adjoint method}

\hypersetup{pageanchor=false}
\maketitle
\hypersetup{pageanchor=true}


\section{Introduction}\label{sec:intro}

Flow control, defined as the ability to finesse a flow into a more desired state, is a field of tremendous societal and economical importance. In applications such as ocean shipping or airline traffic, reducing the overall drag by just a few percent while maintaining lift can help reducing fossil fuel consumption and CO\textsubscript{2} emission while saving several billion dollars annually~\cite{Corbett2003}. Many other scenario relevant to fluid mechanical systems call for similarly improved engineering design, e.g.,
the airline industry is greatly concerned with reducing the structural vibrations and the radiated noise that occur under unsteady flow conditions~\cite{Khorrami2000,Rowley2002}, while microfluidics~\cite{Knight2002} and combustion~\cite{Syred1974} both benefit from enhanced mixing (which can be achieved by promoting unsteadiness in some appropriate manner). All such problems fall under the purview of this line of study. 

Flow control is often benchmarked in the context of bluff body drag reduction. Numerous strategies have been implemented, either open-loop with passive appendices (e.g., end/splitter plates, small secondary cylinder, flexible tail), or open-loop with actuating devices (e.g., plasma actuation, base bleed, rotation) or closed-loop (e.g. transverse motion, blowing/suction, rotation, all relying on appropriate sensing of flow variables); see the comprehensive surveys of recent developments in~\cite{gad96,luml98,glez02,coll04,kim07,choi08,cork10,catt11,seif18}. Nonetheless, most strategies are trial and error and rely on extensive, costly experimental or numerical campaigns, which has motivated the development of rigorous mathematical formalisms capable of achieving optimal design and control with minimal effort. The adjoint method is one family of such algorithms, that has proven efficient at accurately computing the objective gradient with respect to the control variables in large optimization spaces, and has gained prominence in many applications ranging from atmospheric sciences~\cite{Hall1986} to  aerodynamic design~\cite{Jameson1998,Jameson1998b,gunz02,Mohammadi2004}, by way of fresh developments meant to reshape the linear amplification of flow disturbances~\cite{hill92nasa,gian07,marq08cyl,pral10,sipp10,meli14}.

\red{Another promising option for selecting optimal subsets of control parameters is to rely on machine learning algorithms running labeled data through several layers of artificial neural network while providing some form of corrective feedback.}
Neural networks are a family of versatile non-parametric tools that can learn how to hierarchically extract informative features from data, and have gained traction as effective and efficient computational processors for performing a variety of tasks, from exploratory data analysis to qualitative and quantitative predictive modeling. 
The increased affordability of high performance hardware (together with reduced costs for data acquisition and storage) has indeed allowed leveraging the ever-increasing volume of data generated for research and engineering purposes into novel insight and actionable information, which in turn has reshaped entire scientific disciplines such as image analysis~\cite{Krizhevsky2012} 
or robotics~\cite{kober2013reinforcement,Mnih2015}. Since neural networks have produced most remarkable results when applied to stiff large-scale nonlinear problems~\cite{lusch2018deep}, it is only natural to assume that they can successfully tackle the state-space models arising from the high-dimensional discretization of partial differential equation systems. 
Machine learning has thus been making rapid inroads in fluid mechanics,
with consistent efforts aimed at solving the governing equations~\cite{Raissi2018}, predicting closure terms in turbulence models~\cite{Beck2018}, building reduced-order models~\cite{Wang2018}, controlling flows~\cite{Gautier2015,Raibaudo2020}, or performing flow measurements and visualization~\cite{Lee2017,Rabault2017,Cai2020}; see also~\cite{Brunton2020} for an overview of the current developments in this field.

The focus here is on deep reinforcement learning (DRL), an advanced branch of machine learning
in which \red{deep neural networks learn how to behave in an environment so as to maximize some notion of long-term reward (a task compounded by the fact that each action affects both immediate and future rewards).} 
Several notable works using DRL in mastering games (e.g., Go, Poker) have stood out for attaining super-human level~\cite{silver2017mastering,Moravcik2017}, but the approach has also breakthrough potential for practical applications such as robotics~\cite{Schulman2017,hwangbo2019learning}, computer vision~\cite{Bernstein2018}, self-driving cars~\cite{Pan2017} or finance~\cite{deng2017}, to name a few.
There is also great potential for applying DRL to fluid mechanics, for which efforts are ongoing but still at an early stage, with only a handful of pioneering studies providing insight into the performance improvements to be delivered in shape optimization~\cite{Lee2018,Yan2019,viquerat2019direct} and flow control~\cite{ma2018fluid,Biferale2019,Ren2020}. 
Nonetheless, sustained commitment from the machine learning community has allowed expanding the scope from computationally inexpensive, low-dimensional reductions of the underlying fluid dynamics~\cite{Belus2019,Bucci2019,Novati2019} 
to complex Navier--Stokes systems~\cite{Novati2017,Verma2018}. Proximal policy optimization (PPO~\cite{Schulman2017}) has quickly gained momentum as one of the go-to algorithms for this purpose, 
as evidenced by several recent publications assessing relevance for open- and closed-loop drag reduction in cylinder flows at Reynolds numbers in the range of a few hundreds~\cite{rabault2019artificial,Tang2020,Paris2020,Xu2020}.




This research draws on this foundation to further shape the capabilities of PPO (still a newcomer despite its data efficiency, simplicity of implementation and reliable performance) for flow control, and \red{help} narrow the gap between DRL and advanced numerical methods for multiscale, multi-physics computational fluid dynamics (CFD). The main novelty is the use of single-step PPO, a novel ``degenerate'' algorithm intended for open-loop control problems, \red{as the optimal policy to be learnt is then state-independent, and it may be enough for the neural network to get only one attempt per episode at finding the optimal.} 
\red{The objective is twofold: first, to prove feasibility using several prototypical separated flows in two dimensions as testbed. Second, to assess convergence and  relevance in the context of turbulent flows at moderately large Reynolds number (in the range of a few thousands). This is a topic whose surface is barely scratched by the available literature, as our literature review did not reveal any other study considering DRL-based control of turbulent flows besides~\cite{Ren2020b}, another research effort conducted in the same time frame as the present work.
Single-step PPO has been speculated to hold a high potential as a reliable black-box CFD optimizer~\cite{viquerat2019direct}, but we insist that it lies out of the scope of this paper to provide exhaustive performance comparison data against state-of-the art optimization techniques (e.g., evolution strategies or genetic algorithms). This would indeed require a tremendous amount of time and resources even though
the efforts for developing the method remain at an early stage (to the best of our knowledge, no study in the literature has considered using DRL in a similar fashion) and new algorithms cannot be expected to reach right away the level of performance of their more established counterparts.}

The organization is as follows: section~\ref{sec:method} introduces single-step PPO (together with the baseline principles of DRL and PPO), and outlines the main features of the finite element CFD environment used to compute the numerical reward fed to the neural network. 
Two simple lift optimization problems are presented in section~\ref{sec:opt} to assess convergence and accuracy by comparing against in-house DNS data. 
In section~\ref{sec:olc}, the method is applied to two open-loop drag reduction problems whose parameter spaces are large enough to dismiss DNS, \red{namely the placement of a small control cylinder (for which results computed under laminar and turbulent conditions are compared to in-house data obtained by the adjoint method), and the cylinder rotation of a turbulent fluidic pinball. Finally, in section~\ref{sec:disc}, the method is thoroughly compared (in terms of scope, applicability and performances) to the adjoint method. Evolutionary strategies are also briefly reviewed to put our contribution in perspective and discuss the advantages that may be expected once single-step PPO is finely tuned and characterized.}


\section{Methodology}\label{sec:method}

\subsection{Deep reinforcement learning}\label{sec:method:sub:drl}

Reinforcement learning (RL) provides a consistent framework for modeling and solving decision-making problems through repeated interaction between an agent and an environment. 
We consider the standard formulation in which the agent takes an action $a_t$ based on a partial observation of the current state $s_t$ the environment is in. The environment transits to the next state $s_{t+1}$, and the agent is fed with a reward $r_t$ that acts as the quality assessment of the actions recently taken. This repeats until some termination state is reached, the objective of the agent being to determine the succession of actions maximizing its cumulative reward over an episode (\red{this is} the reference unit for agent update, best understood as one instance of the scenario in which it takes actions). 
Deep reinforcement learning (DRL) \red{combines RL and} deep neural networks, i.e., collections of connected units or artificial neurons, that can be trained to arbitrarily well approximate the mapping function between input and output spaces.
We consider here fully connected networks in which neurons are stacked in layers 
and information propagates forward from the input layer to the output layer via ``hidden'' layers. 
Each neuron performs a weighted sum of its inputs to assign significance with regard to the task the algorithm is trying to learn, adds a bias to figure out the part of the output independent of the input, and feeds an activation function that determines whether and to what extent the computed value should affect the outcome.

\subsection{Proximal policy optimization}\label{sec:method:sub:ppo}

Proximal policy optimization (PPO)~\cite{Schulman2017} is a model free, on-policy gradient, advantage actor-critic reinforcement algorithm. The related key concepts can be summarized as follows:
\smallskip

\paragraph*{- \textit{model free:\!\!\!}} the agent interacts with the environment itself, not with a surrogate model of the environment \red{(the corollary here being that it needs no assumptions about the fluid dynamics underlying the control problems to be solved). }
\smallskip

\paragraph*{- \textit{policy gradient:\!\!\!}} the behavior of the agent is entirely defined 
by a probability distribution $\pi(s,a)$ over actions given states, optimized by gradient ascent.
In DRL, the policy is represented by a neural \red{network. The free
parameters learnt from data are the network weights and biases,} with respect to which the gradient is computed backwards from the output to the input layer according to the chain rule, one layer at the time, using the back-propagation algorithm~\cite{Rumelhart1986}. 
\smallskip

\paragraph*{- \textit{on-policy:\!\!\!}} the algorithm improves the policy used to generate the training data (in contrast to off-policy methods that also learn from data generated with other policies). 
\smallskip


\paragraph*{- \textit{advantage:\!\!\!}} the policy gradient is approximated by that of the policy loss
\bal
\label{eq:gradient_loss}
\mathbb{E}_{\tau\sim\pi}\left[ \sum_{t=0}^T \log \left( \pi (a_t \vert s_t) \right) \widehat{A}^{\pi}(s,a)\right]\,,
\eal  
where $\tau=(s_0, a_0, \dots  , s_T, a_T)$ is a trajectory of state and actions with horizon $T$, $A^{\pi}$ is the advantage function measuring the gain associated with taking action $a$ in state $s$, compared to taking the average over all possible actions, and $\widehat{A}^{\pi}$ is some biased estimator of the advantage, here
its normalization to zero mean and unit variance.
\smallskip

\paragraph*{- \textit{actor-critic:\!\!\!}}  the learning performance is improved by updating two different networks, a first one called actor that controls the actions taken by the agent, and a second one called critic, that estimates the advantage as
\bal
\label{eq:advdef}
A^{\pi}(s_t,a_t)=r_t + \gamma V(s_{t+1})-V(s_t)\,,
\eal
where $V(s)$ is the expected value of the return of the policy in state $s$ and $\gamma\in[0,1]$ is a discount factor adjusting the trade-off between immediate and future rewards.
\smallskip

PPO uses conservative policy updates 
to alleviate the issue of performance collapse affecting standard policy gradient implementations\footnote{Large policy updates can cause the agent to fall off the cliff and to restart from a poorly performing state with a locally bad policy, which is all the more harmful as the step size for policy updating cannot be tuned locally (an above average value can speed up learning in regions of the parameter space where the policy loss is relatively flat, but trigger exploding updates in sharper variation regions).}. 
We use here PPO-clip\footnote{As opposed to PPO-Penalty, a variant relying on a penalization on the average Kullback--Leibler divergence between the current and new policies, but that tends to perform less well in practice.} to optimize the surrogate loss
\bal
\label{eq:ppo_loss}
{\mathbb{E}}_{(s,a) \sim \pi} \left[ \min \left( \frac{\pi (a \vert s)}{\pi_{old} (a \vert s)} ,  1+\epsilon\sgn{(\widehat{A}^{\pi}(s,a))}\right) \widehat{A}^{\pi} (s,a)\right]\,,
\eal
where $\epsilon$ is the clipping range defining how far away the new policy is allowed to go from the old.
The general picture is that a positive (resp. negative) advantage increases (resp. decreases) the probability of taking action $a$ in state $s$, but always by a proportion smaller than $\epsilon$, otherwise the min kicks in~\myrefeq{eq:ppo_loss} and its argument hits a ceiling of $1+\epsilon$ (resp. a floor of $1-\epsilon$). 
This prevents stepping too far away from the current policy, and ensures that the new policy will
behave similarly. 
There exist more sophisticated PPO algorithms (e.g., Trust region PPO~\cite{Wang2019}, that determines first a maximum step size relevant for exploration, then adaptively adjusts the clipping range to
find the optimal within this trust region), but standard PPO has simple and effective heuristics. It is computationally inexpensive, easy to implement (as it involves only the first-order gradient of the policy log probability), and remains regarded as one of the most successful RL algorithms, achieving state-of-the-art performance across a wide range of challenging tasks, including flow control~\cite{rabault2019artificial}.



\subsection{Single-step PPO}\label{sec:method:sub:singleppo}

\begin{figure}[!t]
\centering
\begin{tikzpicture}[	every loop/.style={	min distance=20mm,looseness=20},
				fontsize/.style={		font=\footnotesize},
				intnode/.style={		black, fontsize, pos=0.5},
				arrow/.style={		thick,color=mybluegray3, rounded corners,thick,-stealth},
				box/.style={		rectangle,rounded corners, draw=mygray1, very thick, 
								text width=1.5cm,minimum height=1cm,text centered,
								inner sep=2pt, outer sep=0pt, fontsize},
				backbox/.style={	box, opacity=0.5,text width=8cm,minimum height=4.5cm}]

	\node[backbox]			(bb)		at (6,0) {};
	\node[box,fill=mygray4] 	(s0) 		at (0,0) {$s_0$};
	\node[box,fill=myorange4] (agent) 	at (4,0) {Agent};
	\node[box,fill=myblue4] 	(env) 	at (8,0) {Parallel envs.};

	\draw[arrow] 	(s0) 			to 	[out=0,in=180] 						 								           (agent.west);
	\draw[arrow] 	(agent) 		to 	[out=0,in=180] 				node[intnode, above] {$a_t$} 				           (env.west);
	\draw[arrow] 	(env.south) 	to 	[out=-90,in=-90] 			node[intnode, below] {$r_t$} 						   (agent.south);
	\draw[] 		(agent.north) 	edge	[out=150,in=30, loop,arrow] 	node[intnode, above] {$\,\,\theta_t \to \theta_{t+1}$} (agent.north);
\end{tikzpicture}
\caption{Action loop for single-step PPO. At each episode, the same input state $s_0$ is provided to the agent, which in turn provides $n$ actions to $n$ parallel environments. The latter return $n$ rewards, that evaluate the quality of each action taken. Once all the rewards are collected, an update of the agent parameters is made using the PPO loss~\myrefeq{eq:ppo_loss}. The process is repeated until convergence.
}
\label{fig:single_step_PPO}
\end{figure}
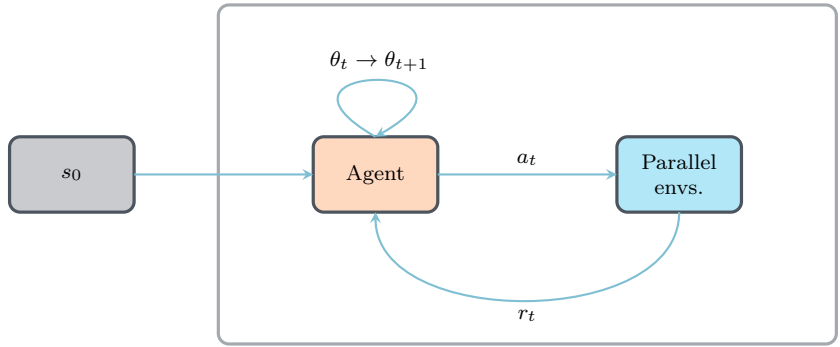

We now come to single-step PPO (hereafter denoted by PPO-1 to ease the reading), a ``degenerate'' version of PPO \red{introduced in~\cite{viquerat2019direct} and} intended for situations where the optimal policy to be learnt by the neural network is state-independent, as is notably the case in open-loop control problems (closed-loop control problems conversely require state-dependent policies for which standard PPO is best suited). The main difference between \red{standard and single-step} PPO can be summed up as follows: where standard PPO seeks the optimal set of actions $a_{opt}$ yielding the largest possible reward, single-step PPO seeks the optimal mapping $f_{\theta_{opt}}$ such that $a_{opt} = f_{\theta_{opt}} (s_0)$, where $\theta$ denotes the network free parameters and $s_0$ is some input state (usually a vector of zeros) consistently fed to the agent for the optimal policy to eventually embody the transformation from $s_0$ to $a_{opt}$. The agent initially implements a random state-action mapping $f_{\theta_0}$ from $s_0$ to an initial policy determined by the free parameters initialization $\theta_0$, after which it gets only one attempt per learning episode at finding the optimal (i.e., it  interacts with the environment only once per episode). This is illustrated in figure~\ref{fig:single_step_PPO} showing the agent draw a population of actions $a_t= f_{\theta_t}(s_0)$ from the current policy, and being returned incentives from the associated rewards to update the free parameters for the next population of actions $a_{t+1} = f_{\theta_{t+1}}(s_0)$ to yield larger rewards.

\red{In practice, the agent outputs a policy} parameterized by the mean and variance of the probability density function of a $d$-dimensional multivariate normal distribution, with $d$ the dimension of the action required by the environment. Actions drawn in $[-1,1]^d$ are then mapped into relevant physical ranges, a step deferred to the environment as being problem-specific. The resolution essentially follows the process described in section~\ref{sec:method:sub:ppo}, only the surrogate loss reads
\bal
\label{eq:ppo_loss_single}
{\mathbb{E}}_{a \sim \pi} \left[ \min \left( \frac{\pi (a)}{\pi_{old} (a)} ,  1+\epsilon\sgn{(\widehat{A}^{\pi}(a))}\right) \widehat{A}^{\pi} (a)\right]\,,
\eal
and the advantage $A^{\pi}$ reduces to the whitened reward $r_t$. 
\red{This is because the trajectory consists of a single state-action pair, so the discount factor can be set to $\gamma=1$ with no loss of generality. In return, the two rightmost terms cancel each other out in~\myrefeq{eq:advdef}, meaning that single-step PPO can do without the value-function evaluations of the critic network (and is thus not actually actor-critic).}

\subsection{Computational fluid dynamics environment} \label{sec:method:sub:cfd}

The CFD resolution framework relies on the in-house, parallel, finite element library CimLIB\_CFD~\cite{coup13cmame}, whose main ingredients are as follows:
\smallskip

\paragraph*{-\!\!\!} the variational multiscale approach (VMS) is used to solve a stabilized weak form of the governing equations using linear approximations (P$_1$ elements) for all variables, which otherwise breaks the Babuska--Brezzi condition. The approach relies on an a priori decomposition of the solution into coarse and fine scale components~\cite{hugh98,codi00,bazi07}. Only the large scales are fully represented and resolved at the discrete level. The effect of the small scales is encompassed by consistently derived source terms proportional to the residual of the resolved scale solution, hence ad-hoc stabilization parameters comparable to local coefficients of proportionality. 
\smallskip

\paragraph*{-\!\!\!} in laminar regimes, velocity and pressure come as solutions to the Navier--Stokes equations. In turbulent regimes, the focus is on phase-averaged velocity and pressure modeled after the unsteady Reynolds averaged Navier--Stokes (uRANS) equations. In order to avoid
transient negative turbulent viscosities, negative Spalart--Allmaras~\cite{Allmaras2012} is used as turbulence model, whose stabilization proceeds from that of the convection-diffusion-reaction equation~\cite{codina1998comparison,badia2006analysis}. 
\smallskip

\paragraph*{- \!\!\!} the immersed volume method (IVM) is used to immerse and represent all geometries inside a unique mesh. The approach combines level-set functions 
to localize the solid/fluid interface, and anisotropic mesh adaptation to refine the mesh interface 
under the constraint of a fixed, number of edges. This ensures that the quality of all actions taken over the course of a PPO optimization is equally assessed, even though the interface can depend on the action. 
\smallskip




Substantial evidence of the flexibility, accuracy and reliability of this numerical framework is documented in several papers to which the reader is referred for exhaustive details regarding the level-set and mesh adaptation algorithms~\cite{Bruchon2004,gruau20053d}, the VMS formulations, stabilization parameters and discretization schemes used in laminar and turbulent regimes~\cite{hach10,coup13jcp,sari2018,guiza2020}, and the mathematical formulation of the IVM in the context of finite element VMS methods~\cite{hachem2012immersed,hach13}.

\subsection{Numerical implementation}\label{sec:method:sub:setup}

\begin{figure}[!t]
\setlength{\unitlength}{1cm}
\begin{picture}(20,4.5)
\put(0.1,-0.75){\includegraphics[trim=175 92.5 155 340pt,clip,height=5.1cm]{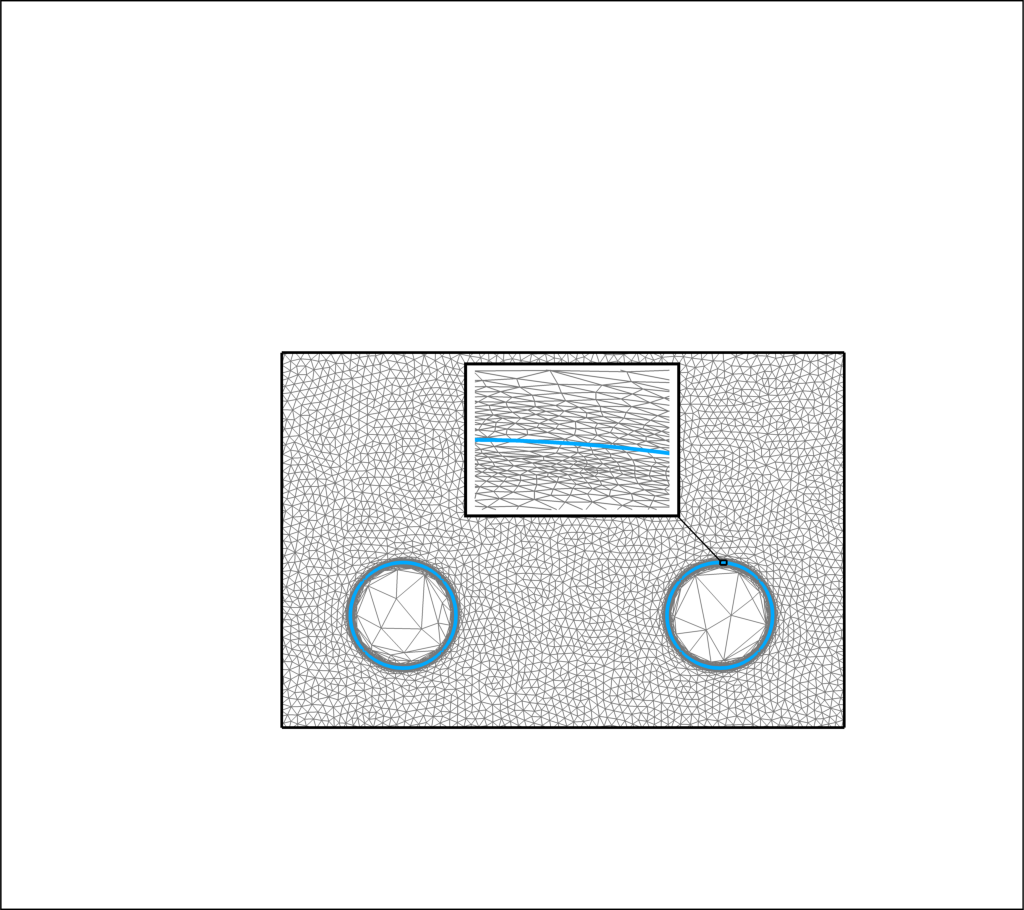}}
\put(7.65,-0.75){\includegraphics[trim=175 92.5 155 340pt,clip,height=5.1cm]{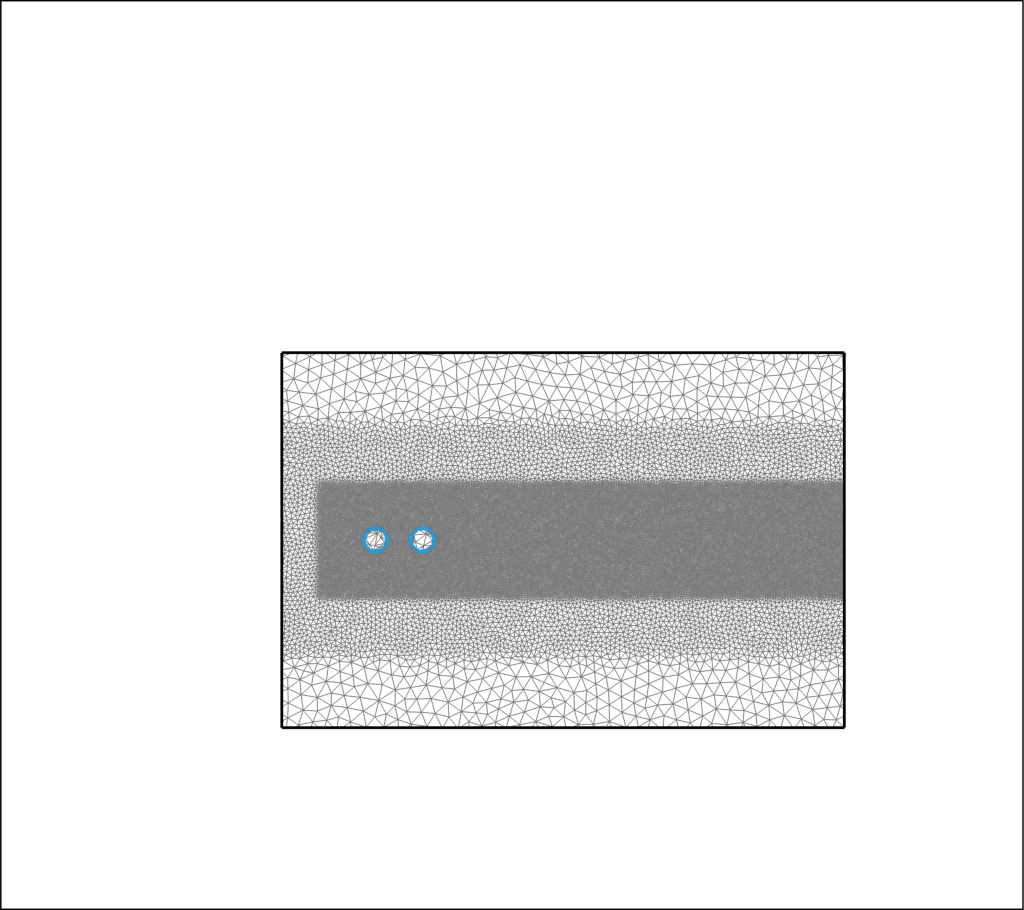}}
\put(0.1,4.5){(a)}
\put(7.65,4.5){(b)}
\end{picture}
\caption{Details of (a) the boundary layer mesh
and (b) successive refinement steps of the background mesh for the flow past a tandem arrangement of two circular cylinders. The blue line in (a) indicates the zero iso-contour of the level set function.}
\label{fig:mesh}
\end{figure}

In practice, actions are distributed to multiple environments running in parallel, each of which executes a self-contained MPI-parallel CFD simulation and feeds data to the DRL algorithm  (hence, two levels of parallelism related to the environment and the computing architecture). \red{Here, all CFD simulations are performed on 12 cores of a workstation of Intel Xeon E5-2640 processors. The algorithm waits for the simulations running in all parallel environments to be completed, then shuffles and splits the rewards data set collected from all environments into several buffers (or mini-batches) used sequentially to compute the loss and perform a network update. }
The process repeats for several epochs, i.e., several full passes of the training algorithm over the entire data set, which ultimately makes the algorithm slightly off-policy (since the policy network ends up being trained on samples generated by older policies, which is customary in standard PPO operation).
This simple parallelization technique is key to use DRL in the context of CFD applications, as a sufficient number of actions drawn from the current policy must be evaluated to accurately estimate the policy gradient.
\red{This comes at the expense of computing the same amount of reward evaluations, and yields a substantial computational cost for high-dimensional fluid dynamics problems (\red{typically from 
a few to several hundred CFD simulations} for the cases considered herein)}. In the same vein, it should be noted that the common practice in DRL studies to gain insight into the performances of the selected algorithm by averaging results over multiple independent training runs with different random seeds is not tractable, as it would trigger a prohibitively large computational burden. The same random seeds have thus been deliberately used over the whole course of study to ensure a minimal level of performance comparison between cases.
The remainder of the practical implementation details are as follows:
\smallskip


\paragraph*{- \!\!\!} the environment 
consists of CFD simulations of two-dimensional (2-D) flows described in a Cartesian coordinate system with drag positive in the $+x$ direction. All equations are discretized on rectangular grids whose side lengths documented in the coming sections have been checked to be large enough not to have a discernible influence on the results (with the exception of the square cylinder flow in section~\ref{sec:olc:sub:control3900}  and the fluidic pinball in section~\ref{sec:olc:sub:pinball}, for which we use respectively the values recommended in~\cite{Rodi1997} and the same values as in~\cite{Raibaudo2020}). Open flow conditions are used, that consist of a uniform inflow in the $x$ direction, together with symmetric lateral, advective outflow and no-slip interface conditions. In turbulent regime, the ambient value of the Spalart--Allmaras variable is three times the molecular viscosity, as recommended to lead to immediate transition. Typical adapted meshes of the interface and wake regions are shown in figure~\ref{fig:mesh}, the latter also being accurately captured via successive refinement of the background elements. 
\smallskip

\begin{table}[!t]
\begin{center}
\begin{tabular}{cc}
\toprule
\multicolumn{2}{r}{Neural network} \\
\cmidrule(lr){1-2}
\multicolumn{1}{r}{2} & \multicolumn{1}{l}{\qquad Nb. hidden layers} \\
\multicolumn{1}{r}{4} & \multicolumn{1}{l}{\qquad Nb. neurons/layer} \\
\multicolumn{1}{r}{TBS} & \multicolumn{1}{l}{\qquad Nb. epochs} \\
\multicolumn{1}{r}{TBS} & \multicolumn{1}{l}{\qquad Nb. environments} \\
\multicolumn{1}{r}{TBS} & \multicolumn{1}{l}{\qquad Size of mini-batches} \\
\cmidrule(lr){1-2}
\multicolumn{2}{r}{PPO} \\
\cmidrule(lr){1-2}
\multicolumn{1}{r}{$5\times10^{-3}$} & \multicolumn{1}{l}{\qquad Learning rate} \\
\multicolumn{1}{r}{0.3} & \multicolumn{1}{l}{\qquad Clipping range} \\
\multicolumn{1}{r}{1} & \multicolumn{1}{l}{\qquad Discount factor} \\
\bottomrule
\end{tabular}
\caption{\label{tab:ppo} Details of the network architecture and PPO hyper parameters. The number of epochs, environments and the size of the mini-batches are provided on a case-by-case basis in sections~\ref{sec:opt} and \ref{sec:olc}.}
\end{center}
\end{table}

\paragraph*{- \!\!\!} the instant reward is (up to a plus/minus sign) either the time-averaged or the root mean square (rms) value of the force coefficient (drag or lift per unit span length), to consider either the 
mean or fluctuating force acting on the immersed body. Instantaneous values are computed with a variational approach featuring only volume integral terms, reportedly less sensitive to the approximation of the body interface than their surface counterparts~\cite{John1997,John2004}.
Time averages are performed over an interval $[t_i;t_f]$ with
edges large enough to dismiss the initial transient and achieve convergence to statistical equilibrium. Moving average rewards and actions are also computed as the sliding average over the $50$ latest values (or the whole sample if it has insufficient size). 
\smallskip

\paragraph*{- \!\!\!} the agent is a fully connected network with 2 hidden layers, each of which holds 4 neurons with hyperbolic tangent activation functions. We use the default online PPO implementation of Stable Baselines, 
a toolset of reinforcement learning algorithms dedicated to the research community and industry~\cite{stable-baselines}, for which a custom OpenAI environment has been designed with the Gym library~\cite{1606.01540}. 
\red{Unlike other RL algorithms, 
PPO does not generally require significant tuning of the hyper parameters (i.e.,  parameters that are not estimated from data). Nonetheless, all values used in this study are documented in table~\ref{tab:ppo} to ease reproducibility, including} the learning rate (the size of the step taken in the gradient direction for policy update), the PPO clipping range (set to the upper edge of the recommended range) 
and the discount factor (set to the default PPO-1 value).


\section{Application to flow optimization} \label{sec:opt} 


\subsection{Flow past a NACA 0012 airfoil} \label{sec:opt:sub:naca}

\begin{figure}[!t]
\setlength{\unitlength}{1cm}
\begin{picture}(20,10.5)
\put(0.1,5.25){\includegraphics[trim=175 92.5 155 340pt,clip,height=5.1cm]{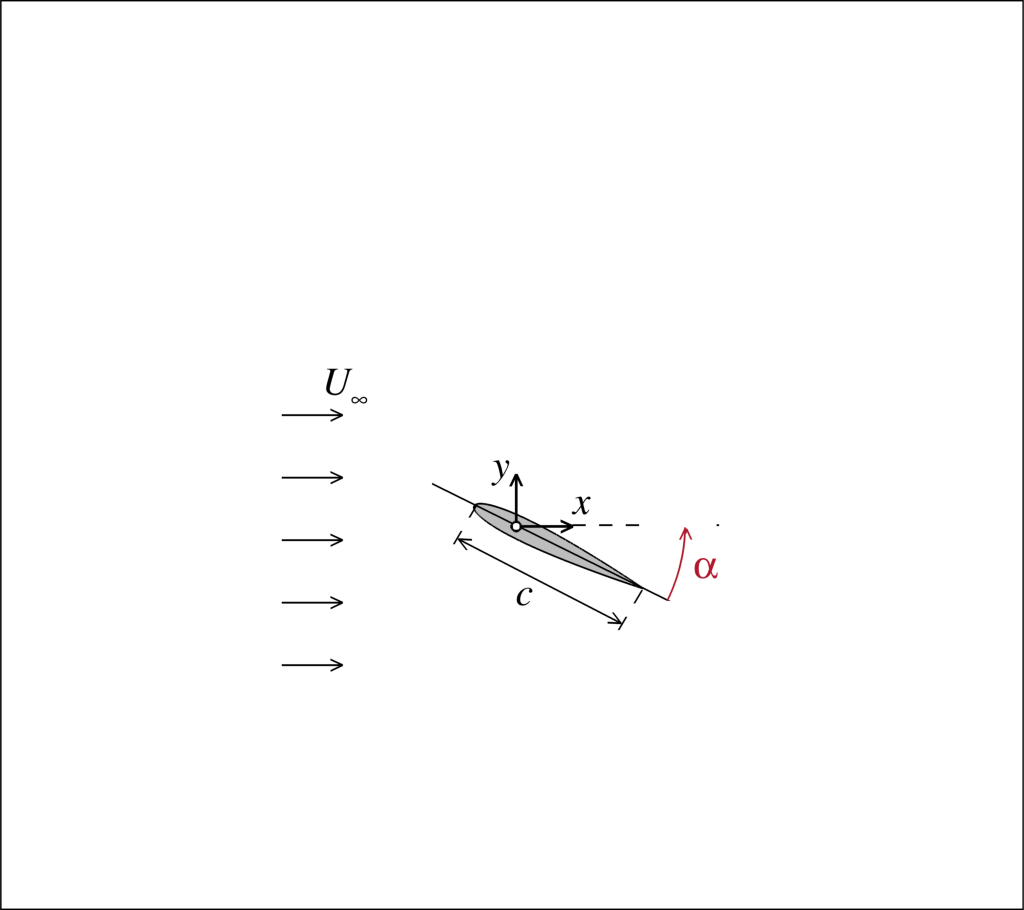}}
\put(7.65,5.25){\includegraphics[trim=175 92.5 155 340pt,clip,height=5.1cm]{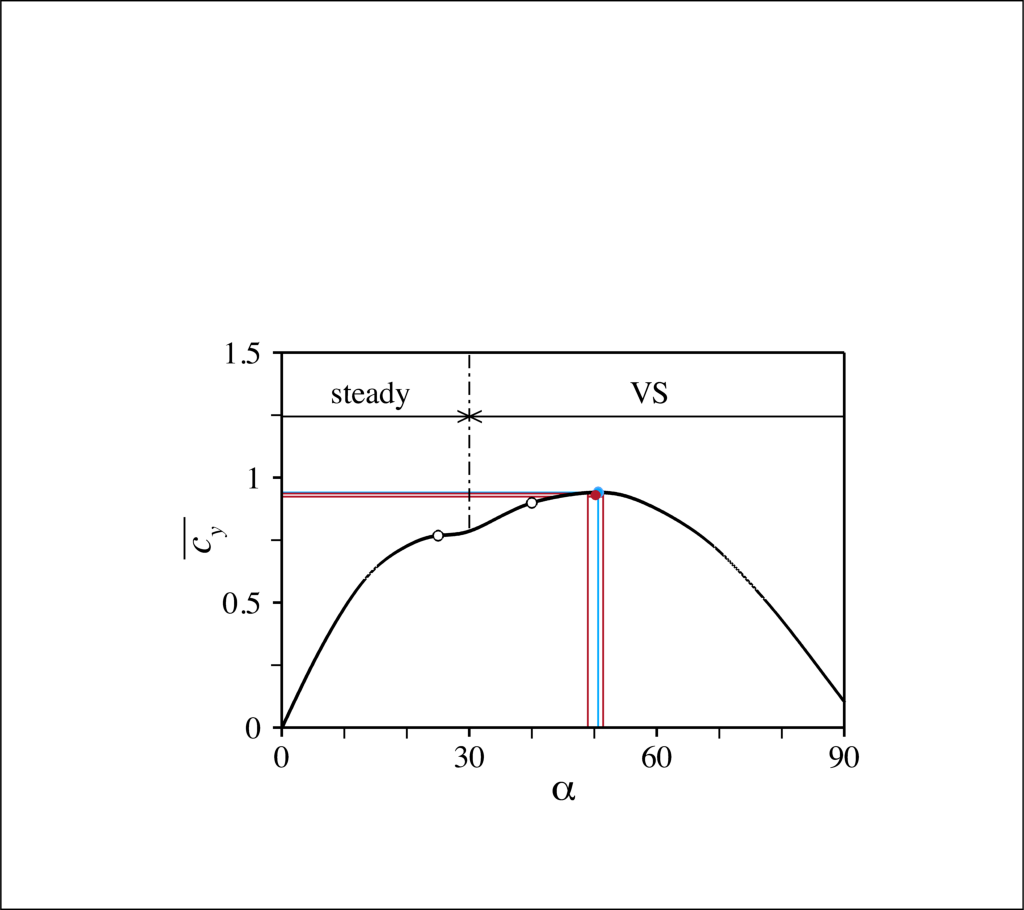}}
\put(0.1,0){\includegraphics[trim=175 92.5 155 340pt,clip,height=5.1cm]{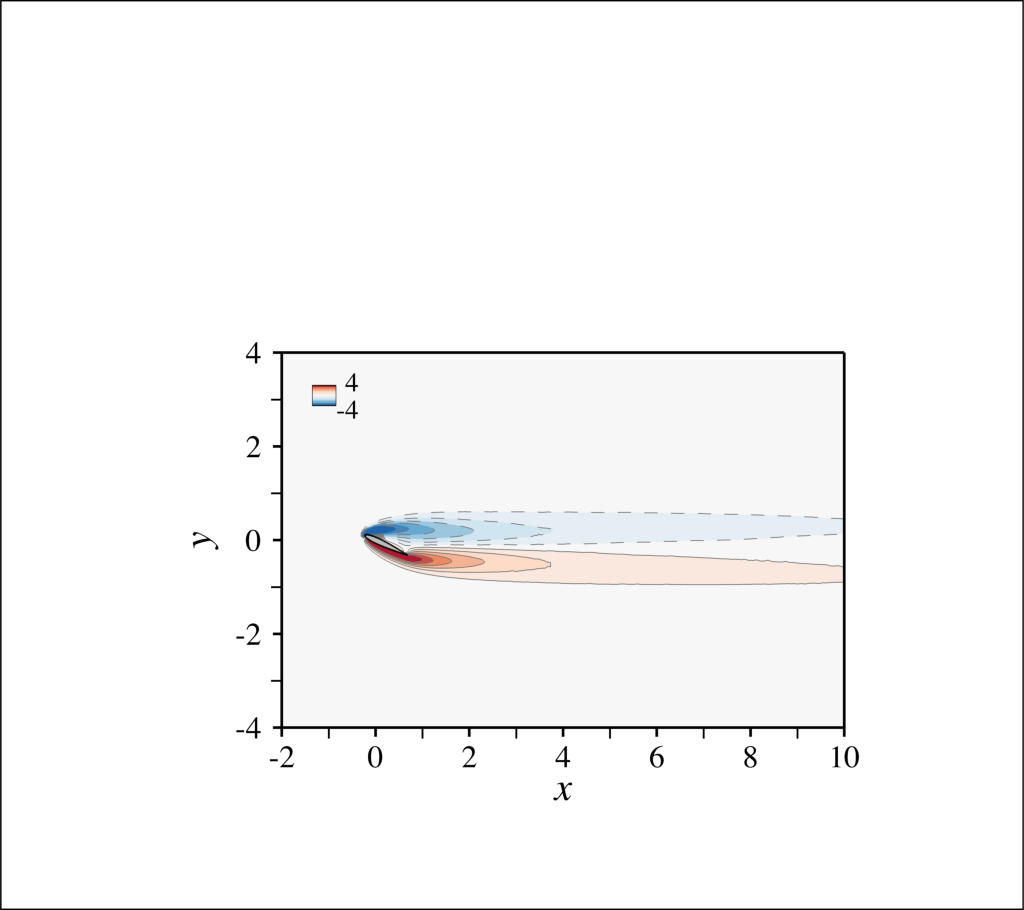}}
\put(7.65,0){\includegraphics[trim=175 92.5 155 340pt,clip,height=5.1cm]{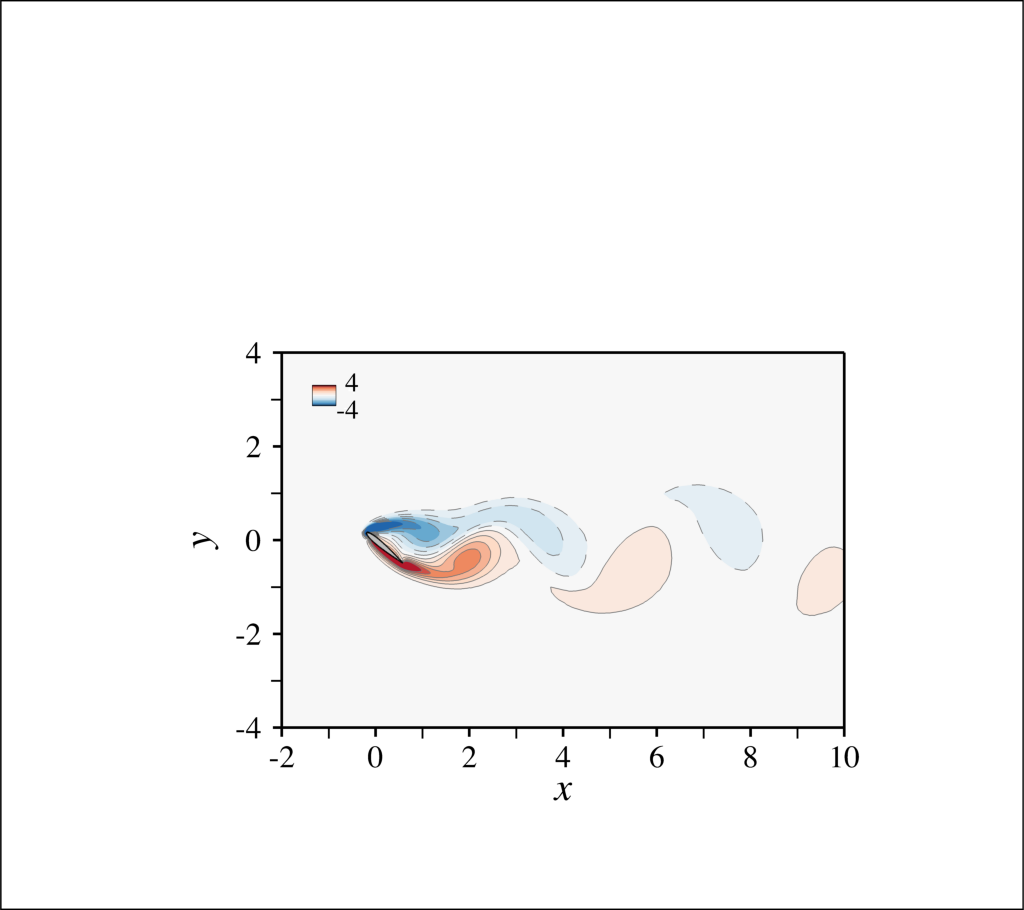}}
\put(0.1,10.5){(a)}
\put(7.65,10.5){(b)}
\put(0.1,5.25){(c)}
\put(7.65,5.25){(d)}
\end{picture}
\caption{Flow past a NACA 0012 - (a) Schematic diagram of the configuration. (b) Mean lift against the angle of attack computed by DNS at $\rey=100$. The VS label indicates the angles for which the flow exhibits unsteadiness in the form of periodic vortex formation and shedding. The blue lines and symbols mark the optimal. The red symbol is the average over the 5 latest single-step PPO episodes, and the red lines 
delimit the corresponding variance intervals. (c-d) Instantaneous vorticity fields computed at $\rey=100$, for values marked by the circle symbols in (b), namely (c) $\alpha=25$ and (d) $\alpha=40$.}
\label{fig:naca}
\end{figure}
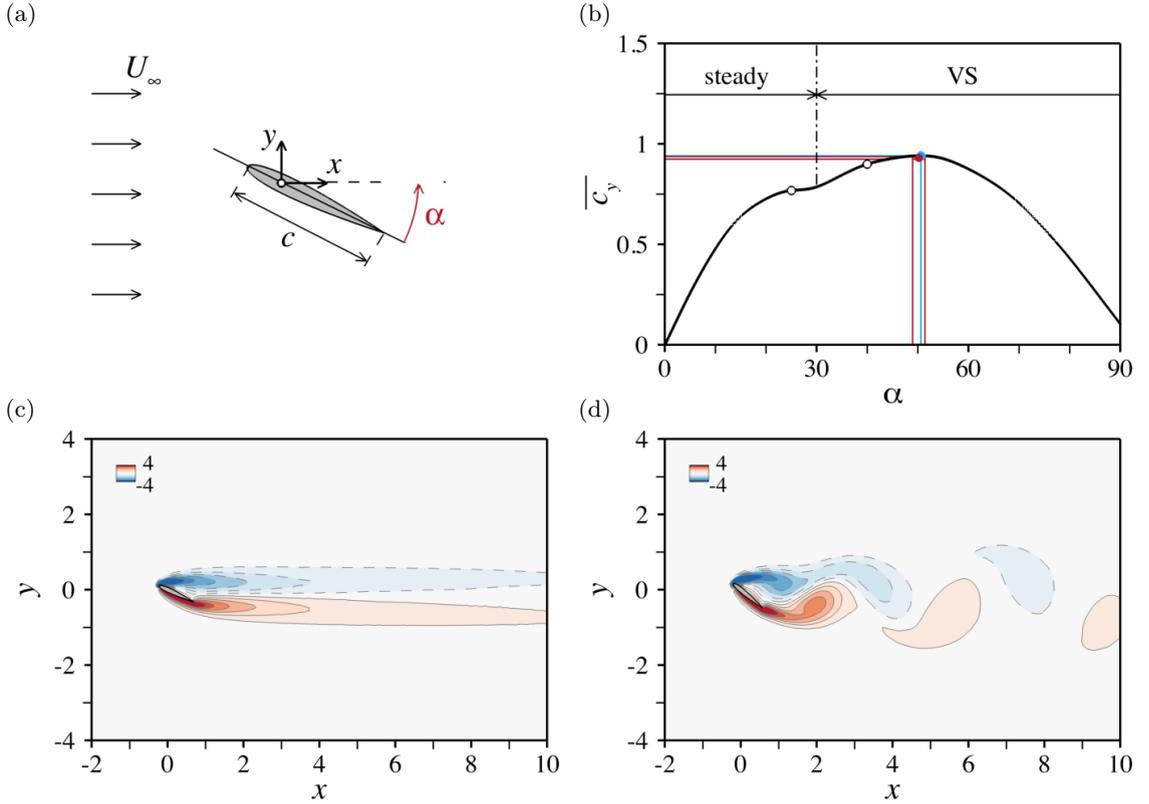

We consider first a NACA 0012 airfoil placed at incidence in a uniform stream, as depicted in figure~\ref{fig:naca}(a). The origin of the coordinate system is at the airfoil pivot-point, set at quarter chord length from the leading edge. 
A laminar, time-dependent case at Reynolds number $\rey=U_\infty c/\nu=100$ 
is modeled after the Navier--Stokes equations, where $U_\infty$ is the inflow velocity, $c$ the straight chord distance and $\nu$ the kinematic viscosity. The objective is to maximize the mean lift $\overline{c_y}$, for which the sole control parameter is the angle of attack $\alpha$ measuring the incidence relative to the chord \red{(in degrees and with the convention that $\alpha> 0$ 
for the airfoil to generate positive lift. Also, we keep in mind that $\alpha$ is rather a state parameter than an adjustable control parameter in practical situations}, but the methodology carries over to related optimization problems such as the design of multi-element high-lift systems). This is a problem simple enough to allow direct comparisons between PPO-1 and DNS (actually VMS, but the difference is clear from context), all the more so as lift varies smoothly with the incidence. This is evidenced in figure~\ref{fig:naca}(b) showing reference data obtained from $15$ DNS runs computing the mean lift to an accuracy of $3\%$ with the simulation parameters documented in table~\ref{tab:opt}. 
The distribution changes slope near $\alpha\sim 30^\circ$ (because the system bifurcates from a steady to a time-periodic vortex-shedding regime; see figure~\ref{fig:naca}(c-d) showing instantaneous vorticity fields computed on either side of the threshold) but otherwise exhibits a well-defined, smooth maximum at $\alpha^\star=50.6$, associated with $\overline{c_y}^{\,\star}=0.94$.

\begin{figure}[!t]
\setlength{\unitlength}{1cm}
\begin{picture}(20,6)
\put(0.1,0){\includegraphics[trim=175 92.5 155 270pt,clip,height=5.85cm]{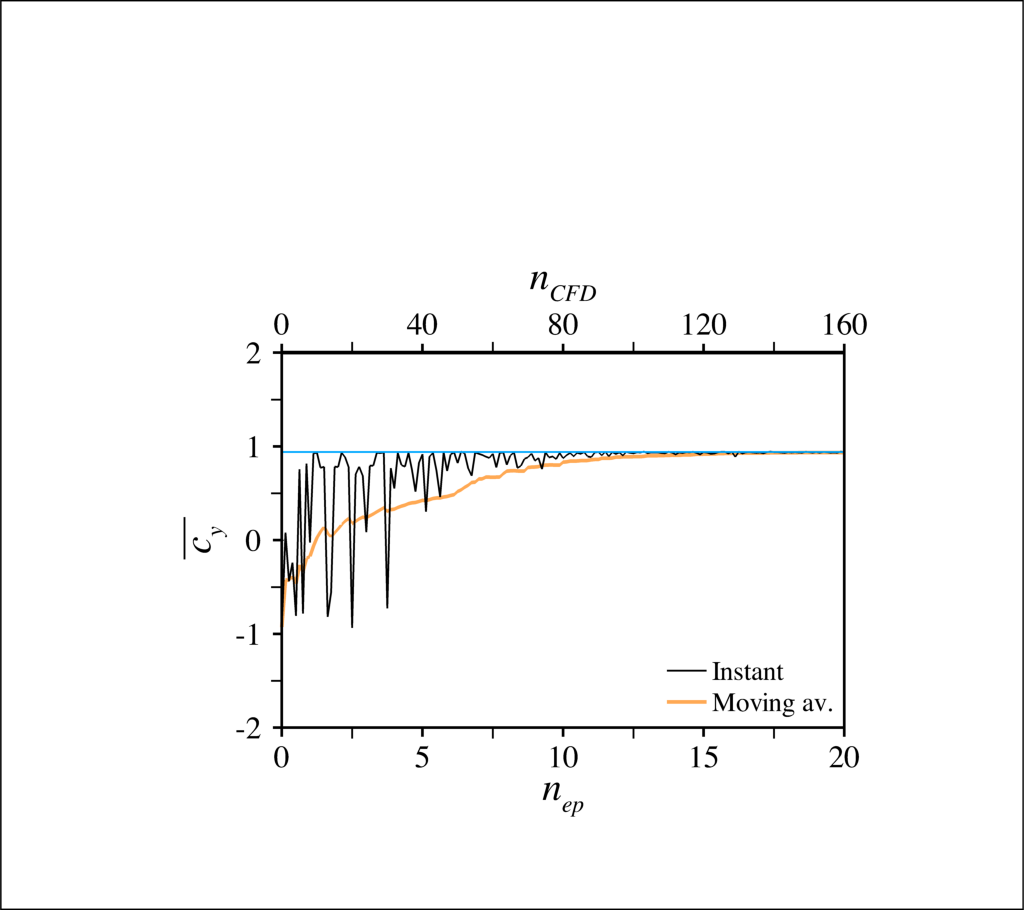}}
\put(7.65,0){\includegraphics[trim=175 92.5 155 270pt,clip,height=5.85cm]{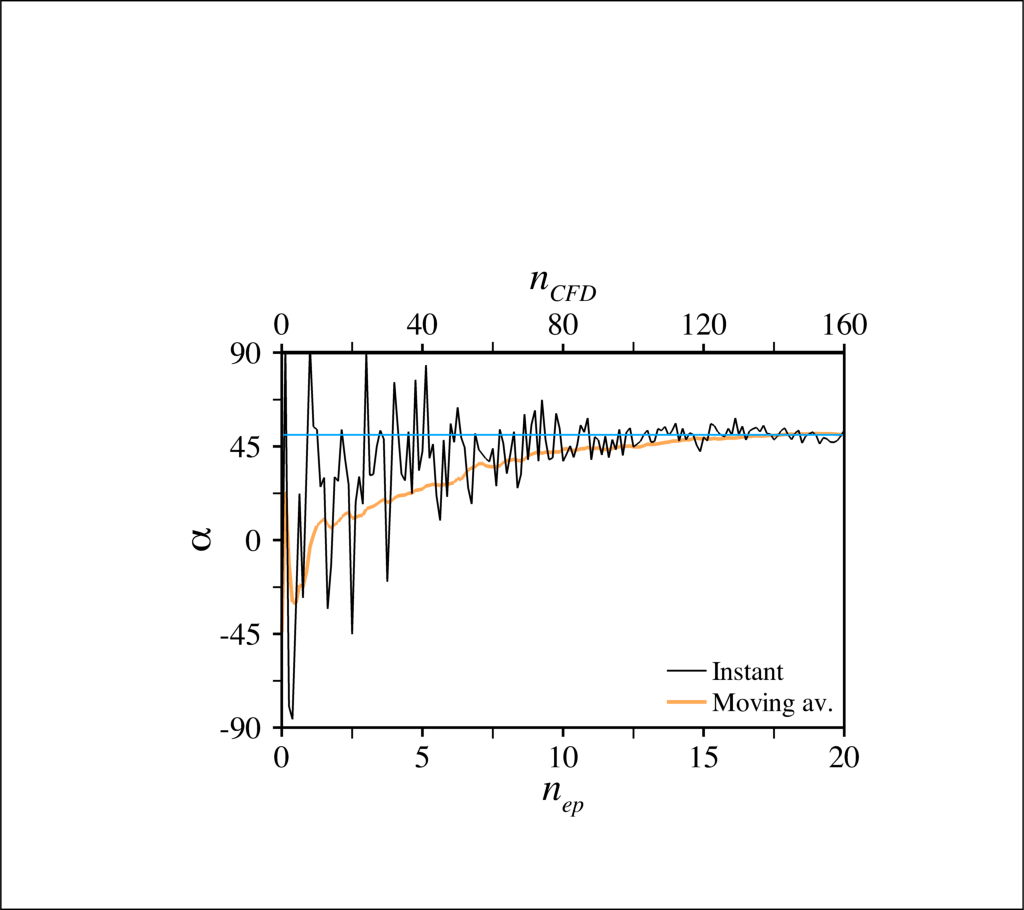}}
\put(4.4,1){\includegraphics[trim=465 275 380 445pt,clip,height=0.8cm]{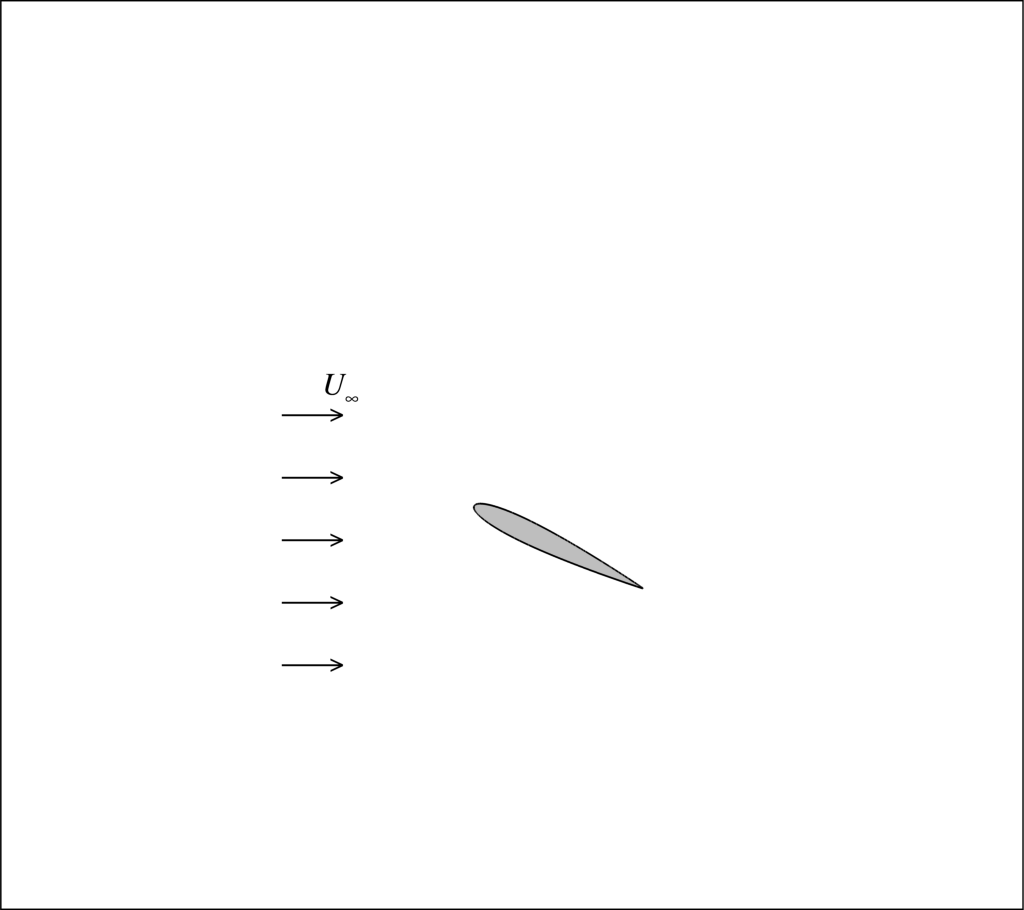}}
\put(11.95,1){\includegraphics[trim=465 275 380 445pt,clip,height=0.8cm]{fig4_thumb.png}}
\put(0.1,6){(a)}
\put(7.65,6){(b)}
\end{picture}
\caption{Flow past a NACA 0012 at $\rey=100$ - (a) Evolution per episode for the instant (black line) and moving average (over episodes, orange line) values of the mean lift (over time). The corresponding number of CFD simulations (obtained multiplying by the number of environments) is displayed on the secondary horizontal axis. (b) Same as (a) for the angle of attack. The blue lines mark the DNS optimal.}
\label{fig:naca_ppo}
\end{figure}

For each PPO-1 learning episode, the network outputs a single value $\xi$ in $[-1; 1]$ mapped into
\bal
\alpha=\xi\alpha_{\text{max}}\,,
\eal
for the angle of attack to vary in $[-\alpha_{\text{max}};\alpha_{\text{max}}]$ 
with $\alpha_{\text{max}}=90^\circ$. The reward $r = \overline{c_y}$ is then computed using the same simulation parameters, after which the network is updated for 32 epochs using 8 environments and 4 steps mini-batches. 20 episodes have been run for this case, which represents 160 simulations, each of which lasts $\sim25$mn using 12 cores,\footnote{\red{This is the time needed to compute periodic vortex shedding solutions. It takes less than 10mn to march the solution to steady state, but this barely affects the total CPU cost, as the time needed to complete an episode is that of completing its longest simulation (so only the cost of those episodes exclusively computing steady state solutions is reduced by a few minutes).}} hence $\sim65$h of total CPU cost (equivalently, $\sim 8$h of resolution time).
We show in figure~\ref{fig:naca_ppo}(a) the evolution of the reward collected over the course of the optimization. The moving average increases almost monotonically and reaches a plateau after about 15 episodes, and the optimal lift computed as the average over the 5 latest episodes is $\overline{c_y}^{\,\star}=0.93\pm 0.01$ (the variations are computed from the rms of the moving average over the same interval,
which is a simple yet robust criterion to assess qualitatively convergence a posteriori). The associated angle 
$\alpha^\star=50.2^\circ\pm 1.2^\circ$ varies by a larger factor, which is because lift is relatively insensitive to the exact incidence in the vicinity of the optimal. This is perfectly in line with the DNS, as illustrated by the red lines in figure~\ref{fig:naca}(b) showing the limits of the so-computed variance intervals.
\red{Nonetheless, PPO-1 turns to be rather inefficient at finding the optimal,
because it must span continuous ranges of angles while the one-dimensionality of the control space and the smoothness of the optimal allow DNS to test only a few discrete values (hence it can converge within $\sim 1$h using the same level of CFD parallelization).}

\subsection{Flow past an arrangement of two side-by-side circular cylinders} \label{sec:opt:sub:tandem} 

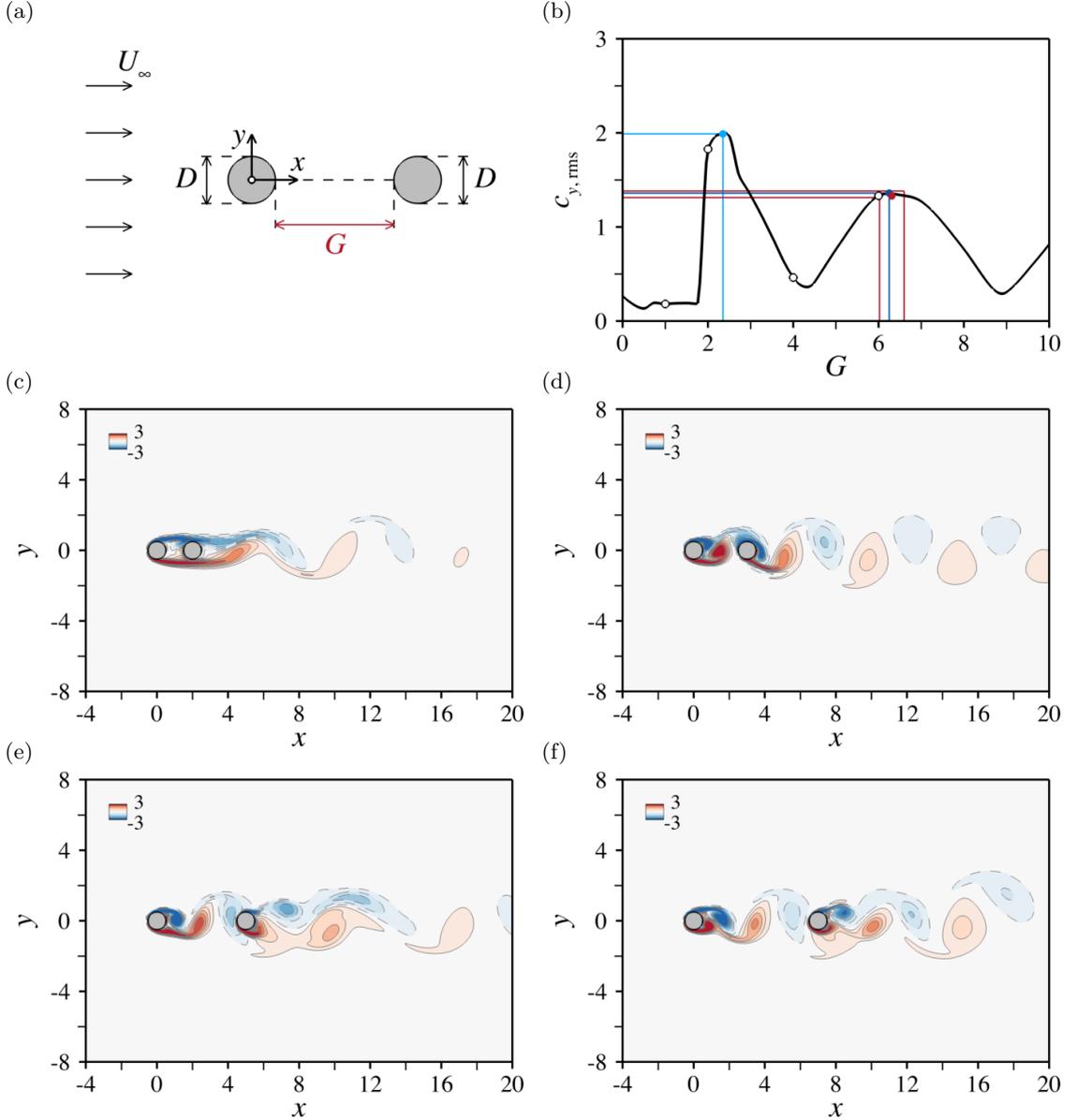
\begin{figure}[!t]
\setlength{\unitlength}{1cm}
\begin{picture}(20,15.75)
\put(0.1,10.5){\includegraphics[trim=175 92.5 155 340pt,clip,height=5.1cm]{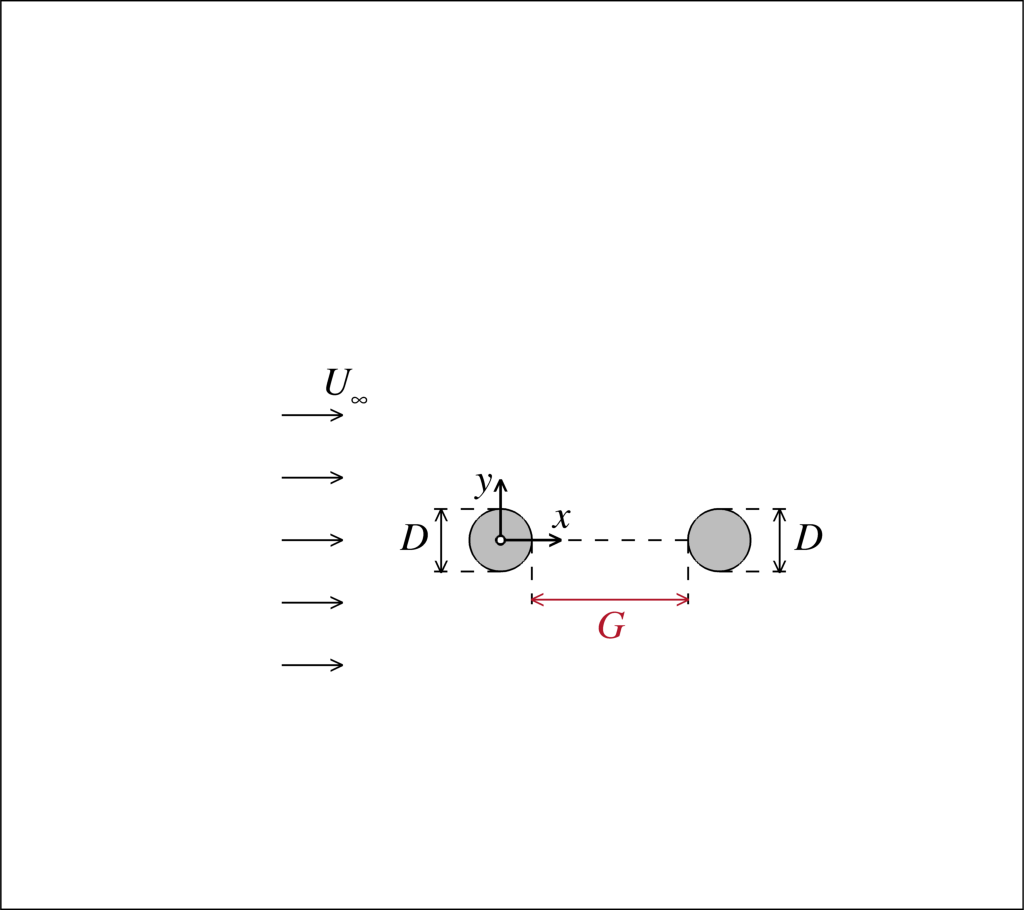}}
\put(7.65,10.5){\includegraphics[trim=175 92.5 155 340pt,clip,height=5.1cm]{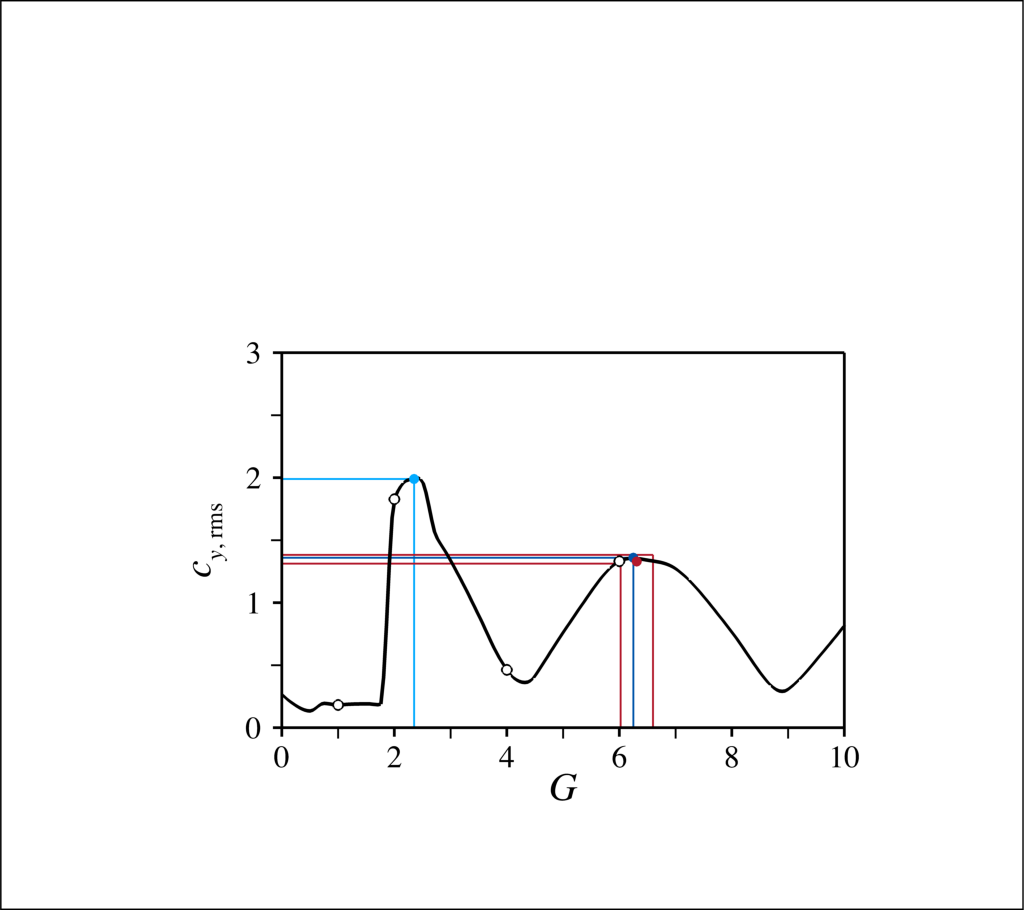}}
\put(0.1,5.25){\includegraphics[trim=175 92.5 155 340pt,clip,height=5.1cm]{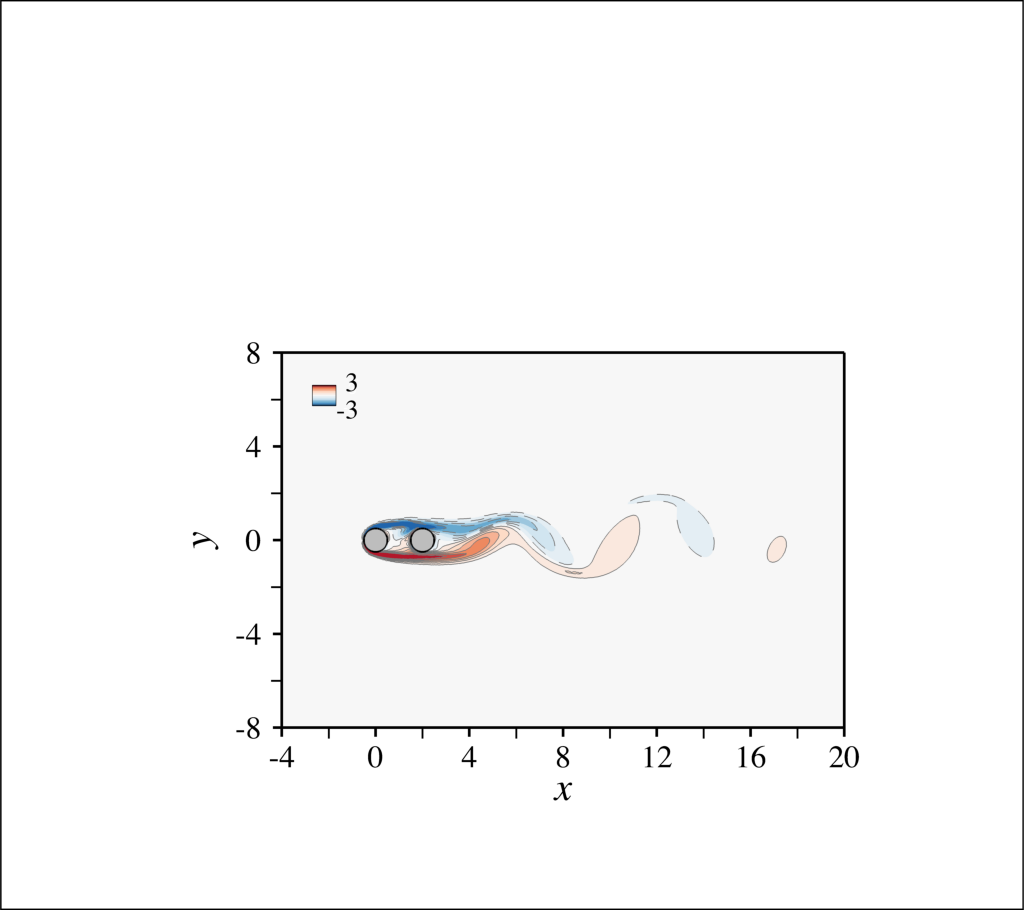}}
\put(7.65,5.25){\includegraphics[trim=175 92.5 155 340pt,clip,height=5.1cm]{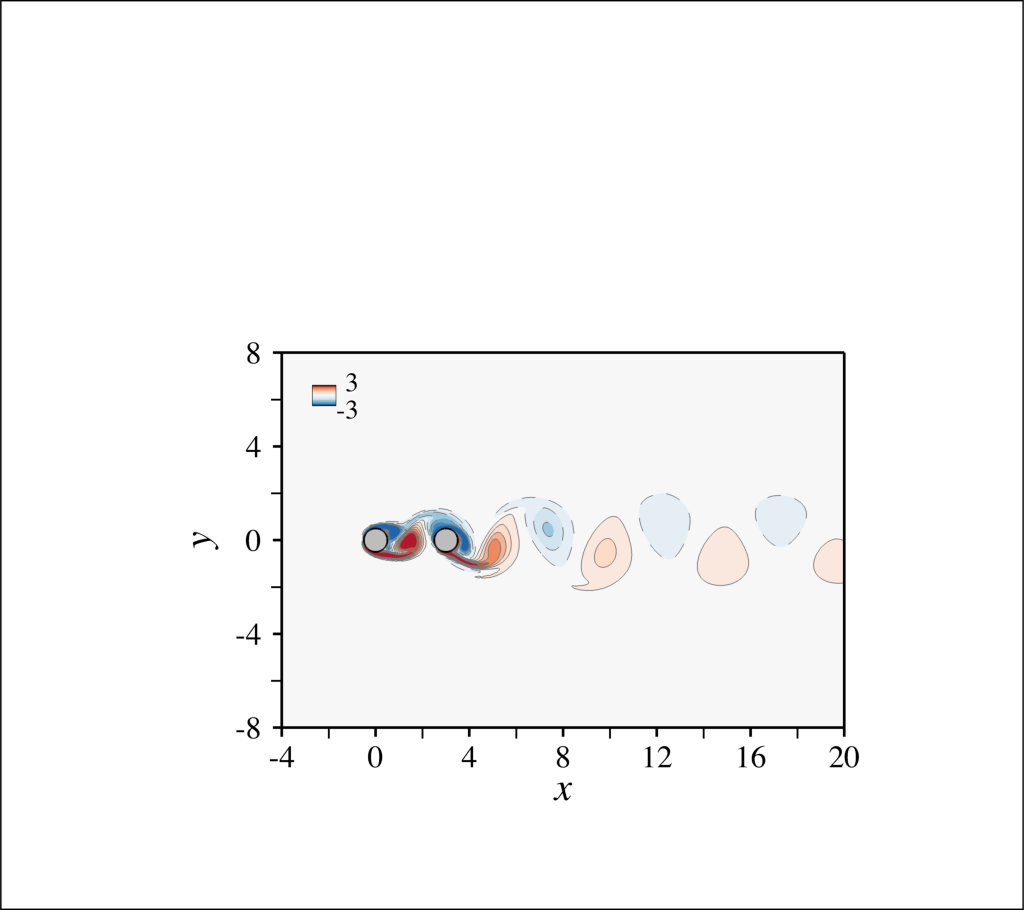}}
\put(0.1,0){\includegraphics[trim=175 92.5 155 340pt,clip,height=5.1cm]{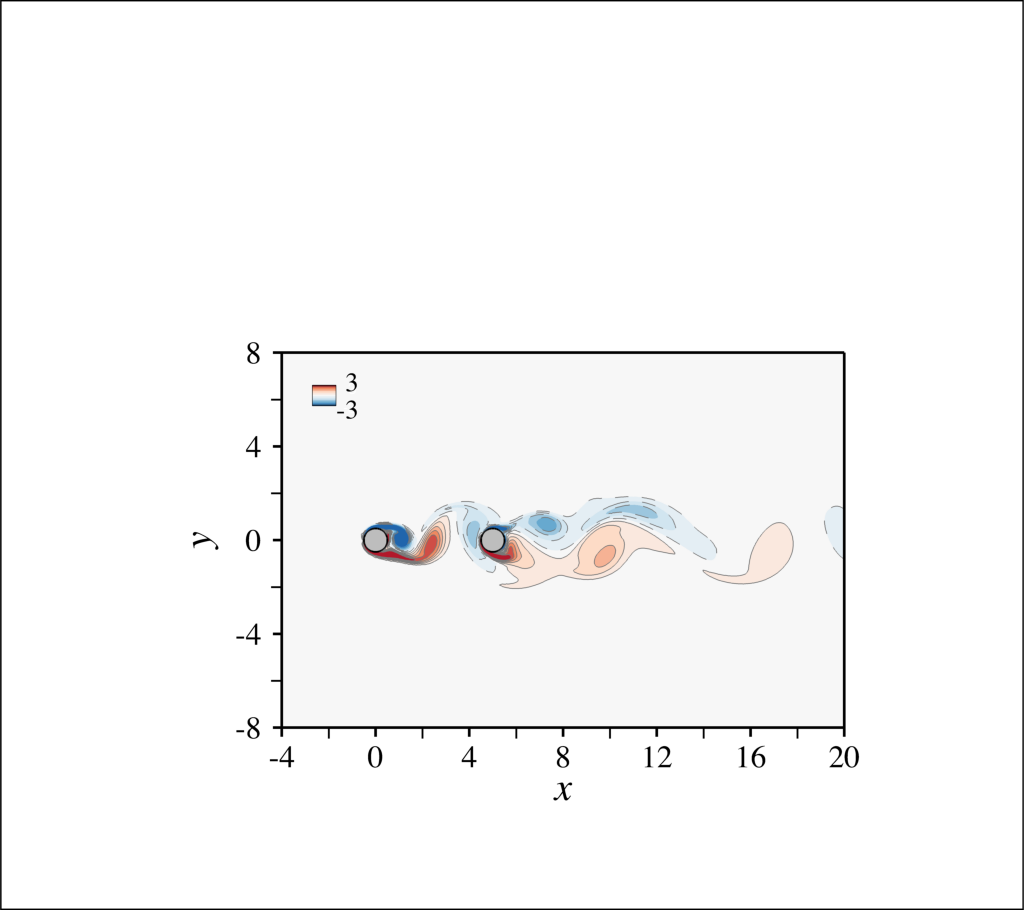}}
\put(7.65,0){\includegraphics[trim=175 92.5 155 340pt,clip,height=5.1cm]{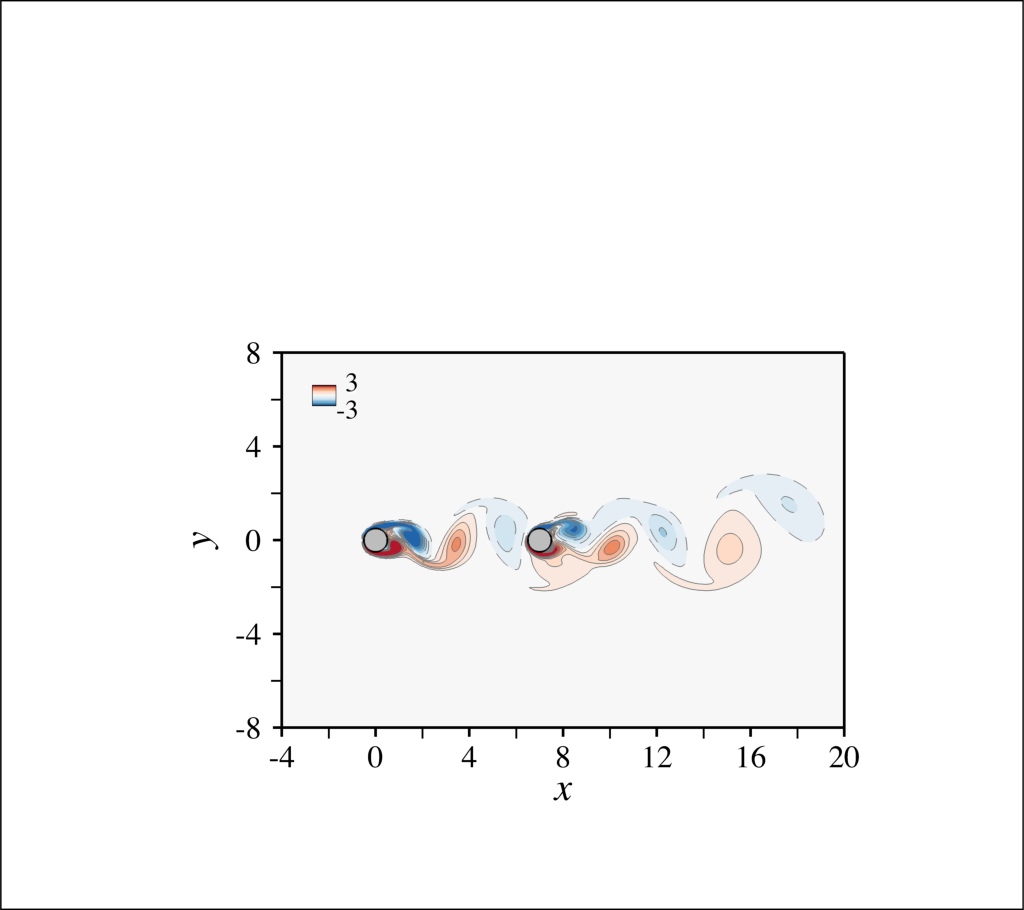}}
\put(0.1,15.75){(a)}
\put(7.65,15.75){(b)}
\put(0.1,10.5){(c)}
\put(7.65,10.5){(d)}
\put(0.1,5.25){(e)}
\put(7.65,5.25){(f)}
\end{picture}
\caption{Flow past the tandem arrangement of two circular cylinders - (a) Schematic diagram of the configuration. (b) Fluctuating (rms) lift against the gap spacing computed by DNS at $\rey=300$. The red symbol is the average over the 5 latest single-step PPO episodes, and the red lines 
delimit the corresponding variance intervals. (c-f) Instantaneous vorticity fields computed at $\rey=300$, for values marked by the circle symbols in (b), namely (c) $G=1$, (d) $G=2$, 
(d) $G=4$, and (e) $G=6$.} 
\label{fig:tandem}
\end{figure}


We examine now the side-by-side tandem arrangement of two identical circular cylinders in a uniform stream, whose configuration is sketched in figure~\ref{fig:tandem}(a). The origin of the coordinate system is at the center of the main cylinder, where we refer to the upstream and downstream cylinders as ``main'' and ``surrounding'', respectively. 
A laminar, time-dependent case at $\rey=U_\infty D/\nu=300$
is modeled after the Navier--Stokes equations, where $D$ is the diameter of either cylinder. The objective is to maximize the rms lift $c_{y,\text{rms}}$ of the two-cylinder system (for instance, to increase the amount of energy available for harnessing from fluid-structure interactions) for which the sole control parameter is the gap spacing $G$, i.e., the side-to-side distance between the two cylinders.
On paper, this is another problem simple enough to allow direct comparisons between PPO-1 and DNS. \red{In practice, the results are not so unequivocal, as evidenced in figure~\ref{fig:tandem}(b) showing} 
reference data obtained from $30$ DNS runs computing the rms lift to an accuracy of $5\%$ with the simulation parameters documented in table~\ref{tab:opt}.
A steep global maximum lies at $G^\star=2.35$, associated with $c_{y,\text{rms}}^{\,\star}=1.99$, but there is a smoother local maximum at $G^{\star\star}=6.25$, associated with $c_{y,\text{rms}}^{\,\star\star}=1.36$,
which reflects the high sensitivity of the pattern of flow unsteadiness to the center distance. 
Namely without going into too much detail (as this has been extensively discussed in the literature~\cite{Mittal1997,Meneghini2001,Sharman2005,Lee2009}), the instantaneous vorticity field computed for $G = 1$ in figure~\ref{fig:tandem}(c) shows that the gap flow between the two cylinders is initially steady, while the shear layers separating from the main cylinder engulf those of the surrounding cylinder and trigger vortex shedding in the far wake. For $G=2$ (close to the global maximum), the gap flow is unsteady, but the gap vortices are not fully developed by the time they impinge on the surrounding cylinder, hence a single vortex street in the far wake; see figure~\ref{fig:tandem}(d). For $G=4$, one pair of gap vortices fully develops, then impinges on the surrounding cylinder, which triggers a complex interaction in the near wake before a vortex street eventually forms further downstream; see figure~\ref{fig:tandem}(e). Finally for $G=6$ (close to the local maximum) the wake of the surrounding cylinder is unsteady again, and both cylinders shed synchronized vortices close to anti-phase; see figure~\ref{fig:tandem}(f).

\begin{figure}[!t]
\setlength{\unitlength}{1cm}
\begin{picture}(20,6)
\put(0.1,0){\includegraphics[trim=175 92.5 155 270pt,clip,height=5.85cm]{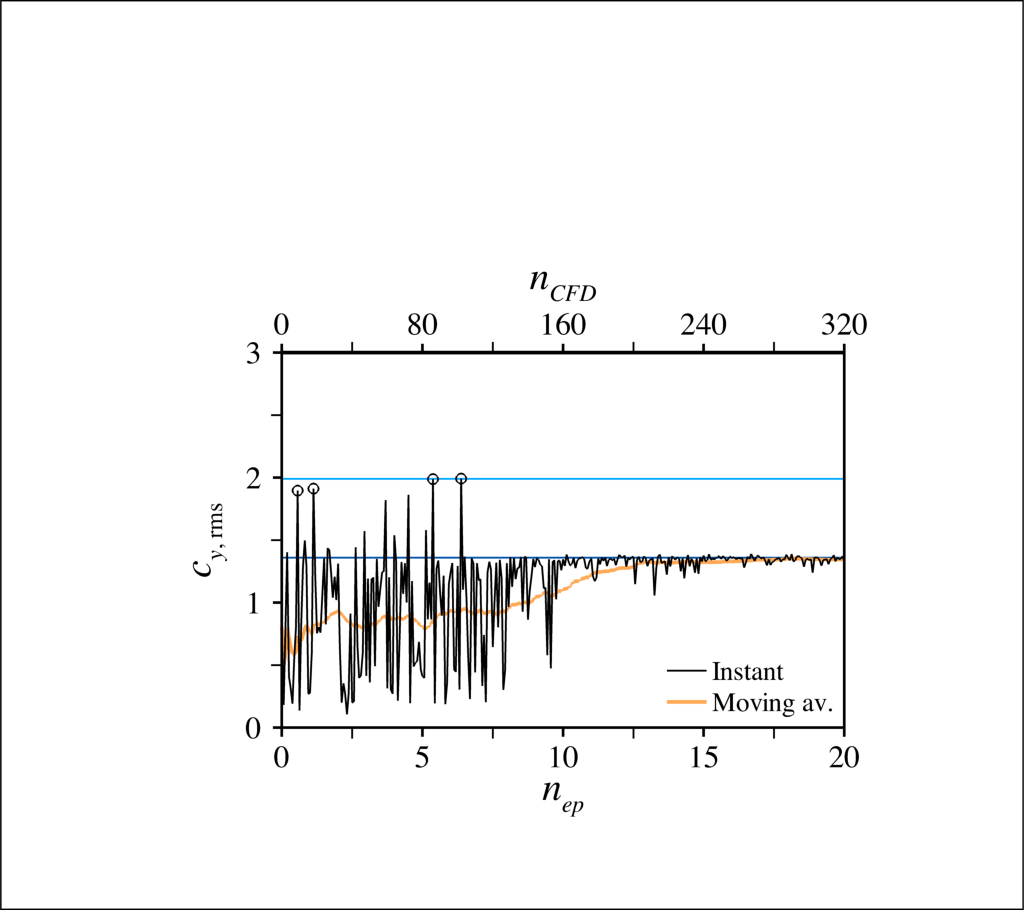}}
\put(7.65,0){\includegraphics[trim=175 92.5 155 270pt,clip,height=5.85cm]{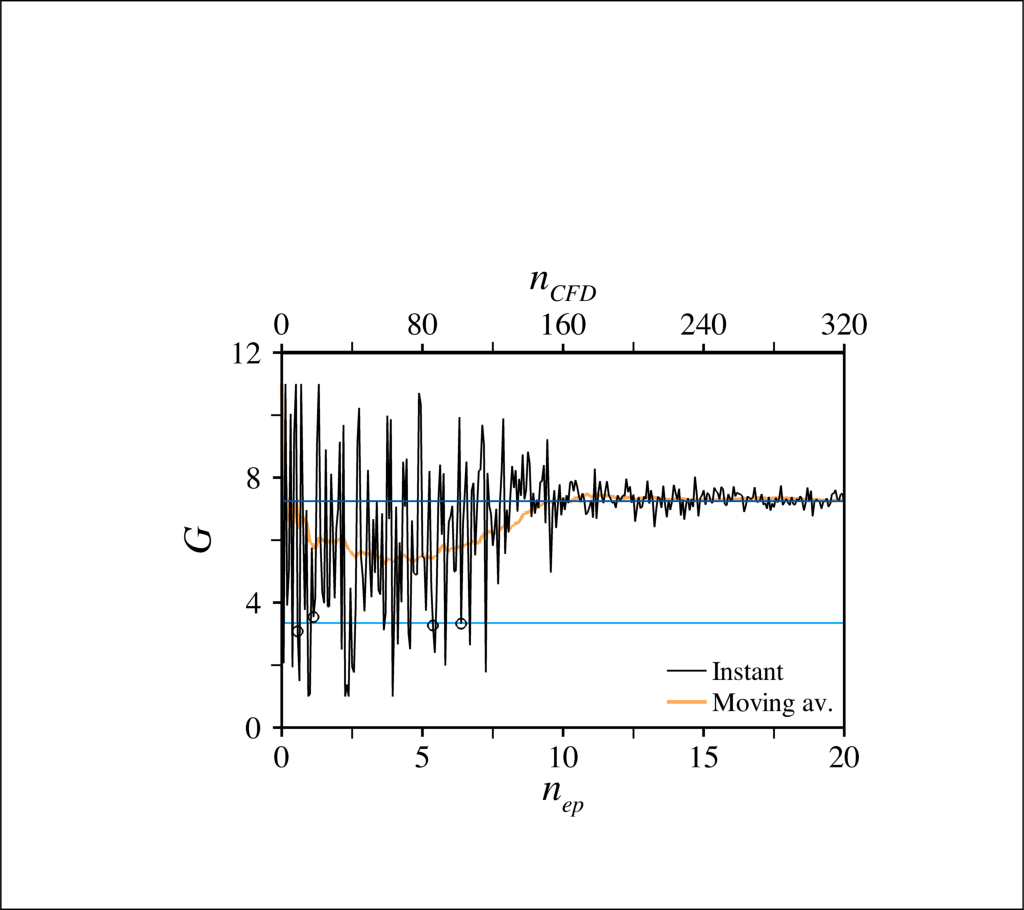}}
\put(4.5,1){\includegraphics[trim=465 275 380 445pt,clip,height=0.8cm]{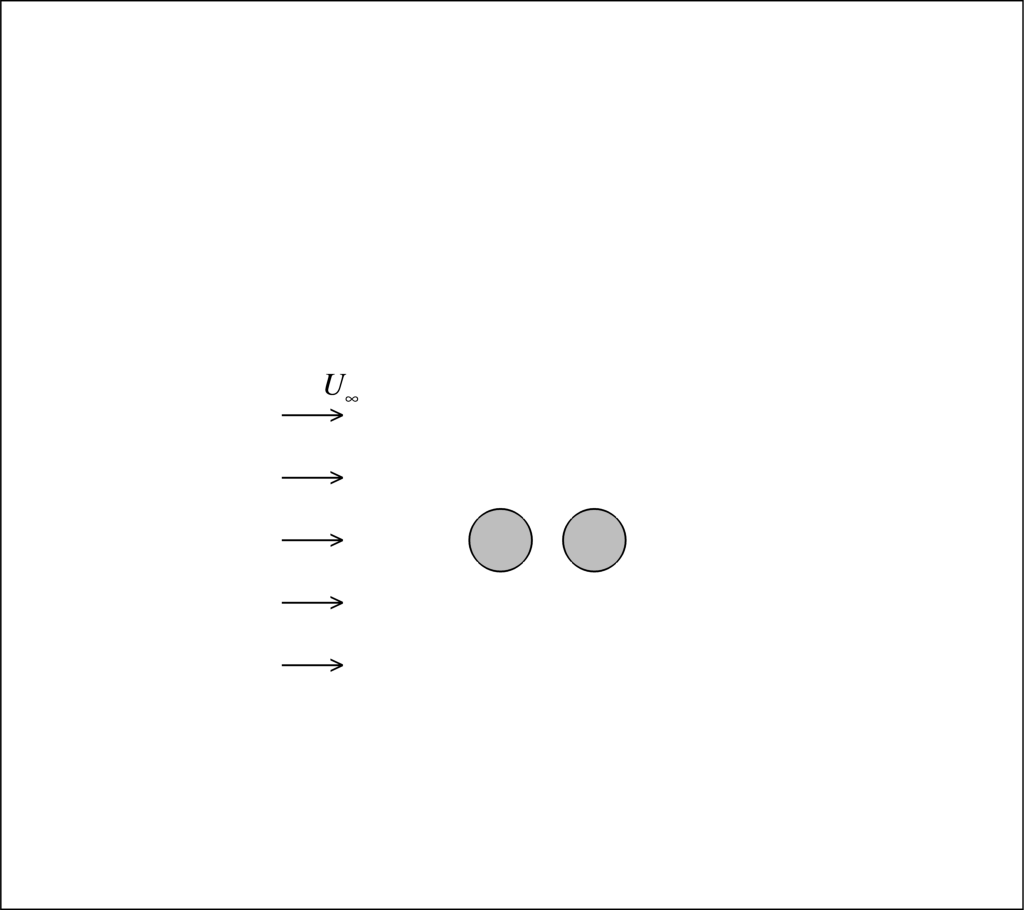}}
\put(12.05,1){\includegraphics[trim=465 275 380 445pt,clip,height=0.8cm]{fig6_thumb.png}}
\put(0.1,6){(a)}
\put(7.65,6){(b)}
\end{picture}
\caption{Flow past the tandem arrangement of two circular cylinders at $\rey=300$ - (a) Evolution per episode for the instant (black line) and moving average (over episodes, orange line) values of the rms lift. (b) Same as (a) for the gap spacing. The light (resp. dark) blue lines mark the DNS global (resp. local) maximum. The circles are high reward parameters close to the DNS global maximum.}
\label{fig:tandem_ppo}
\end{figure}

For each PPO-1 learning episode, the network outputs a single value $\xi$ in $[-1; 1]$ mapped into 
\bal
G=\frac{1+\xi}{2}G_{\text{max}}\,,
\eal
for the gap to vary in $[0;G_{\text{max}}]$ with $G_{\text{max}}=10$. This enables contact between the two cylinders and keeps the computational cost affordable, as pushing the surrounding cylinder further downstream would require extending the computational domain and increasing the numbers of grid points accordingly (all the more so as we do not anticipate such large distances to be relevant from the standpoint of optimization because the interaction between both cylinders will weaken increasingly at some point, although it can take up to several tens of diameters to do so).
The reward $r =c_{y,\text{rms}}$ is then computed using the same simulation parameters, after which the network is updated for 32 epochs using 16 environments and 4 steps mini-batches. 
Another 20 episodes have been run for this case. This represents 320 simulations, each of which lasts 
$\sim60$mn on 12 cores (much longer than in the NACA case due to the increased simulation time), 
hence $\sim 320$h of total CPU cost (equivalently, $\sim20$h of resolution time), \red{still much more than by DNS because DRL keeps spanning continuous ranges of distances while DNS can settle for only a few discrete values despite the sharpness of the global maximum (hence it can converge within $\sim 3$h using the same level of CFD parallelization).}
Figure~\ref{fig:tandem_ppo}(a) shows a plateau in the moving average reward after about 15 episodes. The optimal lift computed as the average over the 5 latest episodes is $c_{y,\text{rms}}^{\,\star}=1.34\pm 0.02$, associated with $G^\star=6.31\pm 0.04$, meaning that the agent misses the global maximum, but converges to a value close to the local maximum; see the red lines in figure~\ref{fig:tandem}(b) 
indicating the limits of the computed variance intervals.

\begin{figure}[!t]
\setlength{\unitlength}{1cm}
\begin{picture}(20,5.5)
\put(0.1,0){\includegraphics[trim=175 92.5 155 340pt,clip,height=5.1cm]{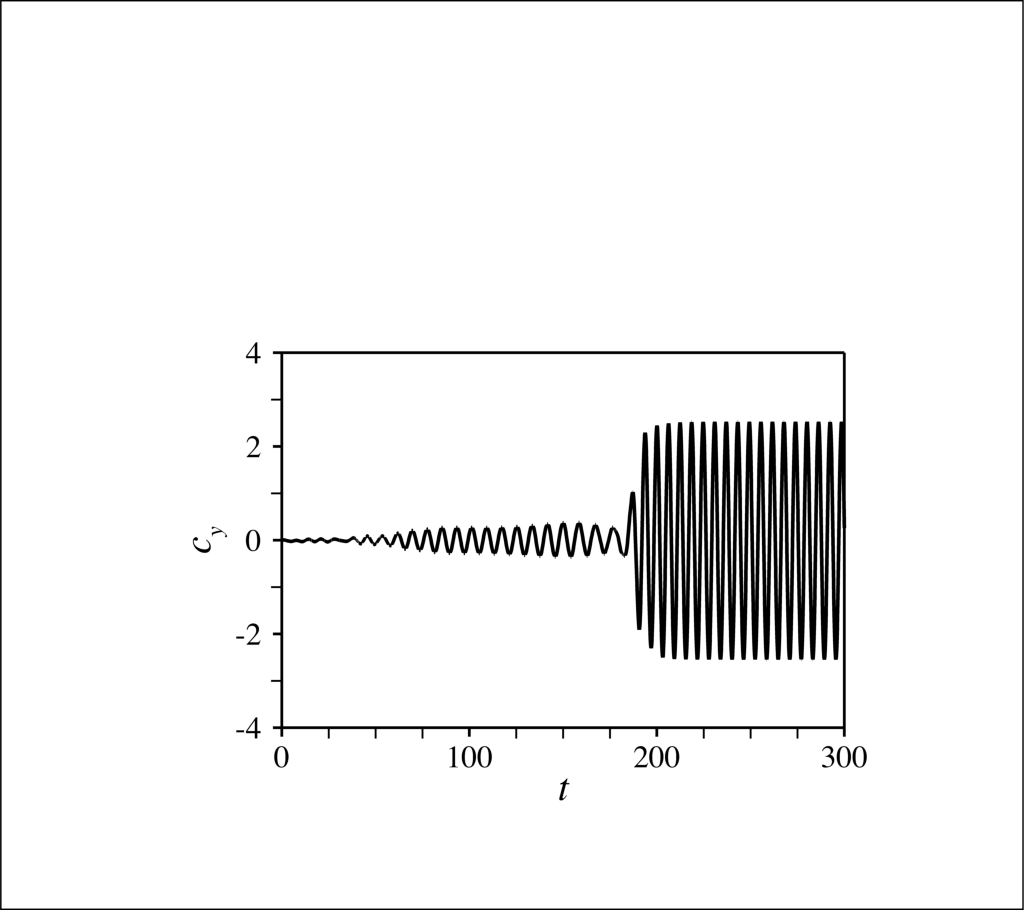}}
\put(7.65,0){\includegraphics[trim=175 92.5 155 340pt,clip,height=5.1cm]{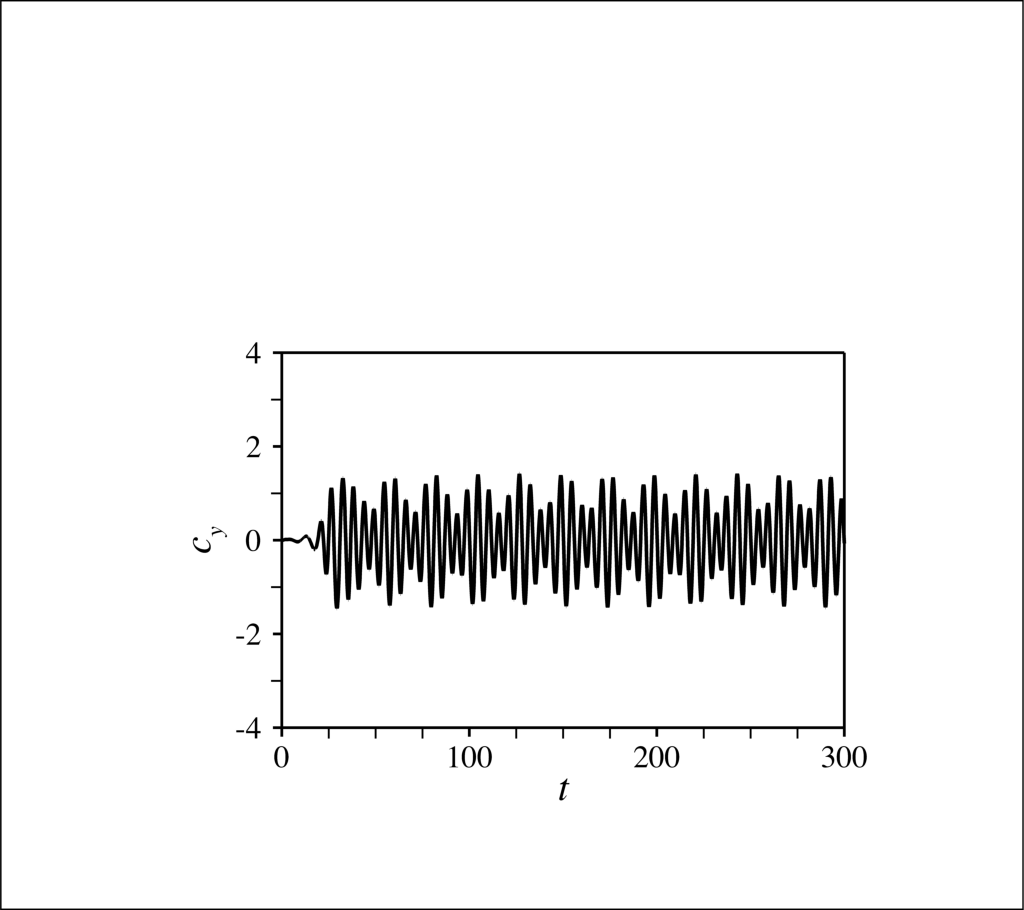}}
\put(0.1,5.5){(a)}
\put(7.65,5.5){(b)}
\end{picture}
\caption{Flow past the tandem arrangement of two circular cylinders at $\rey=300$ - Time history of lift computed by DNS for gap spacings (a) $G=2$ (close to the DNS global maximum) and (b) $G=8$.} 
\label{fig:tandemvst_re300}
\end{figure}


\begin{table}[!t]
\begin{center}
\begin{tabular}{ccccccccc}
\toprule
\multicolumn{1}{r}{\multirow{1}{*}{\makecell[r]{\raisebox{-\totalheight}{\includegraphics[trim=465 225 350 495pt,clip,height=1cm]{fig4_thumb.png}}}}}
& \multicolumn{1}{p{0.9cm}}{\makecell[r]{$\overline{c_y}$}} & \multicolumn{1}{p{0.9cm}}{\makecell[r]{$\alpha$}} & &
\multicolumn{1}{r}{\multirow{3}{*}{\makecell[r]{\raisebox{-\totalheight}{\includegraphics[trim=425 225 350 495pt,clip,height=1cm]{fig6_thumb.png}}}}}
& \multicolumn{1}{p{1cm}}{\makecell[r]{$c_{y,\text{rms}}$}} & \multicolumn{1}{p{1cm}}{\makecell[r]{$G$}} \\
\cmidrule(lr){2-3}\cmidrule(lr){6-7}
\multicolumn{2}{r}{$0.93$} & \multicolumn{1}{r}{$50.2^\circ$} & \multicolumn{1}{p{1.2cm}}{PPO-1} & & \multicolumn{1}{r}{$1.34$} & \multicolumn{1}{r}{$6.31$}  & \multicolumn{1}{p{1.2cm}}{PPO-1}& \multicolumn{1}{l}{\multirow{3}{*}{\qquad Optimal}}\\
\multicolumn{2}{r}{$0.94$} & \multicolumn{1}{r}{$50.6^\circ$} & \multicolumn{1}{p{1.2cm}}{DNS}  & & \multicolumn{1}{r}{$1.99$} & \multicolumn{1}{r}{$2.35$} & \multicolumn{1}{p{1.2cm}}{\multirow{2}{*}{DNS}}\\
\multicolumn{6}{r}{$\gray{1.36}$} & \multicolumn{1}{r}{$\gray{6.25}$} \\
\cmidrule(lr){1-9}
\multicolumn{9}{r}{CFD} \\
\cmidrule(lr){1-9}
\multicolumn{3}{r}{100} & & \multicolumn{3}{r}{300} &&  \multicolumn{1}{l}{\qquad Reynolds number} \\
\multicolumn{3}{r}{0.125} & & \multicolumn{3}{r}{\guillemotright} & & \multicolumn{1}{l}{\qquad Time-step} \\
\multicolumn{3}{r}{$[50;150]$} & & \multicolumn{3}{r}{$[200;300]$} & & \multicolumn{1}{l}{\qquad Averaging time span} \\
\multicolumn{3}{r}{$[-15;40]\times[-15;15]$} & & \multicolumn{3}{r}{\guillemotright} & & \multicolumn{1}{l}{\qquad Mesh dimensions} \\
\multicolumn{3}{r}{$115000$} & & \multicolumn{3}{r}{125000} & & \multicolumn{1}{l}{\qquad Nb. mesh elements} \\
\multicolumn{3}{r}{$0.001$} & & \multicolumn{3}{r}{\guillemotright} & & \multicolumn{1}{l}{\qquad Interface $\perp$ mesh size} \\
\multicolumn{3}{r}{12} & & \multicolumn{3}{r}{\guillemotright} & & \multicolumn{1}{l}{\qquad Nb. Cores} \\
\cmidrule(lr){1-9}
\multicolumn{9}{r}{PPO-1} \\
\cmidrule(lr){1-9}
\multicolumn{3}{r}{20} & & \multicolumn{3}{r}{\guillemotright} & & \multicolumn{1}{l}{\qquad Nb. DRL episodes} \\
\multicolumn{3}{r}{8} & & \multicolumn{3}{r}{16} & & \multicolumn{1}{l}{\qquad Nb. Environments} \\
\multicolumn{3}{r}{32} & & \multicolumn{3}{r}{\guillemotright} & & \multicolumn{1}{l}{\qquad Nb. Epochs} \\
\multicolumn{3}{r}{4} & & \multicolumn{3}{r}{\guillemotright} & & \multicolumn{1}{l}{\qquad Size of mini-batches} \\
\multicolumn{3}{r}{60h} & & \multicolumn{3}{r}{320h} &&  \multicolumn{1}{l}{\qquad CPU time} \\
\multicolumn{3}{r}{7.5h} & & \multicolumn{3}{r}{20h} &&  \multicolumn{1}{l}{\qquad Resolution time} \\
\bottomrule
\end{tabular}
\caption{\label{tab:opt} Simulation parameters and convergence data for the flow past a NACA 0012 at $\rey=100$ and the flow past the tandem arrangement of two circular cylinders at $\rey=300$. 
\red{NACA 0012: the interface mesh size yields $\sim 20$ grid points in the boundary-layer at mid-chord, under zero incidence, and the averaging time-span represents $\sim 15-20$ shedding cycles, depending on the incidence.
Tandem arrangement of two circular cylinders: the interface mesh size yields $\sim 20$ grid points in the boundary-layer of the main cylinder, just prior to separation, and the averaging time-span represents $\sim 20$ shedding cycles.}}
\end{center}
\end{table}

This half-failure can be explained by the steepness of the reward gradients with respect to the control variable in the vicinity of the global maximum. This is due to the existence of a secondary instability mechanism at play in a narrow range of center distances, as illustrated in figure~\ref{fig:tandemvst_re300}(a) 
showing that for $G\sim2$, the flow settles to a first time-periodic solution, 
then bifurcates to a second time-periodic solution associated with increased lift oscillations (hence the large values of $t_i$ used for this case).
Actually, DRL does identify high reward positions close to $G = 2$ (circle symbols in figure~\ref{fig:tandem_ppo}), whose value $c_{y,\text{rms}}\sim 2$ is consistent with
the global maximum, but there are very few times where the global maximum is met during the exploration phase (compared to its local counterpart, again because of the topology of the reward function). Because PPO voluntarily dismisses large policy updates to avoid performance collapse, the clipped policy updates only lead to limited exploration and trap the optimization process into a local maximum. 
\red{Low to moderate Reynolds numbers are likely required for such instability cascade scenario to occur, so such results do not cast doubt on the applicability of single-step PPO to practically meaningful high Reynolds flows. They do stress, however, that the method can benefit from carefully tuning the trade-off between exploration and exploitation, which will be addressed in future work.}

\section{Application to open-loop flow control} \label{sec:olc}

\subsection{Optimal cylinder drag reduction using a smaller control cylinder} \label{sec:olc:sub:control}

\begin{figure}[!t]
\setlength{\unitlength}{1cm}
\begin{picture}(20,10.5)
\put(0.1,5.25){\includegraphics[trim=175 92.5 155 340pt,clip,height=5.1cm]{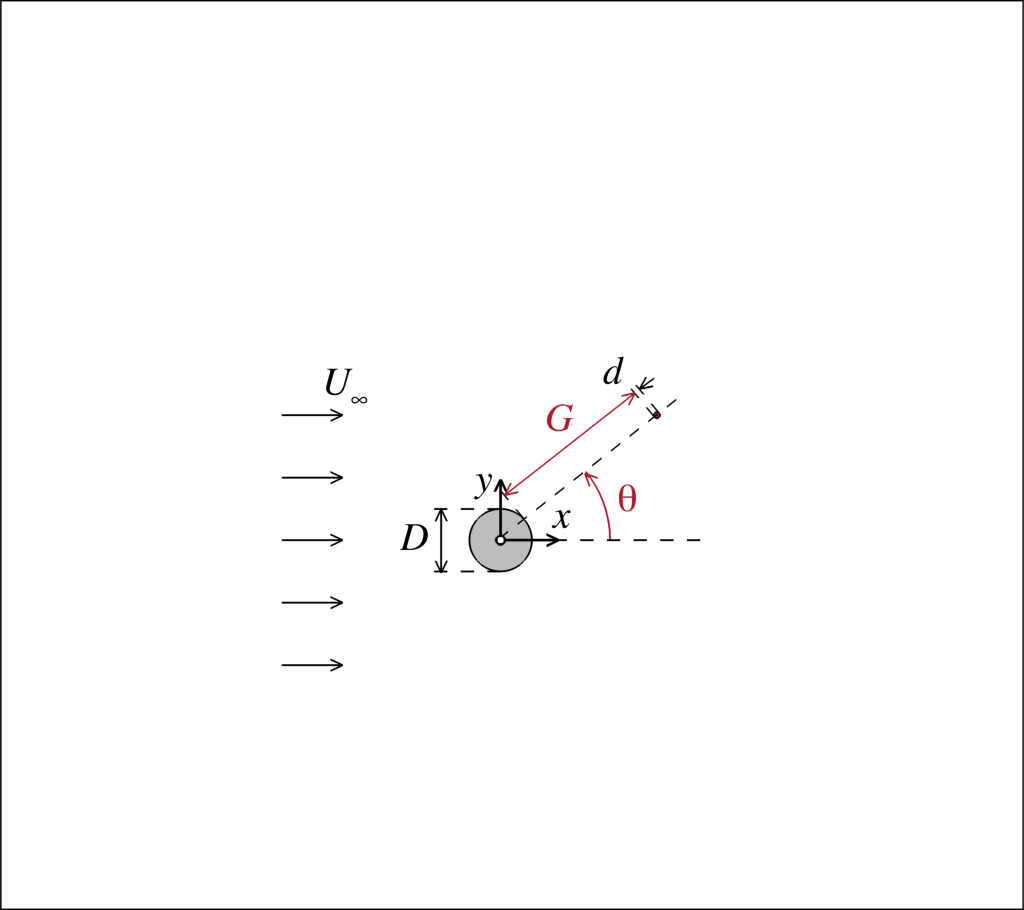}}
\put(7.65,5.25){\includegraphics[trim=175 92.5 155 340pt,clip,height=5.1cm]{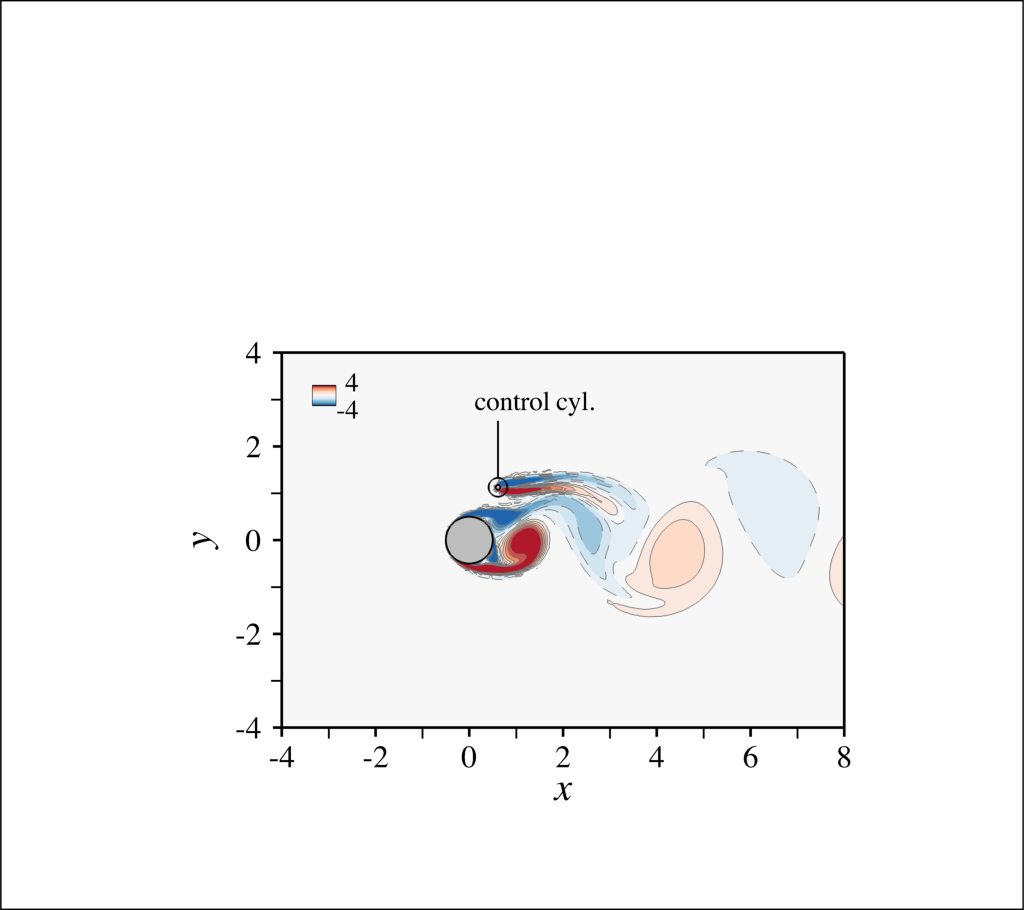}}
\put(0.1,0){\includegraphics[trim=175 92.5 155 340pt,clip,height=5.1cm]{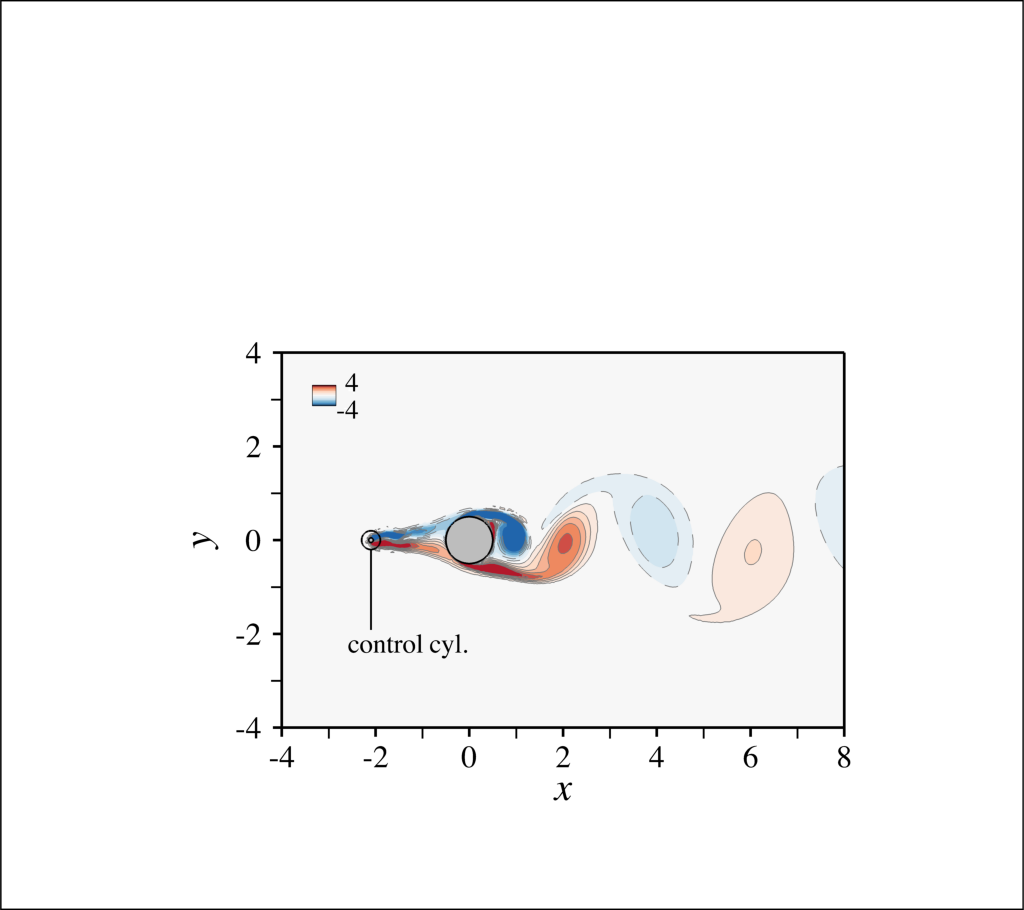}}
\put(7.65,0){\includegraphics[trim=175 92.5 155 340pt,clip,height=5.1cm]{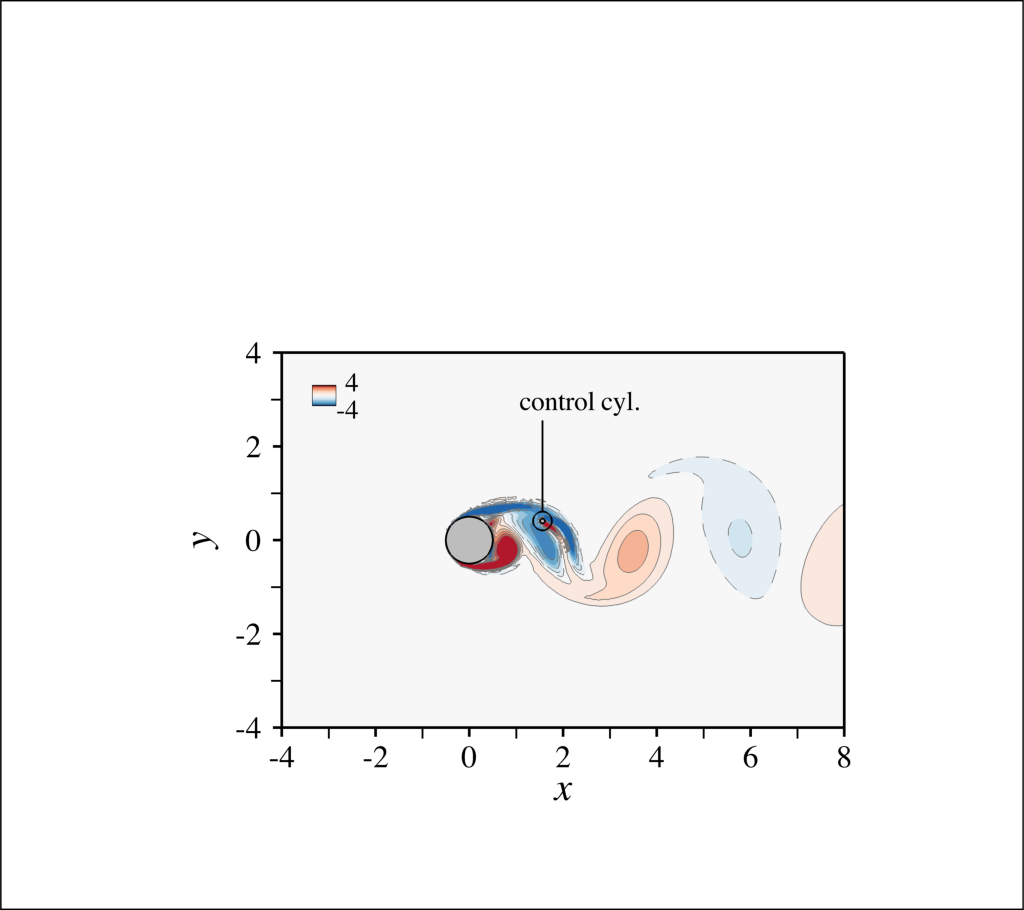}}
\put(0.1,10.5){(a)}
\put(7.65,10.5){(b)}
\put(0.1,5.25){(c)}
\put(7.65,5.25){(d)}
\end{picture}
\caption{Open-loop control of the circular cylinder flow by a small control cylinder of diameter $d=0.1$ - (a) Schematic diagram of the configuration. (b-d) Iso-contours of the vorticity field computed at $\rey=3900$ for representative positions $(x_c,y_c)$ of the control cylinder, namely (b) $(0.61,1.13)$, (c) $(-2.10,0.00)$ and (d) $(1.56,0.41)$.}
\label{fig:control}
\end{figure}
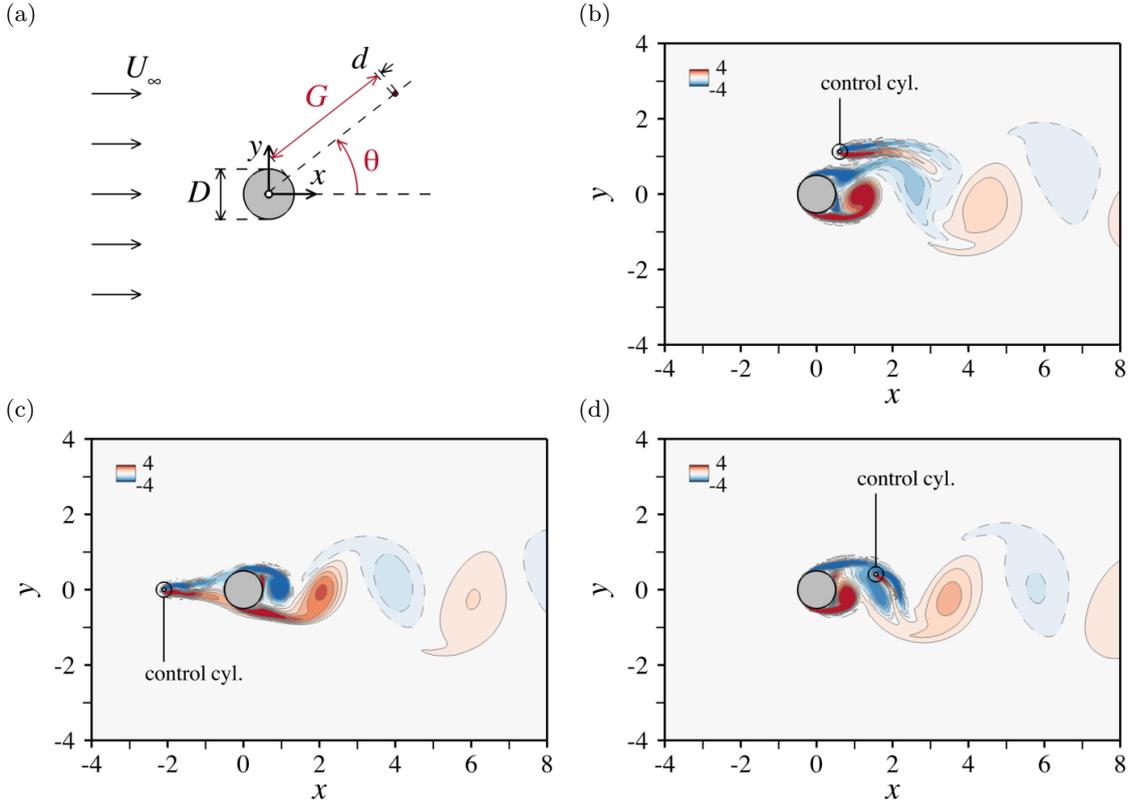

\red{The relevance of single-step PPO is now showcased by tackling various open-loop control problems.}
The first one is that of a cylinder in a uniform stream, controlled open-loop by a much smaller circular cylinder. Figure~\ref{fig:control}(a) presents a sketch of the configuration pertaining to a circular geometry of the main cylinder, where we refer to the large and small cylinders as ``main'' and ``control'', respectively, but section~\ref{sec:olc:sub:control22000} also considers a square geometry. The origin of the coordinate system is at the center of the main cylinder. 
The objective is to minimize the mean drag $\overline{c_x}$ of the two-cylinder system, \red{which requires reducing the drag of the main cylinder sufficiently to compensate for the fact that the control cylinder itself} is a source of drag. Several laminar and turbulent Reynolds numbers $\rey=U_\infty D/\nu$ are considered, where $D$ is the diameter of the main cylinder. 
The diameter of the control cylinder is set to $d=0.1$, therefore the sole control parameter is the 2-D position of the control cylinder center, measured by the gap distance $G$ between the two cylinders and the azimuthal position $\theta$ with respect to the rear stagnation point.
This may not seem overly complicated on paper, but the parameter space is actually large enough to dismiss mapping the best positions for placement of the control cylinder by DNS, as tens of thousands of runs are required to cover merely a few diameters around the main cylinder.
In the following, single-step PPO is thus compared to theoretical predictions obtained by the adjoint method. The latter has proven fruitful to gain insight into the most efficient region from the linear sensitivity of the uncontrolled flow (i.e., the flow past the main cylinder), without ever calculating the controlled states, using 
instead a simple model of the force exerted by the control cylinder on the flow.
We shall not go into the technicalities of how to derive the related adjoint equations,
as the line of thought here is to take the output sensitivity as a given to assess relevance of PPO-1.
Suffice it to say here that we rely on various levels of adjoint modeling whose key assumptions are 
reviewed in appendix~\ref{sec:appadj}. The reader interested in more details is directed to the original literature on this topic~\cite{meli14,mao15,meli18ejmb}, where 
 in-depth technical and mathematical information, together with 
extensive discussions regarding the validity of the approximations are available. From the numerical standpoint, all calculations are performed with the mixed finite elements adjoint solver presented and validated in~\cite{meli14}.

On the CFD side, one of the challenges lies in the fact that the control cylinder acts as a small local disturbance redistributing the vorticity in the separated shear layers; see figures~\ref{fig:control}(b-d) showing instantaneous vorticity fields computed for representative positions of the control cylinder. \red{Accurate numerical methods are thus} mandatory to capture the small drag variations induced by the control. Several values of the Reynolds number are investigated : a laminar, steady case at $\rey=40$, for which the flow remains steady-state regardless of the position of the control cylinder, a laminar, time-dependent case at $\rey=100$, for which vortex shedding consistently develops from the main cylinder but the flow past the control cylinder remains steady, and two turbulent cases at $\rey=3900$ and at $\rey=22000$ (hence modeled after the uRANS equations with negative Spalart--Allmaras as turbulence model), for which vortex shedding develops from both cylinders. This is because the Reynolds number in the wake of the control cylinder must be scaled by the ratio of the cylinder diameters, which yields values below (resp. above) the instability threshold at $\rey=100$ (resp. $\rey=3900$ and $\rey=22000$).

For each PPO-1 episode, the network outputs two values $\xi_{{1, 2}}$ in $[-1; 1]^2$ mapped into 
\bal
G=\frac{1+\xi_1}{2}G_{\text{max}}\,, \qquad\qquad
\theta=\frac{1+\xi_2}{2}\theta_{\text{max}}\,,\label{eq:mapcontrol}
\eal
\red{for the gap to vary in $[0;G_{\text{max}}]$ with $G_{\text{max}}=3$, and the azimuthal position to vary in $[0;\theta_{\text{max}}]$ with $\theta_{\text{max}}=180^\circ$. This enables contact between the two cylinders, and allows taking advantage of the problem symmetry, as it amounts to moving the control cylinder in the upper half of a torus bounded by the surface of the main cylinder and the user-defined exterior radius $G_{\text{max}}$}. In the following, the center position is conveniently presented in terms of the Cartesian coordinates $x_c=\rho\cos\theta$ and $y_c=\rho\sin\theta$, where we note $\rho=G+(1+d)/2$.
Since the aim is to minimize drag, the reward $r=-\overline{D}$ is then computed using the simulation parameters documented in table~\ref{tab:control}, \red{after which the network is updated for
32 epochs using 8 environments and 2 steps mini-batches (note the zero averaging span in table~\ref{tab:control} for $\rey=40$, as this is a steady case for which the steady asymptotic value of total drag can be evaluated at the final time $t_f$, provided it is large enough for the solution to relax to steady-state).}

\subsubsection{Laminar steady regime and circular geometry at \text{Re=40}} \label{sec:olc:sub:control40} 

\begin{figure}[!t]
\setlength{\unitlength}{1cm}
\begin{picture}(20,12.25)
\put(0.1,6.25){\includegraphics[trim=175 92.5 155 270pt,clip,height=5.85cm]{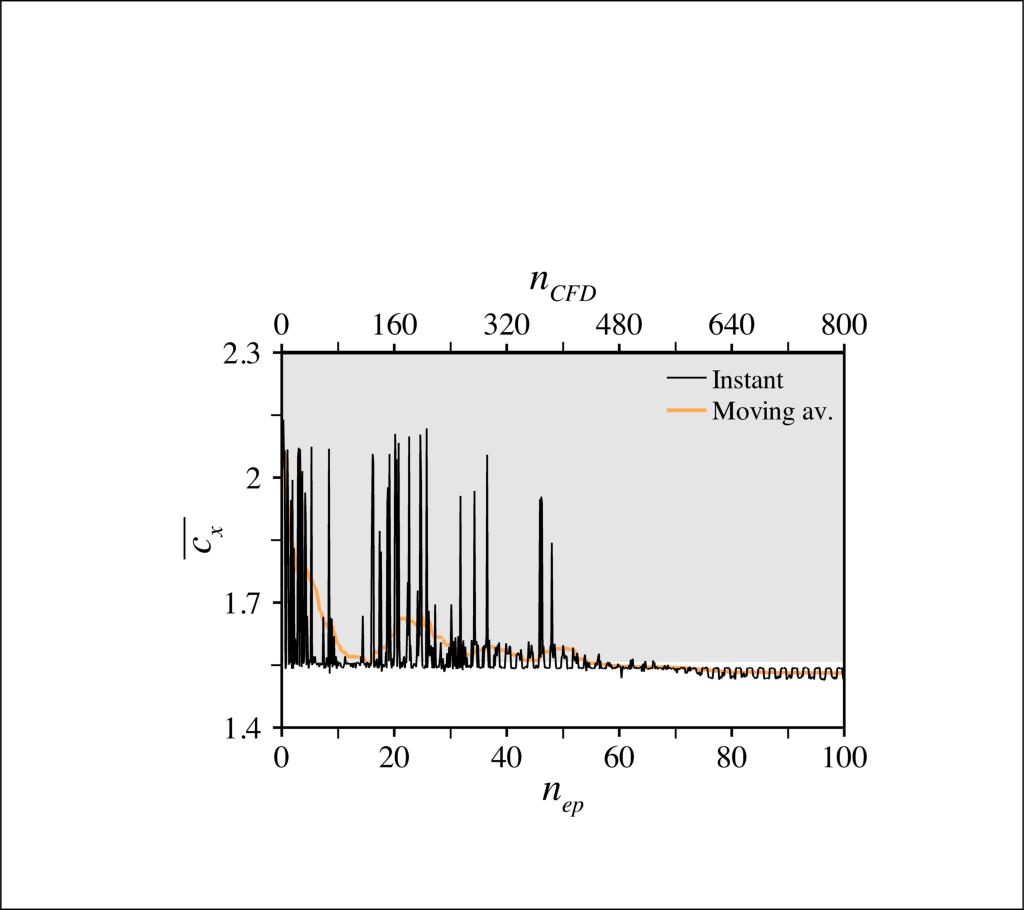}}
\put(7.65,6.25){\includegraphics[trim=175 92.5 155 270pt,clip,height=5.85cm]{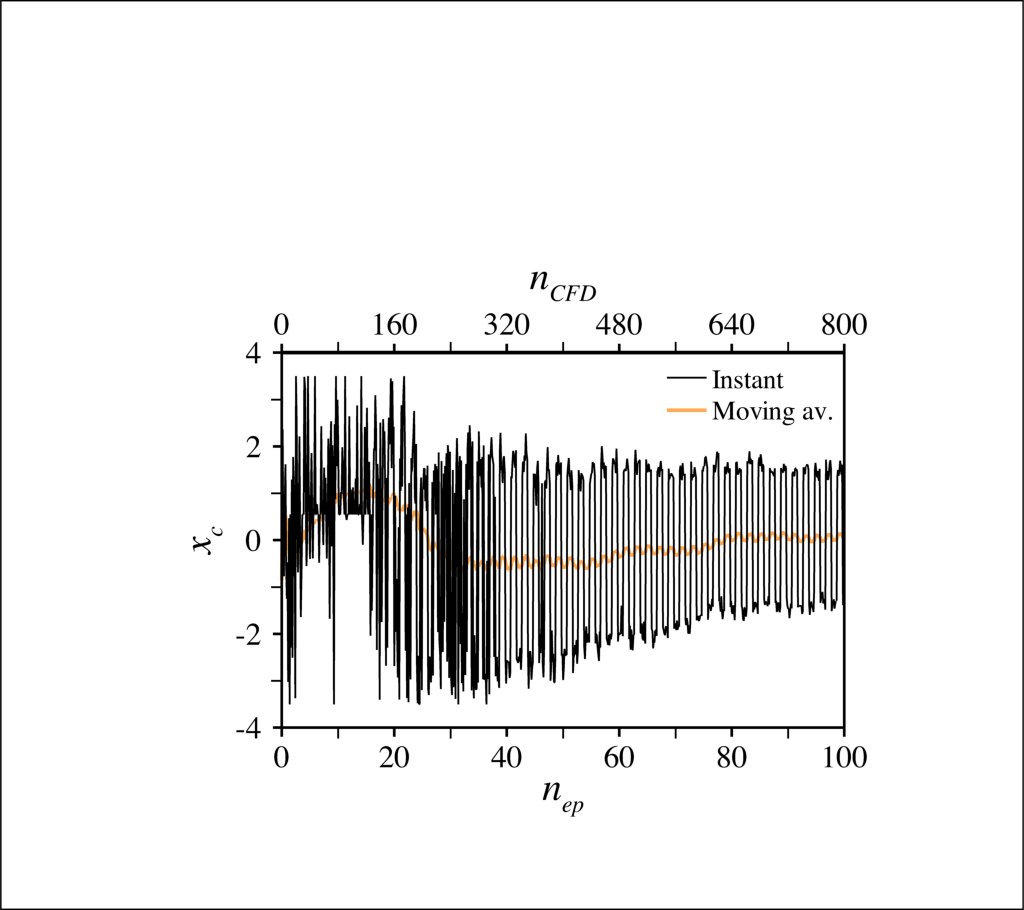}}
\put(0.1,0){\includegraphics[trim=175 92.5 155 270pt,clip,height=5.85cm]{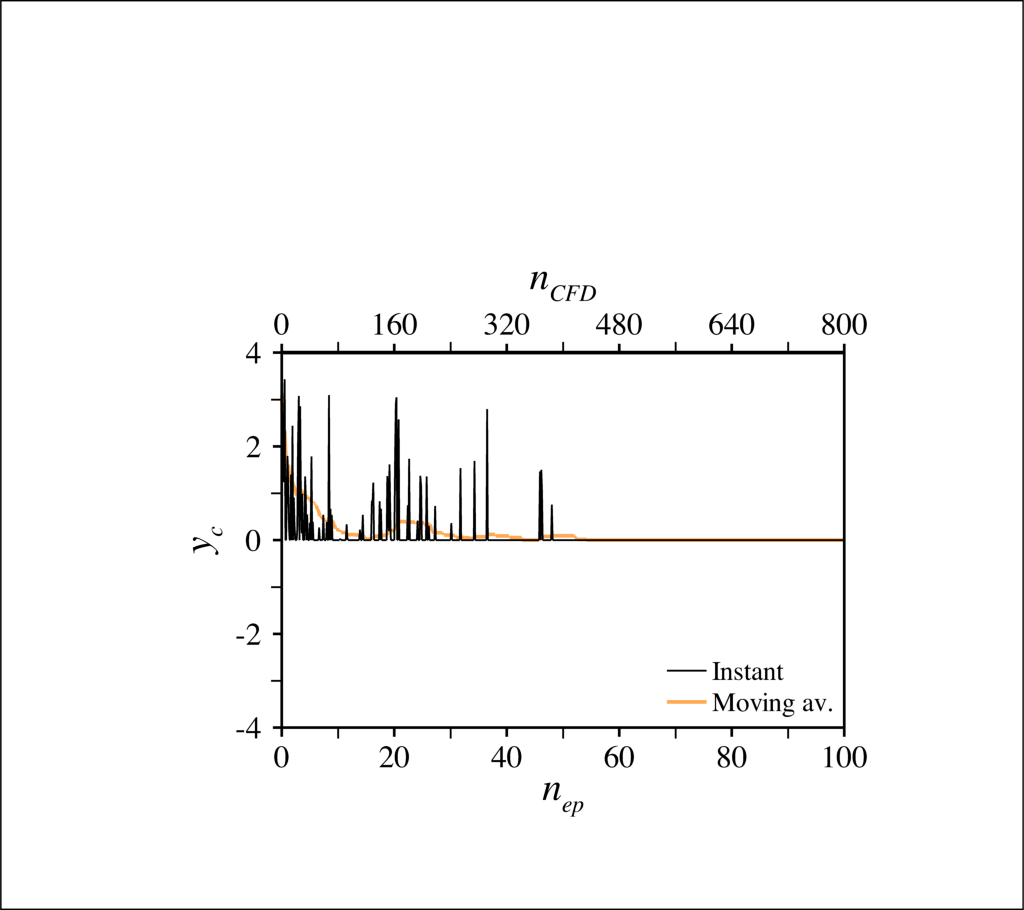}}
\put(7.65,0){\includegraphics[trim=175 92.5 155 340pt,clip,height=5.1cm]{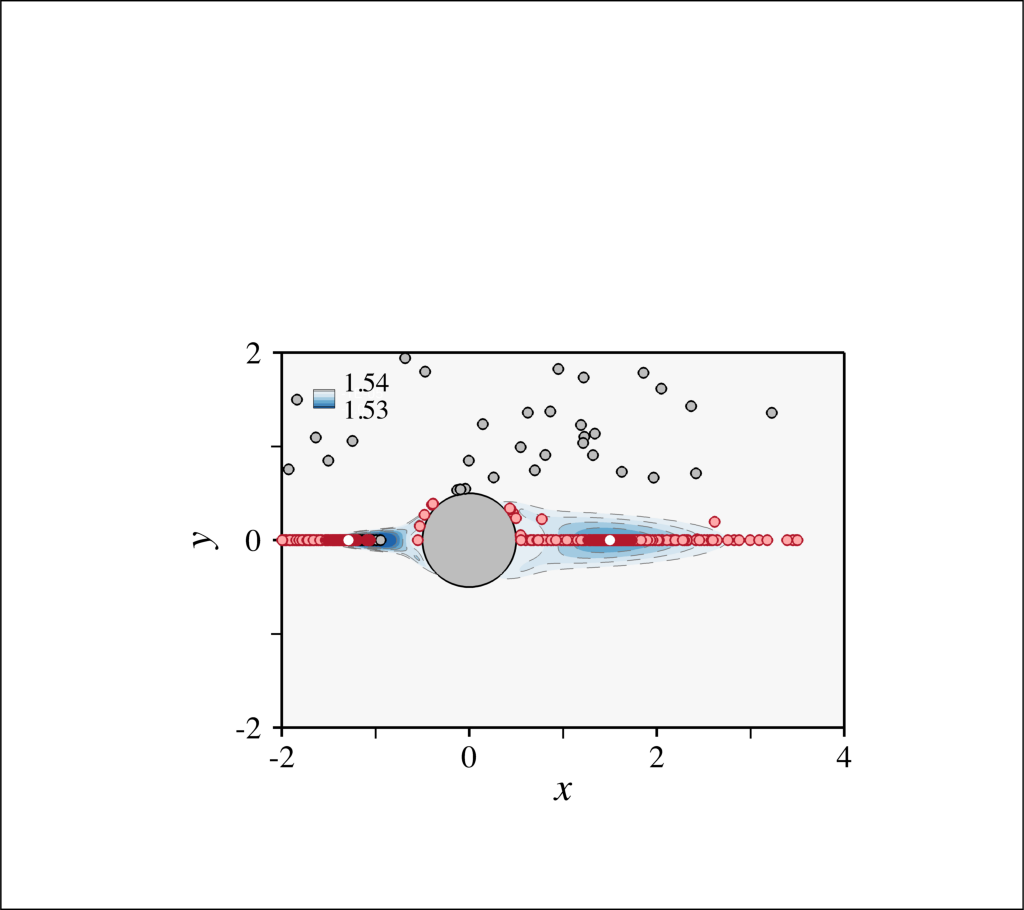}}
\put(4.7,1){\includegraphics[trim=465 275 410 445pt,clip,height=0.8cm]{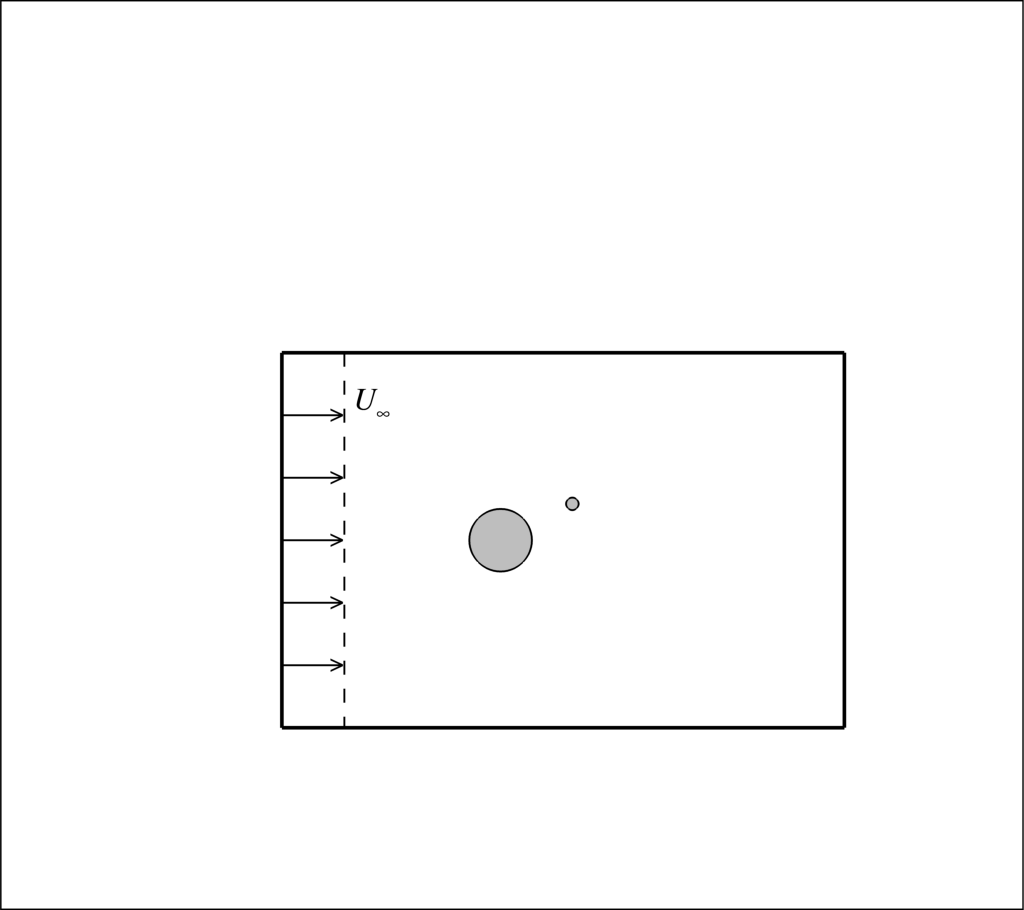}}
\put(4.7,10.37){\includegraphics[trim=465 275 410 445pt,clip,height=0.8cm]{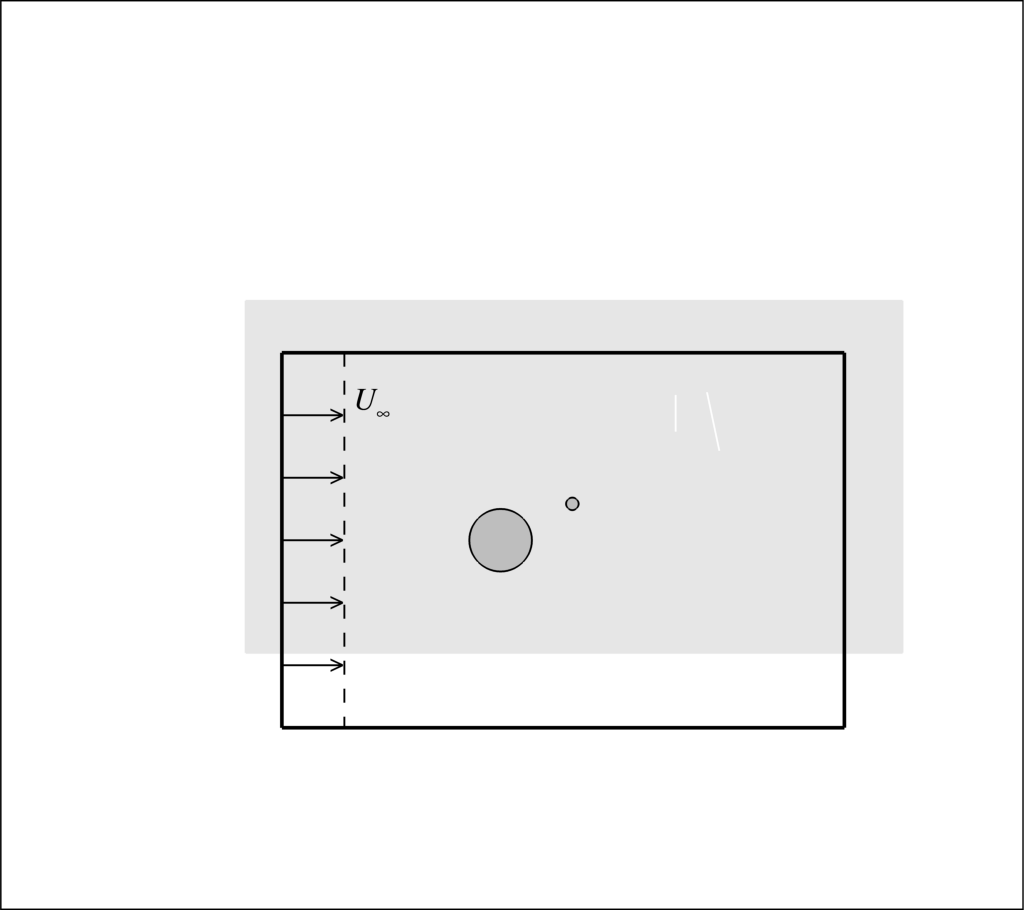}}
\put(12.25,10.37){\includegraphics[trim=465 275 410 445pt,clip,height=0.8cm]{fig9_thumba.png}}
\put(0.1,12.25){(a)}
\put(7.65,12.25){(b)}
\put(0.1,6){(c)}
\put(7.65,6){(d)}
\end{picture}
\caption{Open-loop control of the circular cylinder flow by a small control cylinder of diameter $d=0.1$ at $\rey=40$ - (a) Evolution per episode for the instant (black line) and moving average (over episodes, orange line) values of the mean drag (over time). The uncontrolled drag is at the bottom of the grey shaded area. (b-c) Same as (a) for the $x_c$ and (c) $y_c$ positions of the control cylinder center. (d) Theoretical mean drag variation computed by a steady adjoint method modeling the presence of the control cylinder by a pointwise reacting force localized at the same location where the control cylinder is placed (only the negative iso-contours are reported for clarity). The grey circles are the positions investigated by the DRL. The light red circles are high reward positions spanned over the course of optimization. The dark red circles are those high reward positions spanned over the last 5 episodes. The white circles are the median values reported in the summarizing table~\ref{tab:control}.}
\label{fig:control_re40_ppo}
\end{figure}

For this first case, 100 episodes have been run, which represents 800 simulations, each of which lasts 
$\sim 35$mn on 12 cores, 
hence $\sim 480$h of total CPU cost (equivalently, $\sim60$h of resolution time).
The moving average value of drag reaches a plateau after about 60 episodes in figure~\ref{fig:control_re40_ppo}(a), with the optimal value $\overline{c_x}^{\,\star}=1.53 \pm 0.01$ computed as the average over the 5 latest episodes representing a reduction by roughly $2\%$ with respect to the uncontrolled value $1.56$ (in good agreement with the reference $1.54$ from the literature \cite{Fornberg1980,Henderson1995}). \red{Meanwhile, the instant value of drag actually keeps oscillating over the next 40 episodes with small but finite amplitude}, which is further evidenced in figure~\ref{fig:control_re40_ppo}(b-c) showing the instant and moving average center positions of the control  cylinder. On the one hand, ${y_c}^\star$ quickly settles to zero, i.e., the control cylinder converges to the horizontal centerline. On the other hand, $x_c$ keeps exchanging positions between two regions distributed almost symmetrically on either side of the main cylinder, an upstream region associated with $\overline{c_x}\sim 1.51$ and a slightly less efficient downstream region associated with $\overline{c_x}\sim1.54$, which suggests that the drag functional has global and local minima located in valleys of comparable depth. 
Confirmation comes from the theoretical drag variations computed (in steady mode) from the baseline adjoint method described  
in appendix~\ref{sec:appadj:sub:base}, whose negative iso-values (associated to drag reduction) are mapped in figure~\ref{fig:control_re40_ppo}(d). The latter unveil two regions nestled against either side of the main cylinder and achieving similar drag reduction by $\sim 2\%$, a first one extending upstream over approximately $1$ diameter, and a second one, slightly less efficient and extending downstream and along the outer boundary of the recirculation over $3$ diameters. DRL manages to find high-reward positions in both, which is best seen from the various symbols in figure~\ref{fig:control_re40_ppo}(d) showing the complete set of PPO-1 positions investigated over the course of optimization (grey circles) together with those positions achieving optimal drag reduction within $5\%$ (light red circles), including a few non-centerline positions along the edge of both drag reduction regions. 
Nonetheless, the algorithm ultimately converges to almost symmetrical core positions, as evidenced by the dark red circles in figure~\ref{fig:control_re40_ppo}(d) showing the positions spanned over the 5 latest episodes. \red{Despite limited discrepancies regarding the exact position of the upstream region (slightly shifted upstream in the present approach), this is consistent with the adjoint-based results and} clearly assesses the ability of single-step PPO to identify both regions of interest and to accurately predict the drag reduction achieved in these regions.
 
\subsubsection{Laminar time-dependent regime and circular geometry at \text{Re=100}} \label{sec:olc:sub:control100} 

\begin{figure}[!t]
\setlength{\unitlength}{1cm}
\begin{picture}(20,12.25)
\put(0.1,6.25){\includegraphics[trim=175 92.5 155 270pt,clip,height=5.85cm]{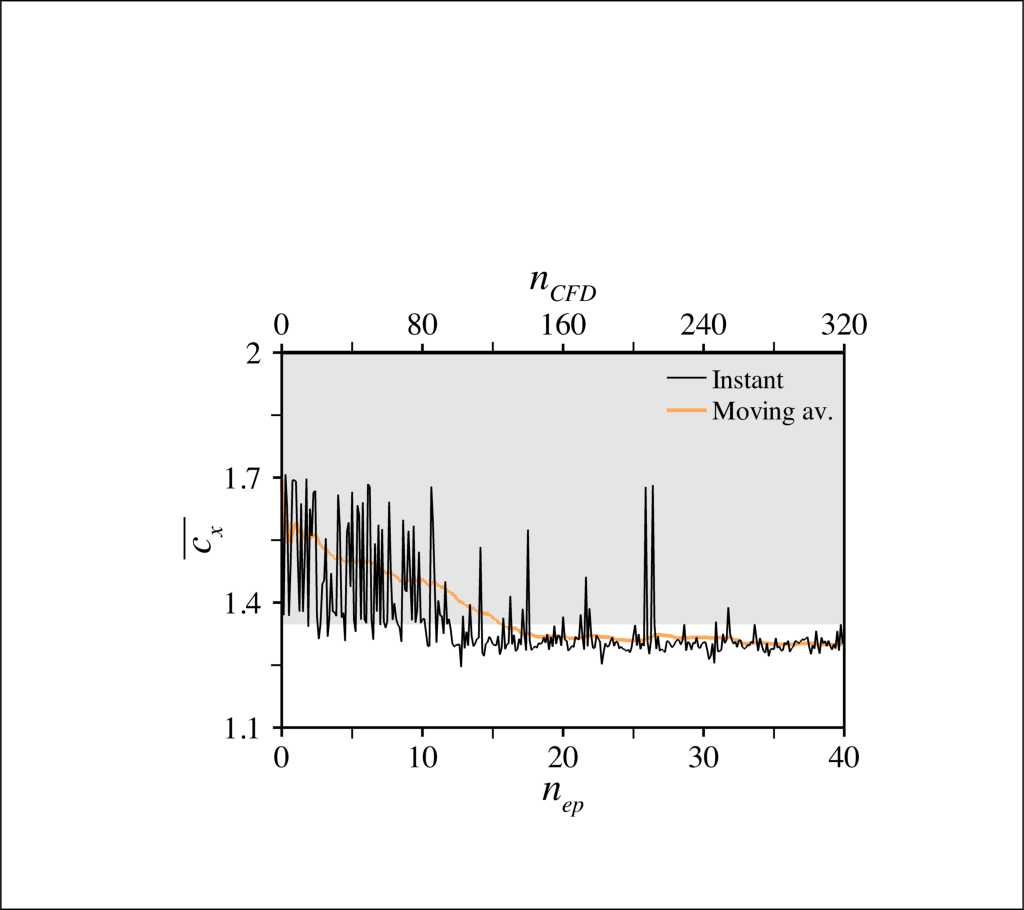}}
\put(7.65,6.25){\includegraphics[trim=175 92.5 155 270pt,clip,height=5.85cm]{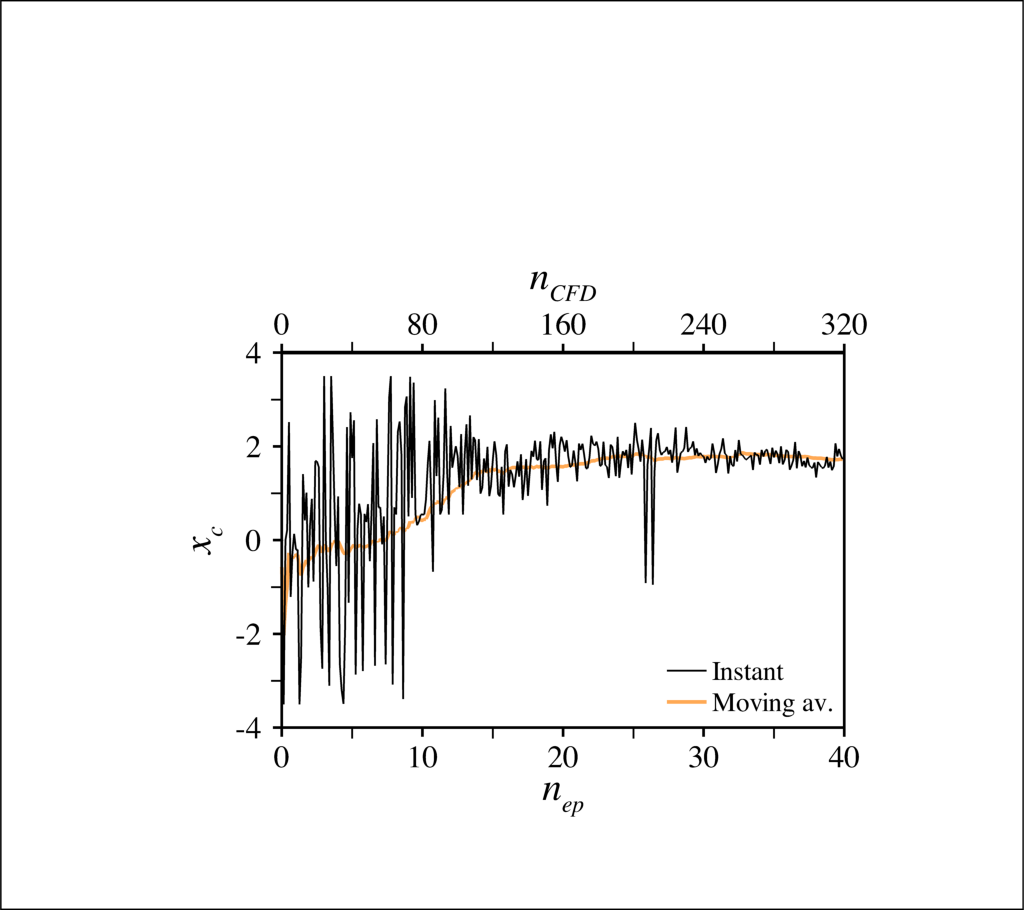}}
\put(0.1,0){\includegraphics[trim=175 92.5 155 270pt,clip,height=5.85cm]{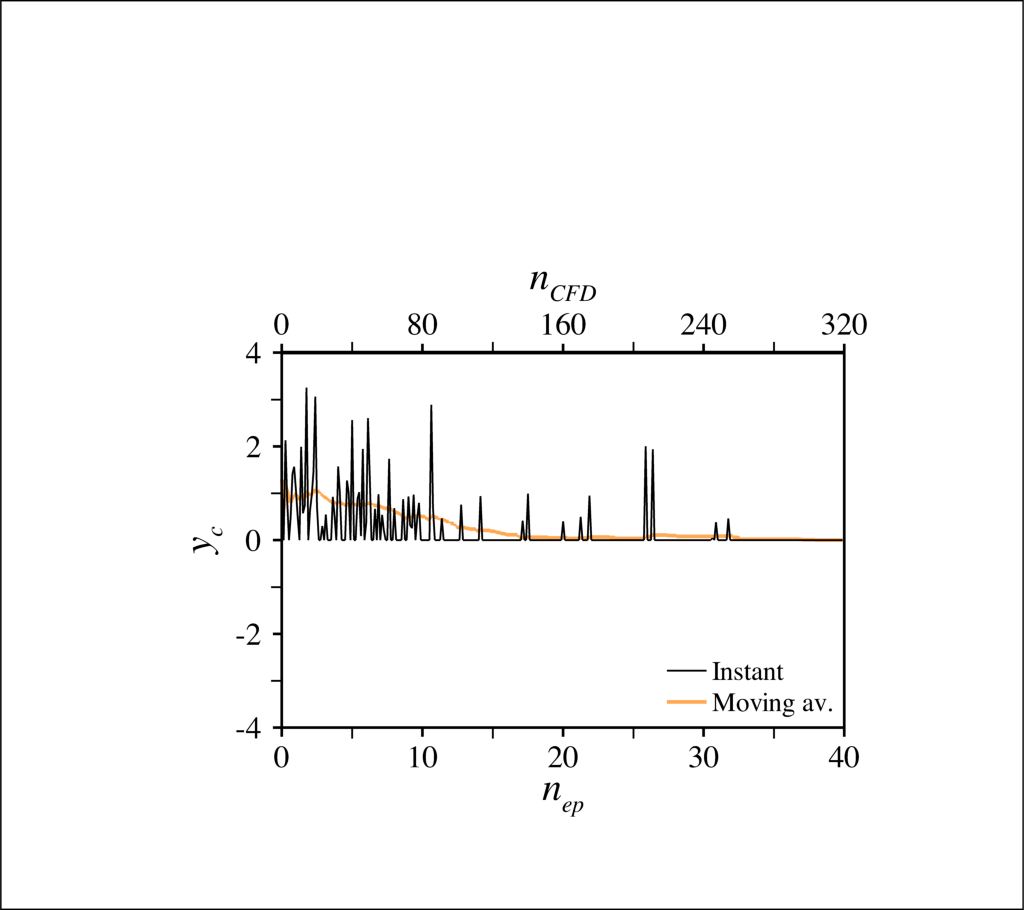}}
\put(7.65,0){\includegraphics[trim=175 92.5 155 340pt,clip,height=5.1cm]{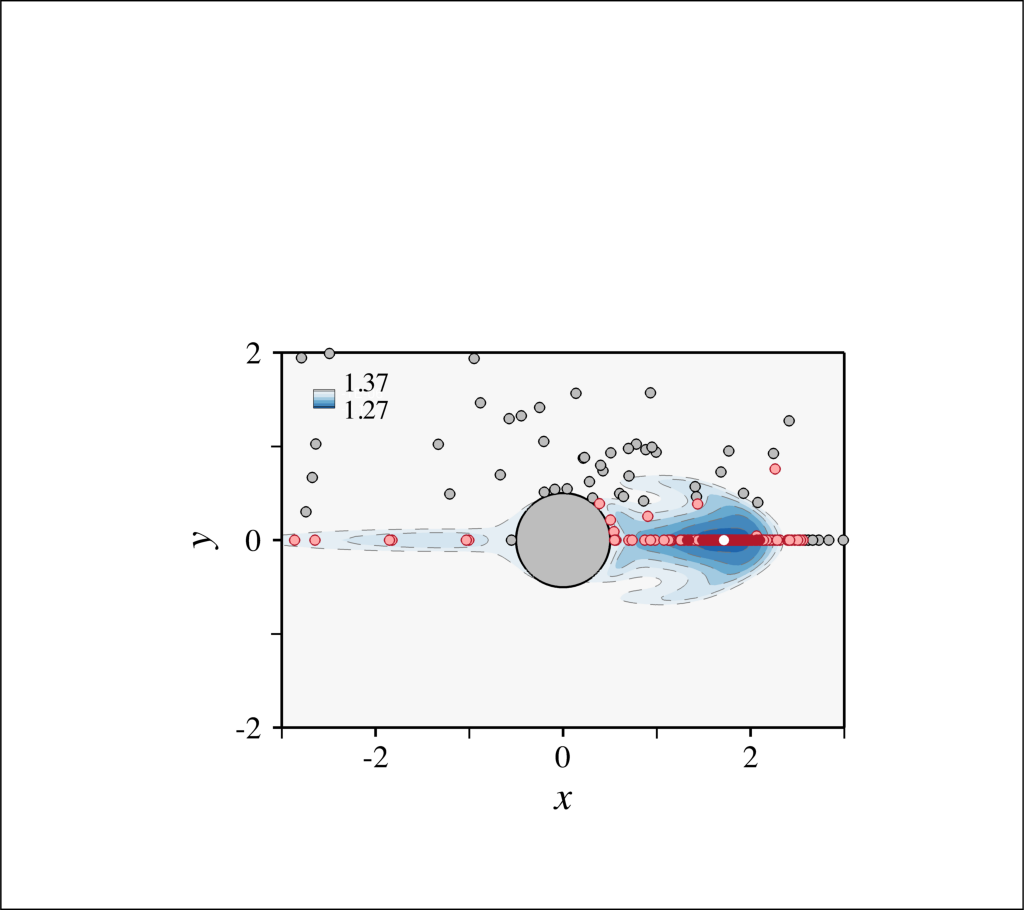}}
\put(4.7,1){\includegraphics[trim=465 275 410 445pt,clip,height=0.8cm]{fig9_thumba.png}}
\put(4.7,10.37){\includegraphics[trim=465 275 410 445pt,clip,height=0.8cm]{fig9_thumbb.png}}
\put(12.25,7.25){\includegraphics[trim=465 275 410 445pt,clip,height=0.8cm]{fig9_thumba.png}}
\put(0.1,12.25){(a)}
\put(7.65,12.25){(b)}
\put(0.1,6){(c)}
\put(7.65,6){(d)}
\end{picture}
\caption{Open-loop control of the circular cylinder flow by a small control cylinder of diameter $d=0.1$ at $\rey=100$ - Same as figure~\ref{fig:control_re40_ppo}, only the theoretical variations in (d) have been computed by the time-varying adjoint method presented in~\cite{meli14}. }
\label{fig:control_re100_ppo}
\end{figure}

For this case, 40 episodes have been run, which represents 320 simulations, 
each of which lasts $\sim 1$h on 12 cores, hence $\sim 320$h of total CPU cost (equivalently, $\sim40$h of resolution time).
The moving average reward plateaus after about 25 episodes in figure~\ref{fig:control_re100_ppo}(a), with the optimal drag $\overline{c_x}^{\,\star}=1.30 \pm 0.01$ computed as the average over the 5 latest episodes representing a reduction by roughly $5\%$ with respect to the uncontrolled value $1.37$ (close to the reference $1.35$ from the literature \cite{Henderson1995}).
Unlike the previous steady case at $\rey=40$, the center position of the control cylinder
exhibits a similarly converging behavior \red{in figure~\ref{fig:control_re100_ppo}(b-c) with ${x_c}^\star=-1.76 \pm 0.03$ and ${y_c}^\star=0$,}
which suggests that the drag functional now has a well-defined global minimum.
Confirmation comes from the theoretical drag variations computed (in unsteady mode) from the baseline adjoint method, whose negative iso-values mapped in figure~\ref{fig:control_re100_ppo}(d) are reproduced from~\cite{meli17prf}. The latter unveil again two regions nestled against either side of the main cylinder, a first one extending upstream over approximately $2$ diameter (more than at $\rey=40$), 
 and a second one extending downstream and along the outer boundary of the mean recirculation over $2$ diameters (less than at $\rey=40$). Drag is reduced by roughly $2\%$ upstream, but almost $8\%$ downstream, meaning that the drag functional has global and local minima in valleys of different depth, in line with the DRL results. 
Again, DRL finds high-reward positions in both regions, as evidenced in figure~\ref{fig:control_re100_ppo}(d) by the complete set of PPO-1 positions investigated over the course of optimization (small grey circles) and the positions achieving optimal drag reduction within $5\%$ (light red circles), including a few centerline upstream positions. The algorithm however quickly settles for the most efficient downstream region, as the positions spanned over the 5 latest episodes (dark red circles) all lie in the core of the mean recirculation region, in striking agreement with the adjoint-based results. 

\subsubsection{Turbulent regime and circular geometry at \text{Re=3900}} \label{sec:olc:sub:control3900} 

\begin{figure}[!t]
\setlength{\unitlength}{1cm}
\begin{picture}(20,12.25)
\put(0.1,6.25){\includegraphics[trim=175 92.5 155 270pt,clip,height=5.85cm]{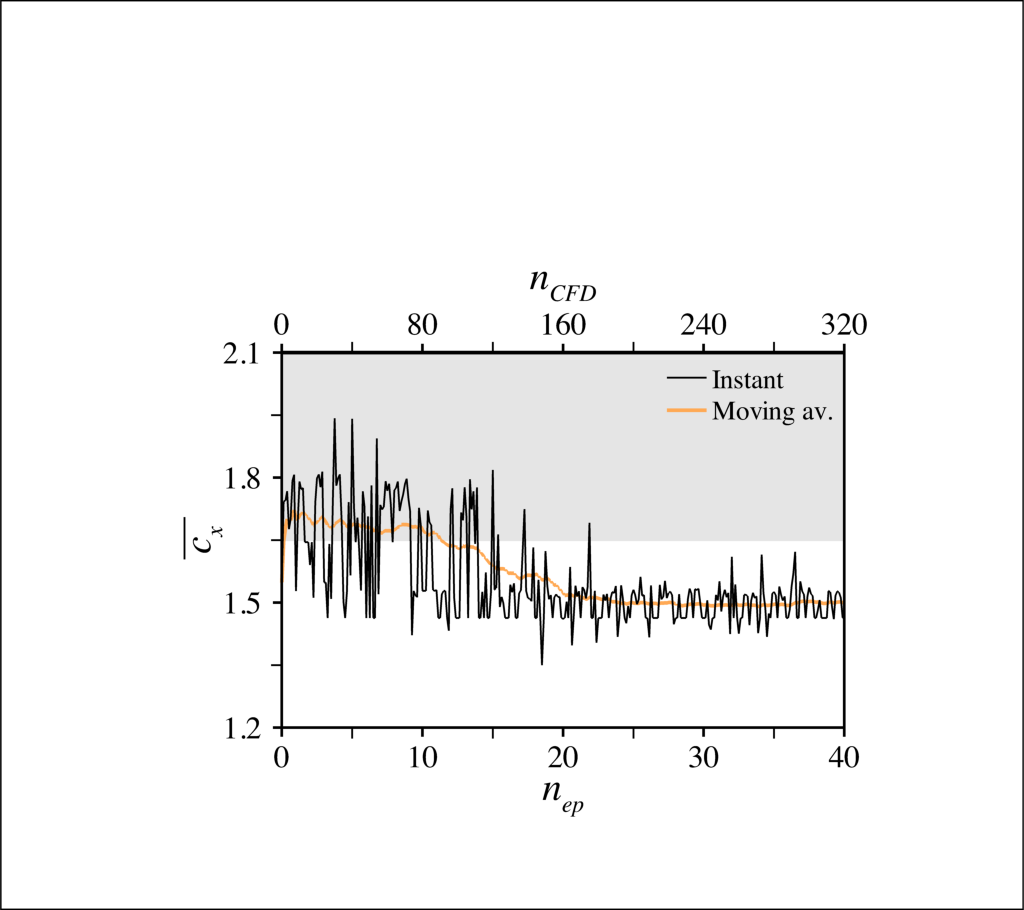}}
\put(7.65,6.25){\includegraphics[trim=175 92.5 155 270pt,clip,height=5.85cm]{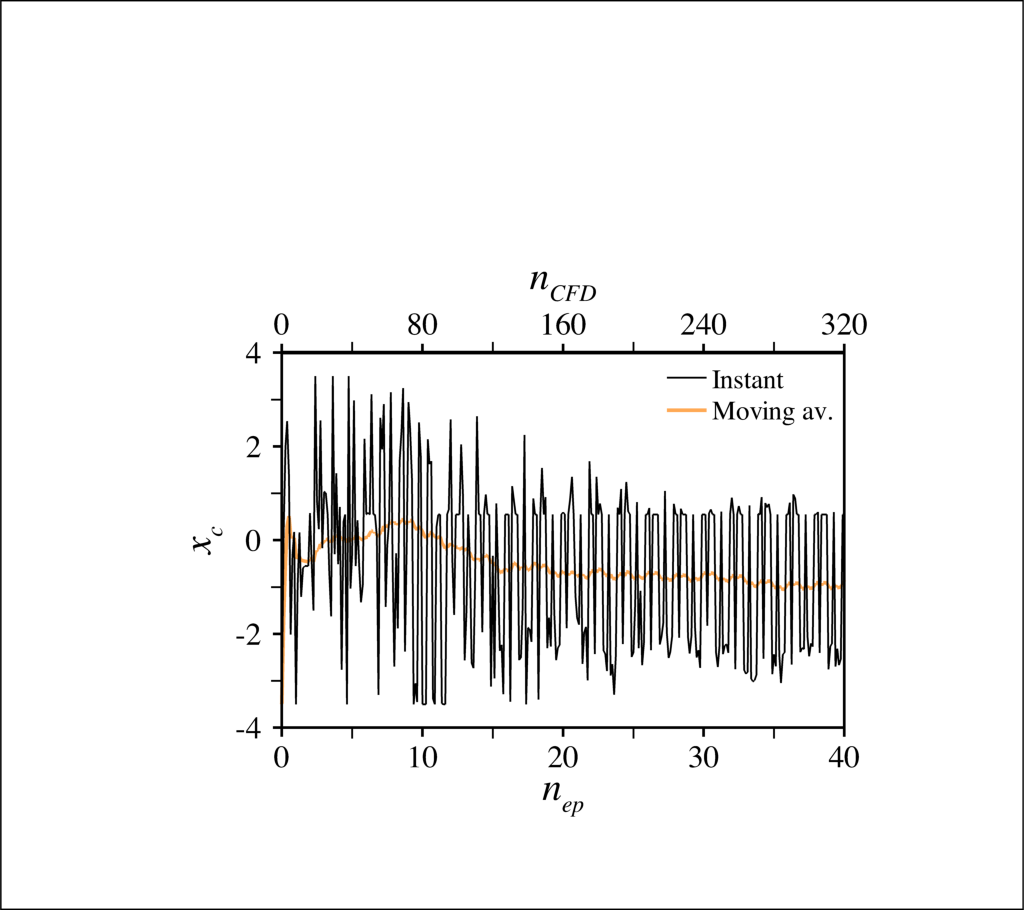}}
\put(0.1,0){\includegraphics[trim=175 92.5 155 270pt,clip,height=5.85cm]{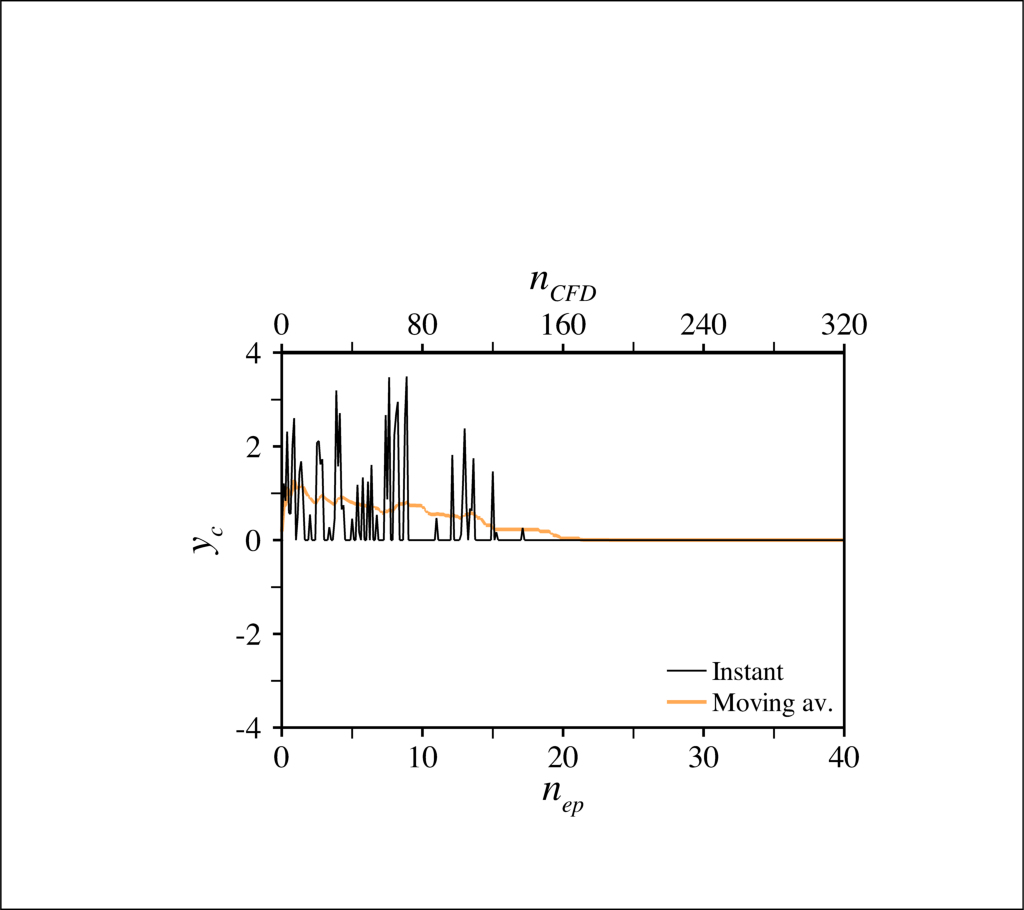}}
\put(7.65,0){\includegraphics[trim=175 92.5 155 340pt,clip,height=5.1cm]{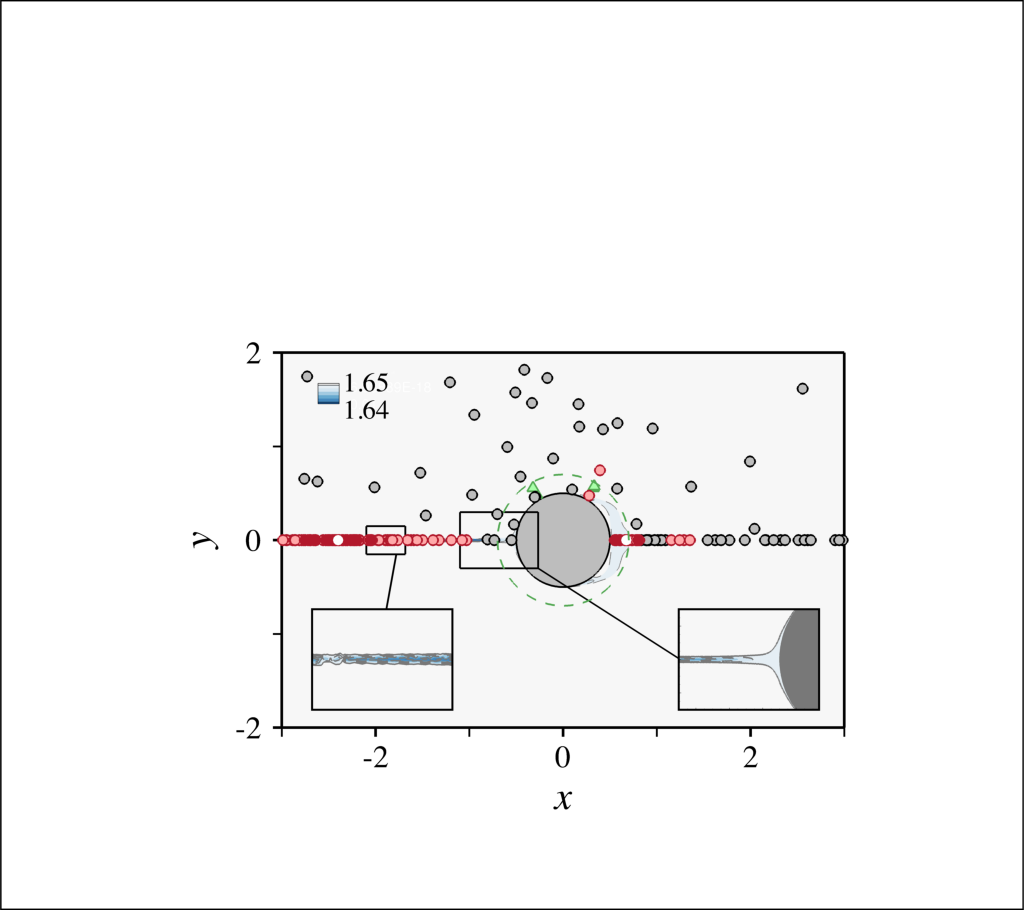}}
\put(4.7,1){\includegraphics[trim=465 275 410 445pt,clip,height=0.8cm]{fig9_thumba.png}}
\put(4.7,10.37){\includegraphics[trim=465 275 410 445pt,clip,height=0.8cm]{fig9_thumbb.png}}
\put(12.25,10.37){\includegraphics[trim=465 275 410 445pt,clip,height=0.8cm]{fig9_thumba.png}}
\put(0.1,12.25){(a)}
\put(7.65,12.25){(b)}
\put(0.1,6){(c)}
\put(7.65,6){(d)}
\end{picture}
\caption{Open-loop control of the circular cylinder flow by a small control cylinder of diameter $d=0.1$ at $\rey=3900$ - Same as figure~\ref{fig:control_re40_ppo}, only the theoretical variations in (d) have been computed by the (steady) simplified mean-flow adjoint method presented in~\cite{meli14}. 
The green dashed circle in (d) indicates the range of center positions spanned experimentally in~\cite{saka94}, with the green triangles marking the sets of positions found to optimally reduce the drag of the main cylinder only.}
\label{fig:control_re3900_ppo}
\end{figure} 

Another 40 episodes have been run for this case, which represents 320 simulations, 
each of which lasts $\sim2$h$30$ on 12 cores (much longer than at $\rey=100$ due to the halved time step), hence $\sim 800$h of total CPU cost (equivalently, $\sim100$h of resolution time).
After about 20 episodes, the moving average reward in figure~\ref{fig:control_re3900_ppo}(a) converges to $\overline{c_x}^{\,\star}=1.50\pm0.01$, which represents a reduction of drag by $9\%$ with respect to the uncontrolled value $1.65$ (in good agreement with reference 2-D RANS data from the literature~\cite{Pereira2015}). 
\red{The center position of the control cylinder however keeps oscillating over the next 15 episodes in figure~\ref{fig:control_re3900_ppo}(b-c), as ${y_c}^\star$ goes to zero but 
$x_c$ exchanges positions 
 between two regions located on either side of the main cylinder, an upstream region associated with $\overline{c_x}\sim1.52$ and a downstream
region associated with $\overline{c_x}\sim1.46$. This suggests that the drag functional has global and local minima in valleys of comparable depth, which is reminiscent of the steady case at $\rey=40$, only the deepest valley is 
now downstream, not upstream.}
\red{Interestingly, Ref.~\cite{saka94} determines experimentally different optimal positions
 $(G,\theta)=(0.14-0.16,60^\circ)$ and $(0.06-0.14,115^\circ)$, shown as the green triangles in figure~\ref{fig:control_re3900_ppo}(d). Additional DNS runs have thus been carried out to confirm sub-optimality for our case, although the algorithm does identify a couple of high-reward positions in the vicinity of the downstream experimental region. This probably stems from the noticeable differences between both studies, as the Reynolds number in~\cite{saka94} is larger by one order of magnitude ($\rey=65000$), the control cylinder is almost twice as small ($d=0.06$), and the experiments focus on the drag of the main cylinder (not the total drag) while spanning a much smaller range of center positions (indicated by the green dashed circle in figure~\ref{fig:control_re3900_ppo}(d)).}

The DRL results are conversely qualitatively in line with the negative iso-values of the adjoint-based drag variations shown in figure~\ref{fig:control_re3900_ppo}(d). Those indicate that drag is reduced in two distinct regions nestled against either side of the main cylinder, a first narrow one extending upstream along the centerline over approximately $2$ diameters, and a second one extending downstream over a half-diameter and in the vicinity of the mean separation points. Nonetheless, the agreement is not quantitative, as the theoretical variations are by a mere $1\%$ upstream (and even lower downstream).
This is most likely because all theoretical variations have been modeled after a 
simplified adjoint method intended to guide near-optimal design with marginal computational effort (as it requires knowledge of the sole mean uncontrolled solution, as explained in appendix~\ref{sec:appadj:sub:mean}), that ends up miscalculating the effect of the control cylinder because of an insufficient level of sophistication. On the one hand, the marginal size of the downstream region (as well as the marginal drag reduction predicted in this region) is ascribed to the fact that the approach has been shown to possibly miss out on sensitivity regions 
involving strong interactions of the mean and fluctuating solution components via the formation of Reynolds stresses~\cite{meli17prf}: the mean recirculation is one such region where reducing the drag of the main cylinder, even by a small amount, suffices to reduce the total drag because the $x$ velocity is negative and the control cylinder is thus a source of thrust, not drag.
On the other hand, the outcome in the upstream region is sensitive to the force model used to mimic the effect of the control cylinder, as it turns out its drag balances almost exactly the amount by which it reduces the drag of the main cylinder.
\red{The weak upstream control efficiency may thus be due to the fact that the simplified adjoint method considers only the mean component of the force acting on the control cylinder, but overlooks the potential for additional drag reduction via the fluctuating component.}
Moreover, this is a region where the control cylinder likely induces strong mean flow modifications because the local inhomogeneity length scale becomes smaller than the diameter of the control cylinder, which in turn may invalidate the linear assumption inherent to the adjoint method (the retained diameter $d=0.1$ is a compromise between smallness and cost control, as implementing a smaller control cylinder would require increasing the number of grid points and decreasing the time-step to capture properly the wake of the control cylinder).

\begin{figure}[!t]
\setlength{\unitlength}{1cm}
\begin{picture}(20,12.25)
\put(0.1,6.25){\includegraphics[trim=175 92.5 155 270pt,clip,height=5.85cm]{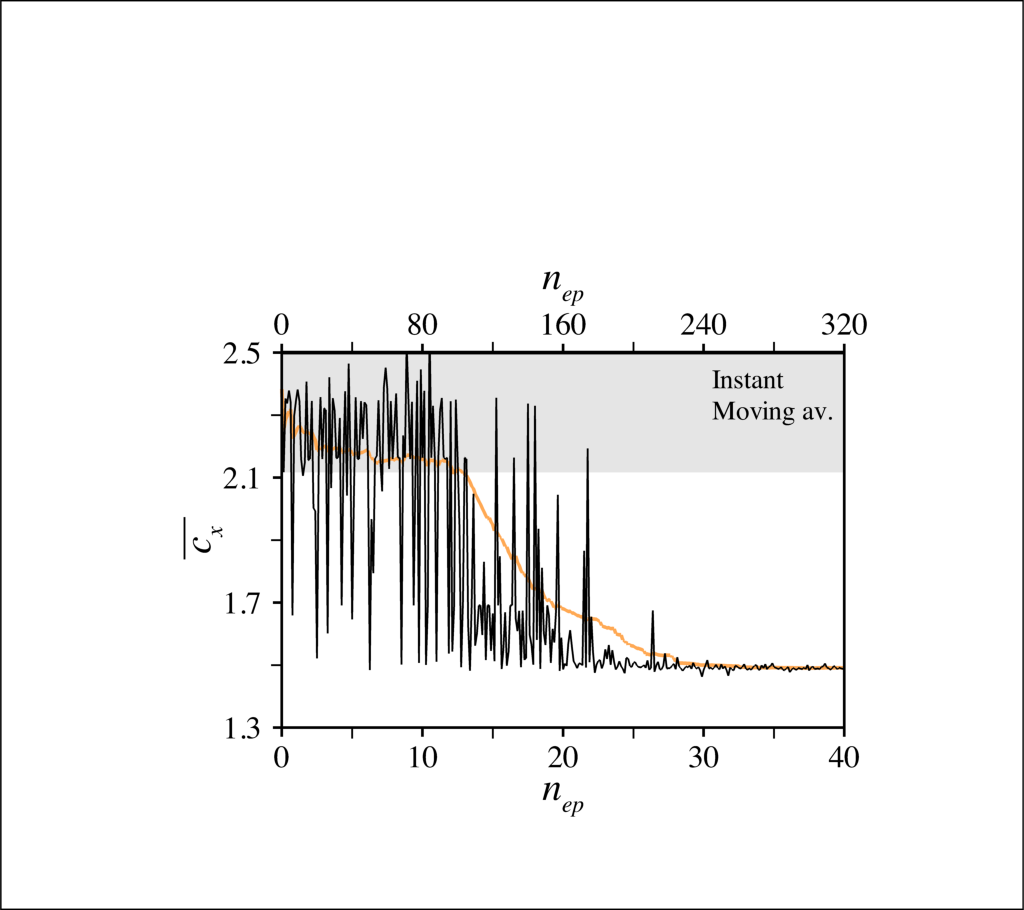}}
\put(7.65,6.25){\includegraphics[trim=175 92.5 155 270pt,clip,height=5.85cm]{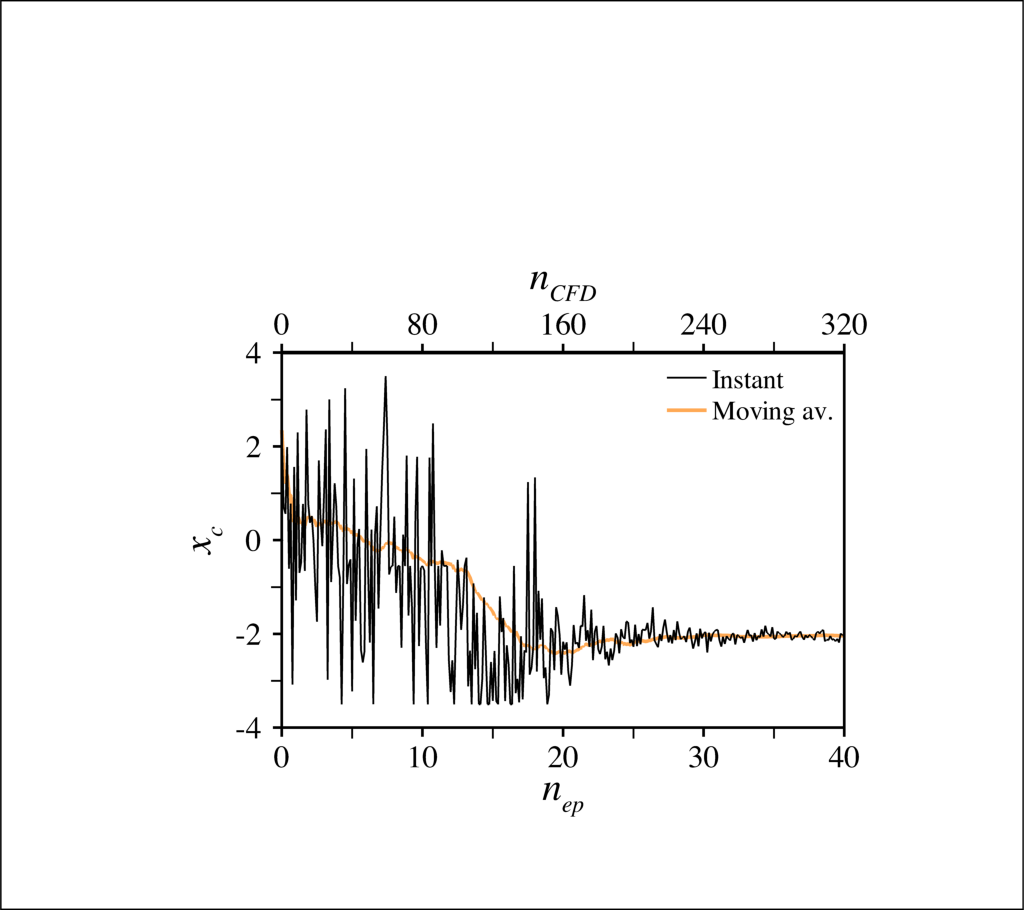}}
\put(0.1,0){\includegraphics[trim=175 92.5 155 270pt,clip,height=5.85cm]{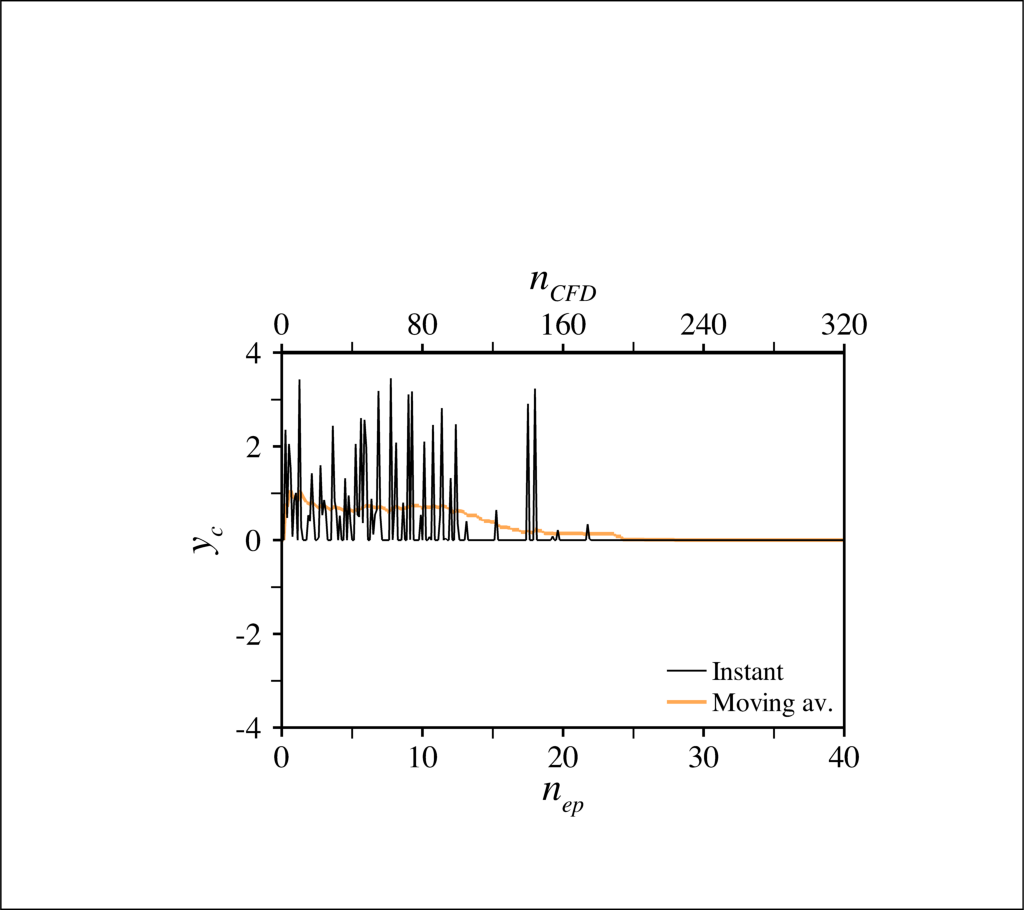}}
\put(7.65,0){\includegraphics[trim=175 92.5 155 340pt,clip,height=5.1cm]{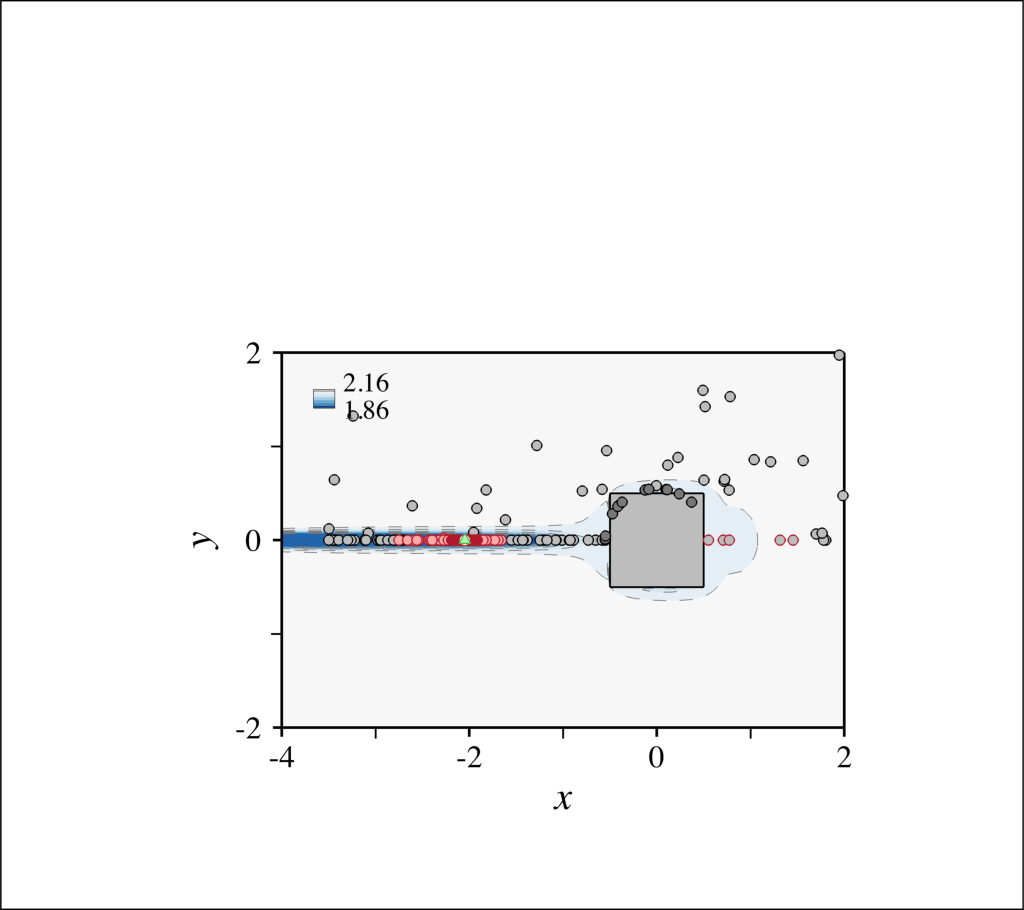}}
\put(4.7,1){\includegraphics[trim=465 275 410 445pt,clip,height=0.8cm]{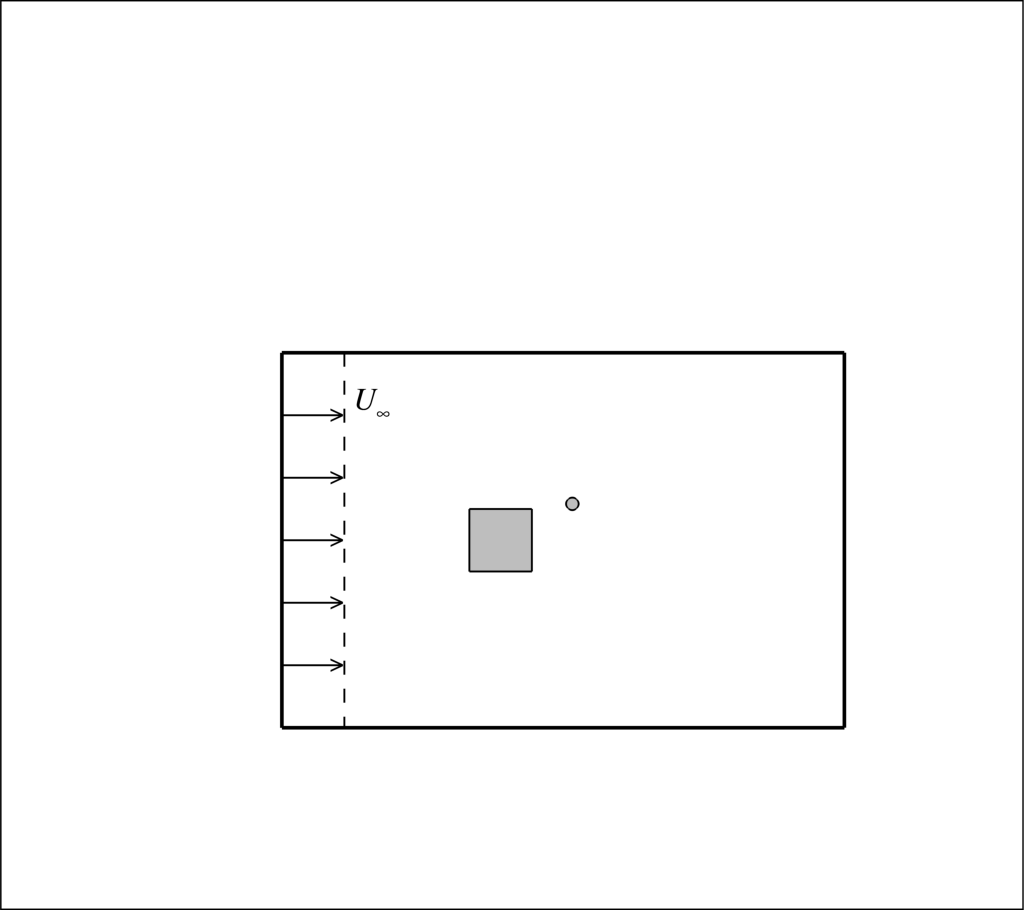}}
\put(4.7,10.37){\includegraphics[trim=465 275 410 445pt,clip,height=0.8cm]{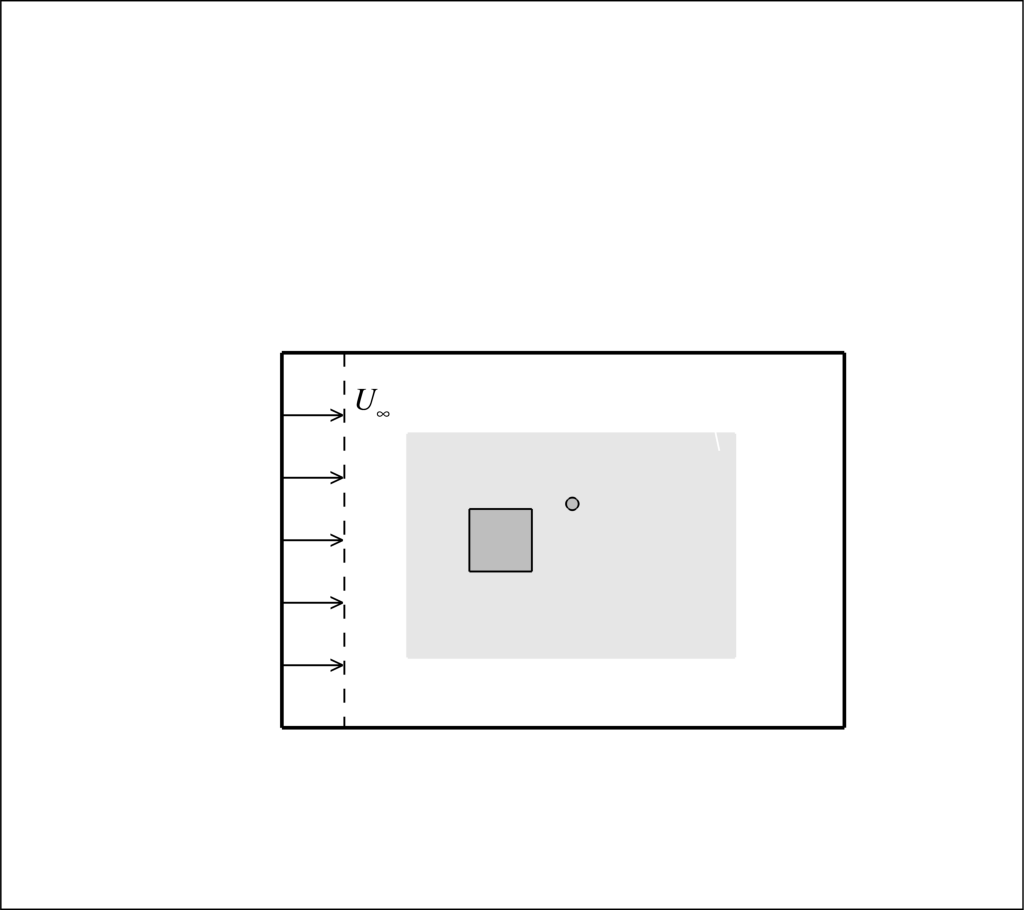}}
\put(12.25,10.37){\includegraphics[trim=465 275 410 445pt,clip,height=0.8cm]{fig12_thumba.png}}
\put(0.1,12.25){(a)}
\put(7.65,12.25){(b)}
\put(0.1,6){(c)}
\put(7.65,6){(d)}
\end{picture}
\caption{Same as figure~\ref{fig:control_re3900_ppo} for the open-loop control of the square cylinder flow by a small control cylinder of diameter $d=0.1$ at $\rey=22000$. In (d), the grey symbols circled in red are downstream sub-optimal position. The dark circles are dismissed positions for which the control and main cylinders intersect, and the green triangle marks the experimental position found in~\cite{Igarashi1997} to optimally reduce the drag of the two-cylinder system at $\rey=32000$.}
\label{fig:control_re22000_ppo}
\end{figure} 

\subsubsection{Turbulent regime and square geometry at \text{Re=22000}} \label{sec:olc:sub:control22000}

In order to push the comparison further, additional calculations have been undertaken for a square geometry of the main cylinder, whose larger upstream sensitivity yields more clear-cut control efficiency, as can be inferred from the results in~\cite{meli14,meli17prf}. This is because 
\red{the blunt square geometry strengthens the upstream pressure gradient (compared to its bluff circular shape). In return, the gap flow velocity between the two cylinders decreases and so does the drag of the control cylinder, hence a boost in efficiency that helps mitigate the issue of sensitivity to the force model.}

Another 40 episodes have been run for this case, which represents 320 simulations, 
each of which lasts $\sim3$h$20$ on 12 cores, hence $\sim 1020$h of total CPU cost (equivalently, $\sim130$h of resolution time).
One difficulty for this case is that the main and control cylinders can intersect each other under mapping~\myrefeq{eq:mapcontrol}, in which case it has been found relevant to simply discard the CFD and force the reward to its uncontrolled value. 
The moving average reward plateaus after about 30 episodes in figure~\ref{fig:control_re22000_ppo}(a), with the optimal drag $\overline{c_x}^{\,\star}=1.49 \pm 0.01$ computed as the average over the 5 latest episodes representing a reduction by $30\%$ with respect to the uncontrolled value $2.16$ (close to the reference $2.1-2.2$ from the literature \cite{Iaccarino2003,Rodi1997Workshop}). The center position of the control cylinder
exhibits a similarly converging behavior in figures~\ref{fig:control_re22000_ppo}(b-c) with ${x_c}^\star=-2.04 \pm 0.02$, and ${y_c}^\star=0$, 
which suggests that the drag functional has a well-defined global minimum.
This is in excellent agreement with~\cite{Igarashi1997} reporting
experimental reduction of the total drag by $30\%$ 
inserting control cylinders of comparable sizes upstream of the main cylinder at a slightly different Reynolds number $\rey = 32 000$ (the optimal reported position for $d=0.1$ being ${x_c}^\star\sim-2.0$).
This is also in line with the theoretical drag variations computed from the same simplified adjoint method as in section~\ref{sec:olc:sub:control3900}, whose negative iso-values mapped in figure~\ref{fig:control_re22000_ppo}(d) are reproduced from~\cite{meli14}. 
The latter unveil a main region of interest, that extends upstream over approximately $4$ diameters, and in which 
drag is reduced by almost $20\%$, which represents a satisfactory qualitative and quantitative compliance with the present PPO-1 results.
Drag is also reduced in a second region originating from the separation points (pinned here at the front edges), that extends downstream and along the outer boundary of the mean recirculation over 1 diameter (similar to what has been found using a circular geometry of the main cylinder).
It is worth noticing that the algorithm does identify sub-optimal positions in this region (shown in figure~\ref{fig:control_re22000_ppo}(d) as the grey symbols circled in red). Also, a couple of other low-efficiency PPO-1 positions lie further downstream, which is consistent with the idea that the simplified adjoint method may miss on additional drag reduction occurring via the formation of Reynolds stresses (this is not true of the upstream drag reduction region, whose flow is essentially steady, except for low-amplitude oscillations in the gap flow between the two cylinders).

\begin{table}[!t]
\begin{center}
\begin{threeparttable}
\begin{tabular}{ccccccccccccccccc}
\toprule
\multicolumn{1}{r}{\multirow{1}{*}{\makecell[r]{\raisebox{-\totalheight}{\includegraphics[trim=415 220 425 500pt,clip,height=1cm]{fig9_thumba.png}}}}}& \multicolumn{1}{r}{\makecell[r]{$\overline{c_x}$}} & \multicolumn{1}{p{0.6cm}}{\makecell[r]{$x_c$\tnote{a}}} & \multicolumn{1}{p{0.3cm}}{\makecell[r]{$y_c$}} & \multicolumn{1}{p{0.2cm}}{\makecell[r]{}} & \multicolumn{1}{r}{\makecell[r]{$\overline{c_x}$}} & \multicolumn{1}{p{0.6cm}}{\makecell[r]{$x_c$\tnote{a}}} & \multicolumn{1}{p{0.3cm}}{\makecell[r]{$y_c$}} & \multicolumn{1}{p{0.2cm}}{\makecell[r]{}} & \multicolumn{1}{r}{\makecell[r]{$\overline{c_x}$}} & \multicolumn{1}{p{0.6cm}}{\makecell[r]{$x_c$\tnote{a}}} & \multicolumn{1}{p{0.3cm}}{\makecell[r]{$y_c$}} & \multicolumn{1}{r}{\multirow{3}{*}{\makecell[r]{\raisebox{-\totalheight}{\includegraphics[trim=390 220 425 500pt,clip,height=1cm]{fig12_thumba.png}}}}} & \multicolumn{1}{r}{\makecell[r]{$\overline{c_x}$}} & \multicolumn{1}{p{0.6cm}}{\makecell[r]{$x_c$\tnote{a}}} & \multicolumn{1}{p{0.3cm}}{\makecell[r]{$y_c$}}\\
\cmidrule(lr){2-4}\cmidrule(lr){6-8}\cmidrule(lr){10-12}\cmidrule(lr){14-16}
& \multicolumn{1}{l}{$1.51$} & \multicolumn{1}{r}{$-1.29$} & \multicolumn{1}{l}{$0$} & & \multicolumn{1}{l}{\multirow{2}{*}{$1.30$}} & \multicolumn{1}{l}{\multirow{2}{*}{\white{-}$1.72$}} & \multicolumn{1}{l}{\multirow{2}{*}{$0$}} & & \multicolumn{1}{l}{$1.46$} & \multicolumn{1}{r}{$0.68$}  & \multicolumn{1}{l}{$0$} & & \multicolumn{1}{l}{\multirow{2}{*}{$1.49$}} & \multicolumn{1}{l}{\multirow{2}{*}{$-2.0$}} & \multicolumn{1}{l}{\multirow{2}{*}{$0$}} & \multicolumn{1}{l}{\multirow{2}{*}{\quad Optimal}}\\
& \multicolumn{1}{l}{$\gray{1.54}$} & \multicolumn{1}{r}{$\gray{1.50}$} & \multicolumn{1}{l}{$\gray{0}$} & & & & & & \multicolumn{1}{l}{$\gray{1.52}$} & \multicolumn{1}{r}{$\gray{-2.40}$}  & \multicolumn{1}{l}{$\gray{0}$} \\
\cmidrule(lr){1-17}
\multicolumn{17}{r}{CFD} \\
\cmidrule(lr){1-17}
\multicolumn{4}{r}{40} & \multicolumn{4}{r}{100} & \multicolumn{4}{r}{3900} & \multicolumn{4}{r}{22000} & \multicolumn{1}{l}{\quad Reynolds number} \\
\multicolumn{4}{r}{0.1} & \multicolumn{4}{r}{\guillemotright} & \multicolumn{4}{r}{0.05} & \multicolumn{4}{r}{\guillemotright} & \multicolumn{1}{l}{\quad Time step} \\
\multicolumn{4}{r}{$[150;150]$} & \multicolumn{4}{r}{$[100;200]$} & \multicolumn{4}{r}{\guillemotright} & \multicolumn{4}{r}{\guillemotright} & \multicolumn{1}{l}{\quad Averaging time span} \\
\multicolumn{4}{r}{$[-15;40]\times[-15;15]$} & \multicolumn{4}{r}{\guillemotright} & \multicolumn{4}{r}{\guillemotright} & \multicolumn{4}{r}{$[-6;15]\times[-7;7]$} & \multicolumn{1}{l}{\quad Mesh dimensions} \\
\multicolumn{4}{r}{150000} & \multicolumn{4}{r}{\guillemotright} & \multicolumn{4}{r}{\guillemotright} & \multicolumn{4}{r}{190000} & \multicolumn{1}{l}{\quad Nb. mesh elements} \\
\multicolumn{4}{r}{0.001} & \multicolumn{4}{r}{\guillemotright} & \multicolumn{4}{r}{\guillemotright} & \multicolumn{4}{r}{\guillemotright} & \multicolumn{1}{l}{\quad Interface $\perp$ mesh size} \\
\multicolumn{4}{r}{12} & \multicolumn{4}{r}{\guillemotright} & \multicolumn{4}{r}{\guillemotright} & \multicolumn{4}{r}{\guillemotright} & \multicolumn{1}{l}{\quad Nb. Cores} \\
\cmidrule(lr){1-17}
\multicolumn{17}{r}{PPO-1} \\
\cmidrule(lr){1-17}
\multicolumn{4}{r}{100} & \multicolumn{4}{r}{40} & \multicolumn{4}{r}{\guillemotright} & \multicolumn{4}{r}{\guillemotright} & \multicolumn{1}{l}{\quad Nb. episodes} \\
\multicolumn{4}{r}{8} & \multicolumn{4}{r}{\guillemotright} & \multicolumn{4}{r}{\guillemotright} & \multicolumn{4}{r}{\guillemotright} & \multicolumn{1}{l}{\quad Nb. environments} \\
\multicolumn{4}{r}{32} & \multicolumn{4}{r}{\guillemotright} & \multicolumn{4}{r}{\guillemotright} & \multicolumn{4}{r}{\guillemotright} & \multicolumn{1}{l}{\quad Nb. epochs} \\
\multicolumn{4}{r}{1} & \multicolumn{4}{r}{\guillemotright} & \multicolumn{4}{r}{2} & \multicolumn{4}{r}{\guillemotright} & \multicolumn{1}{l}{\quad Size of mini-batches} \\
\multicolumn{4}{r}{480h} & \multicolumn{4}{r}{320h} & \multicolumn{4}{r}{800h} & \multicolumn{4}{r}{1020h} & \multicolumn{1}{l}{\quad CPU time} \\
\multicolumn{4}{r}{60h} & \multicolumn{4}{r}{40h} & \multicolumn{4}{r}{100h} & \multicolumn{4}{r}{130h} & \multicolumn{1}{l}{\quad Resolution time} \\
\bottomrule
\end{tabular}
\begin{tablenotes}
\item [a] Only the median value of the optimal interval is reported to ease the presentation.
\end{tablenotes}
\end{threeparttable}
\caption{\label{tab:control} Open-loop control of circular and square cylinder flows by a small control cylinder of diameter $d=0.1$ - Simulation parameters and convergence data. \red{The interface mesh size yields $\sim 25-40$ grid points in the boundary-layer of the control cylinder, just prior to separation, and the averaging time-span in unsteady flow regimes represents $\sim 15-25$ shedding cycles, depending on the geometry of the main cylinder, the position of the control cylinder and the Reynolds number.}}
\end{center}
\end{table}

\subsection{Optimal drag reduction of a triangular bluff-body using rotating cylinders} \label{sec:olc:sub:pinball}

\begin{figure}[!t]
\setlength{\unitlength}{1cm}
\begin{picture}(20,15.75)
\put(0.1,10.5){\includegraphics[trim=175 92.5 155 340pt,clip,height=5.1cm]{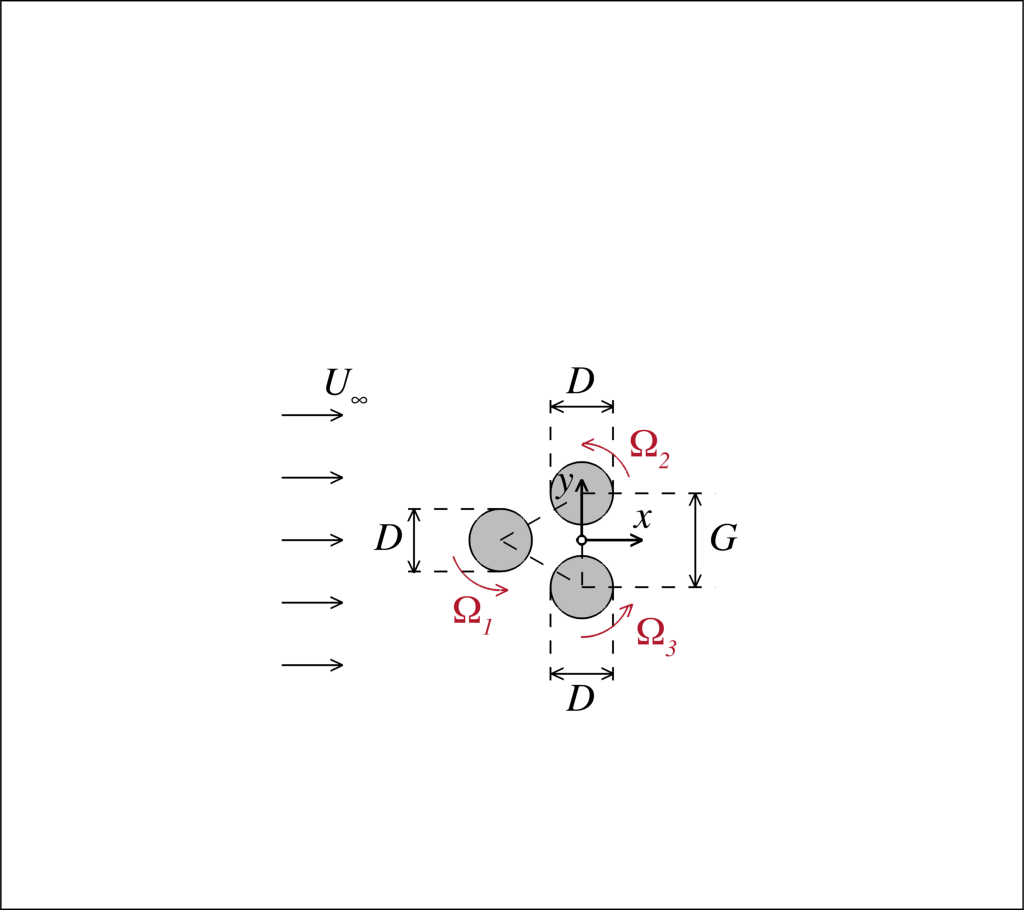}}
\put(7.65,10.5){\includegraphics[trim=175 92.5 155 340pt,clip,height=5.1cm]{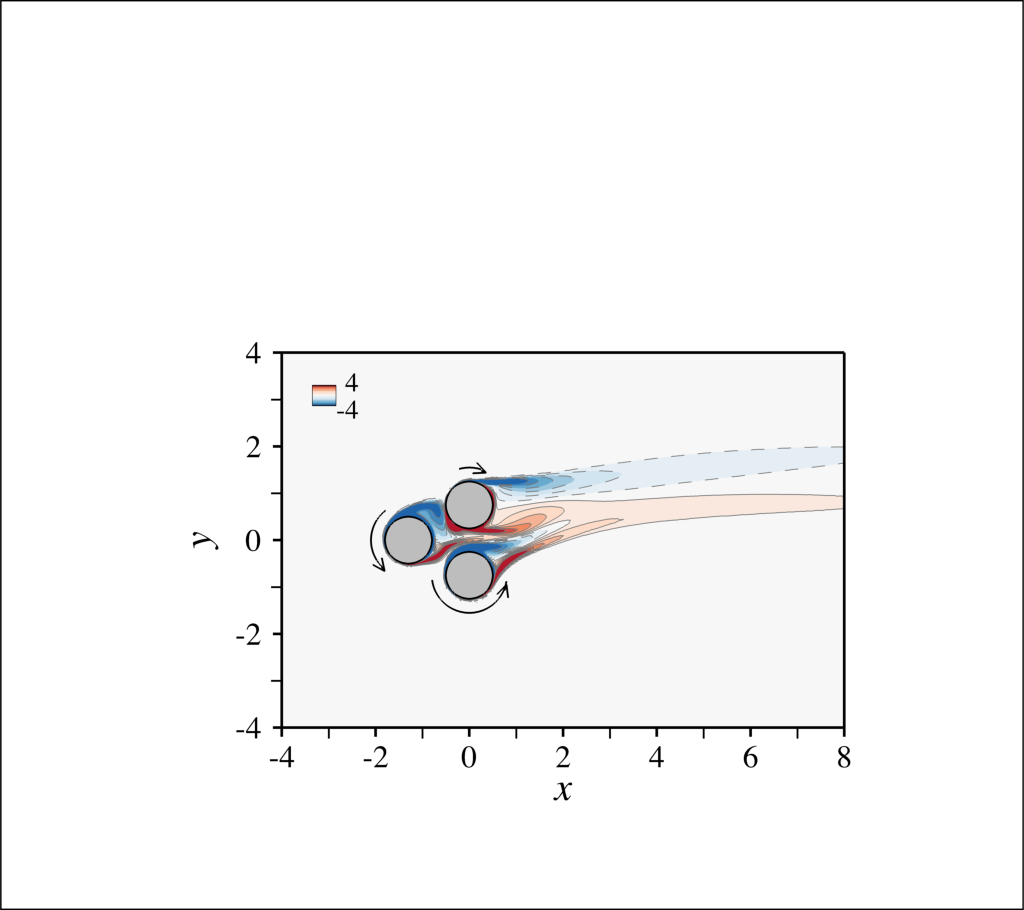}}
\put(0.1,5.25){\includegraphics[trim=175 92.5 155 340pt,clip,height=5.1cm]{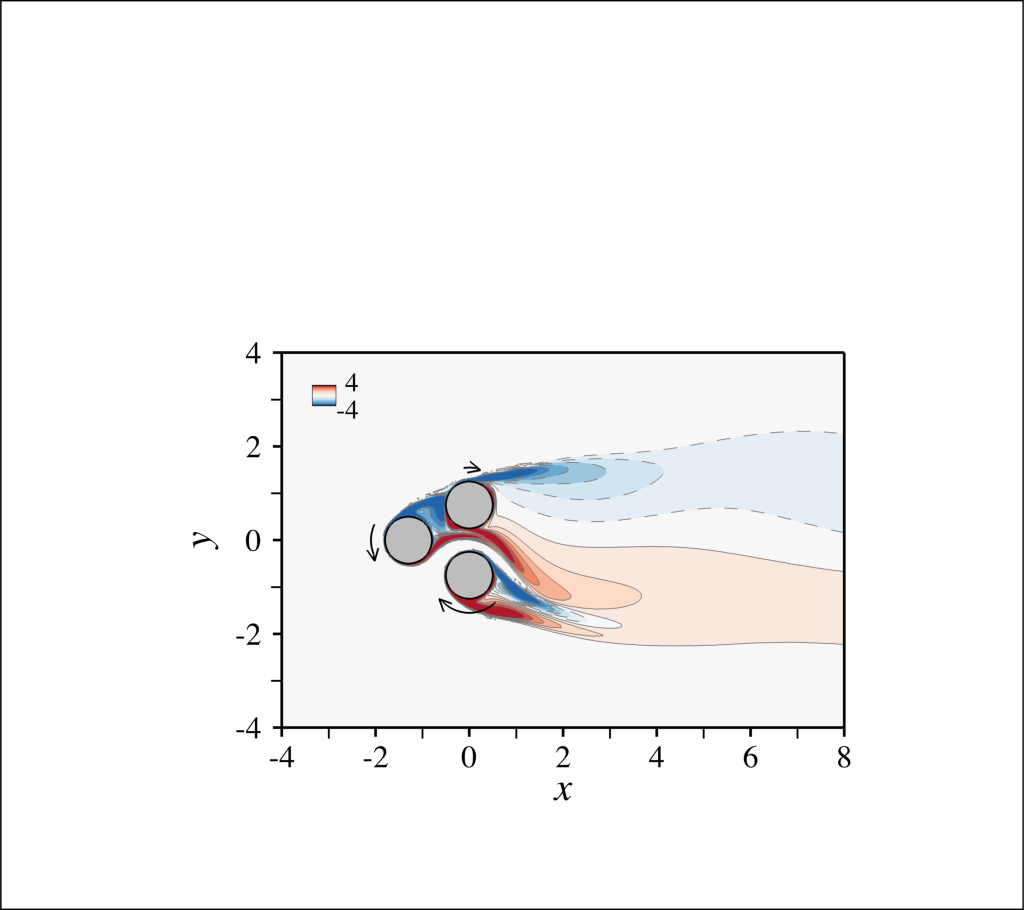}}
\put(7.65,5.25){\includegraphics[trim=175 92.5 155 340pt,clip,height=5.1cm]{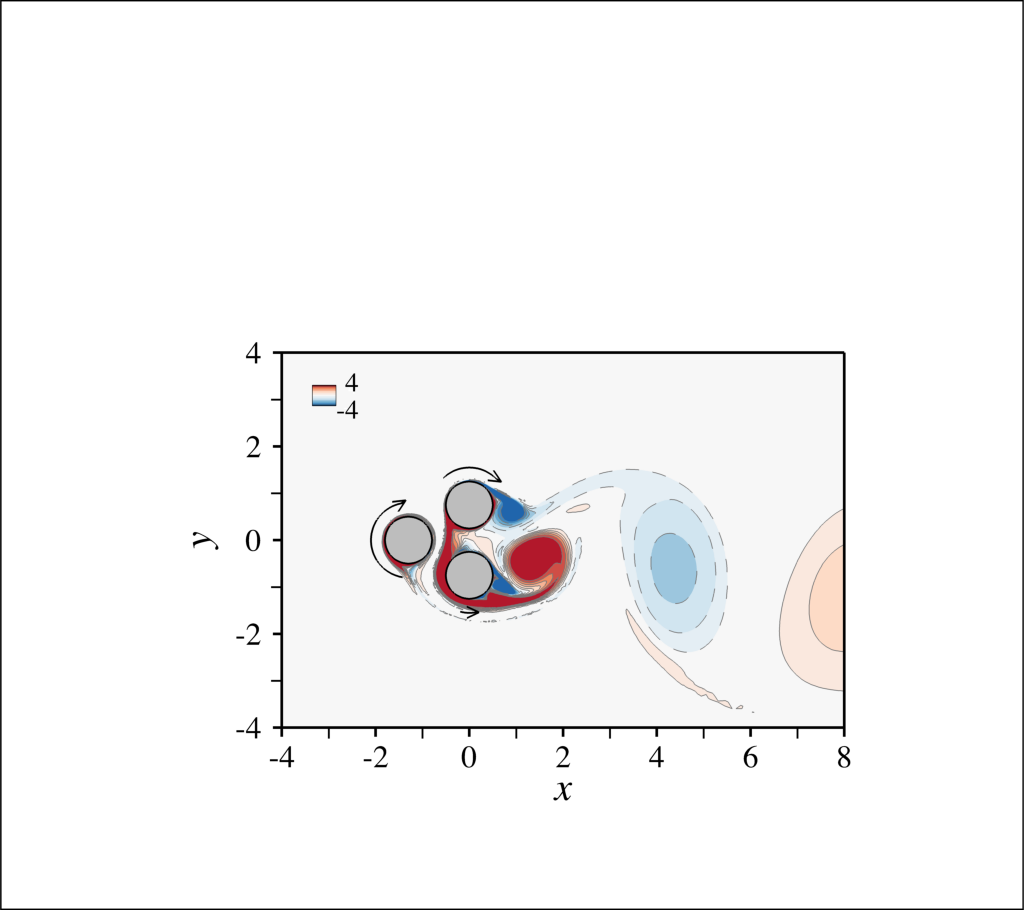}}
\put(0.1,0){\includegraphics[trim=175 92.5 155 340pt,clip,height=5.1cm]{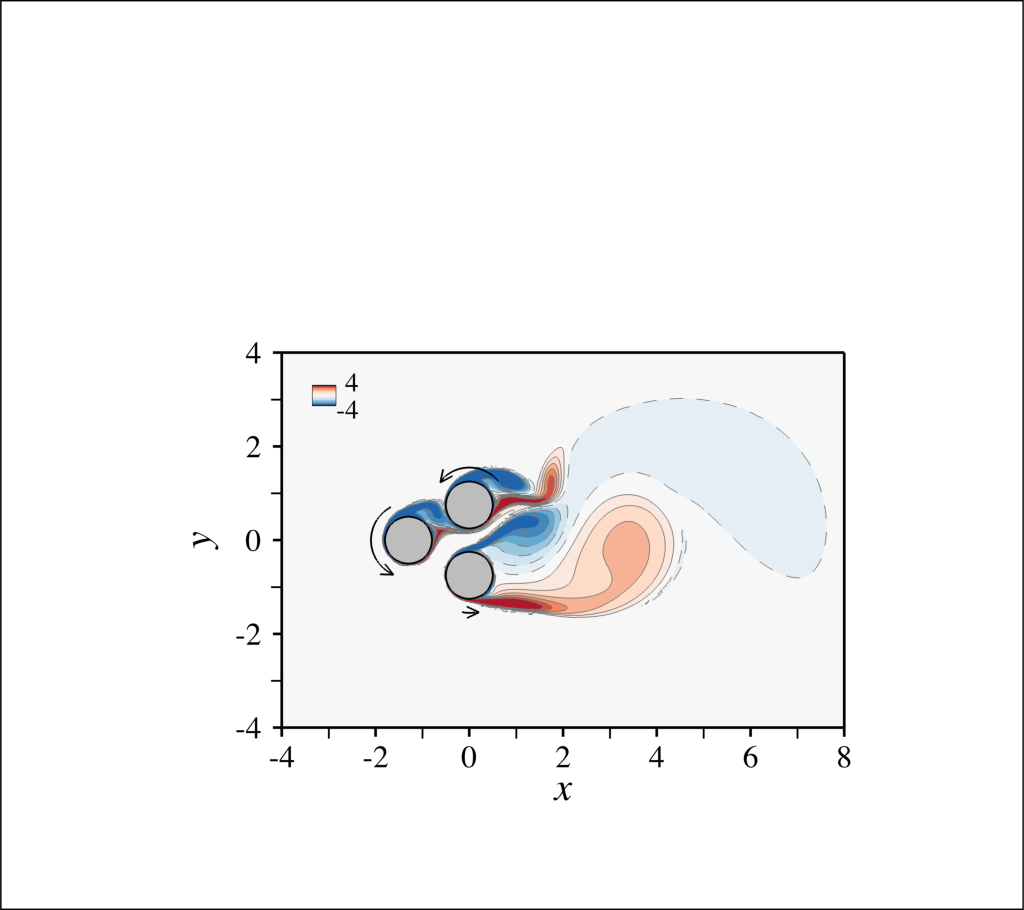}}
\put(7.65,0){\includegraphics[trim=175 92.5 155 340pt,clip,height=5.1cm]{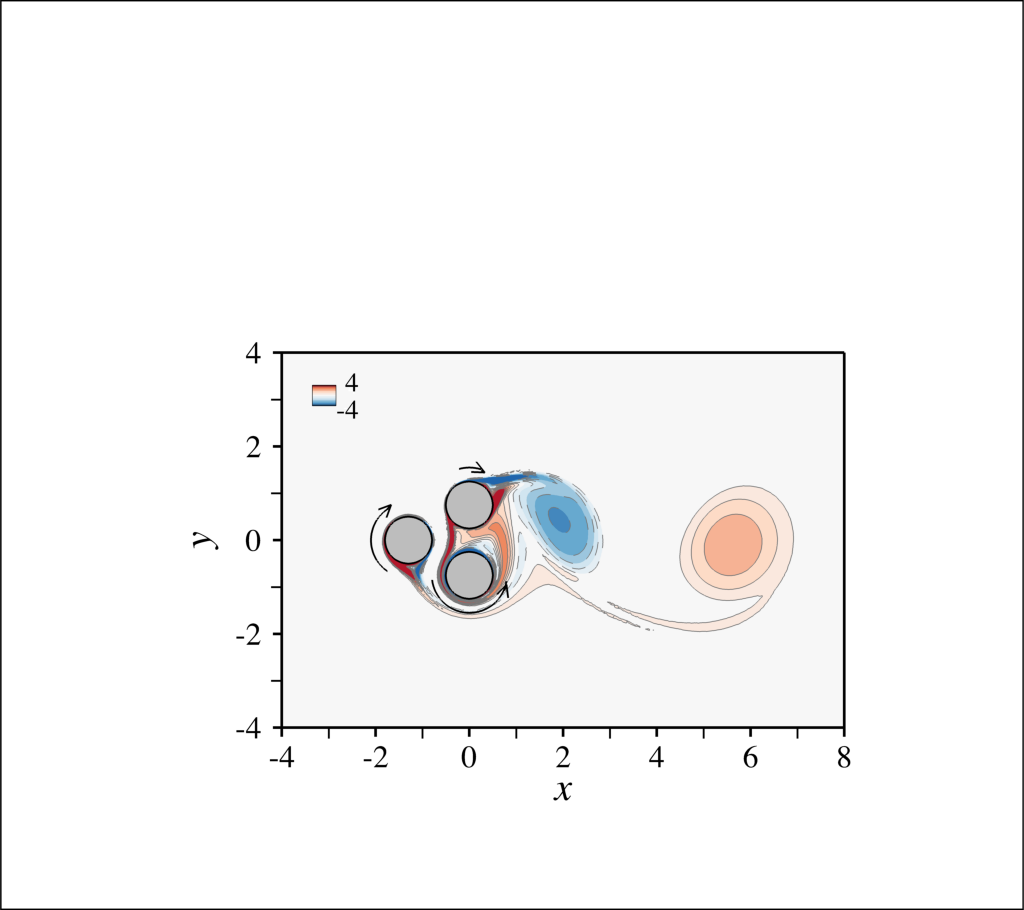}}
\put(0.1,15.75){(a)}
\put(7.65,15.75){(b)}
\put(0.1,10.5){(c)}
\put(7.65,10.5){(d)}
\put(0.1,5.25){(e)}
\put(7.65,5.25){(f)}
\end{picture}
\caption{Open-loop control of a fluid pinball - (a) Schematic diagram of the configuration. (b-f) Iso-contours of the vorticity field computed at $\rey=2200$ for steady angular velocities $(\Omega_1,\Omega_2,\Omega_3)$ of the individual cylinders, namely (b) $(3.09,-1.05,5.00)$, (c) $(1.46,-0.62,-2.75)$ (d) $(-5.00,-2.74,0.81)$, (e) $(3.70,3.00,0.61)$ and (f) $(-3.48,-1.04,5.00)$. The rotation directions are marked by the various arrows whose length is proportional to the angular velocity.}
\label{fig:pinball}
\end{figure}

\red{The second control problem presented in figure~\ref{fig:pinball}(a) is the fluidic pinball~\cite{Raibaudo2017},} an equilateral triangle arrangement of three identical circular cylinders oriented against a uniform stream (i.e., the leftmost triangle vertex points upstream, and the rightmost side is orthogonal to the on-coming flow), controlled open-loop via user-defined angular velocities. The origin of the coordinate system is between the top and bottom cylinders, where we refer to the upstream and downstream cylinders as ``front'', ``top'', and ``bottom'', respectively (also labeled 1, 2 and 3 to ease the notation). The gap spacing $G=1.5$ between cylinders yields a master cross-section of $2.5$. 
A turbulent case at $\rey=U_\infty D/\nu=2200$ is modeled after the negative Spalart--Allmaras uRANS equations, where $D$ is the diameter of either cylinder. The objective is to minimize the mean drag $\overline{D}$ of the three-cylinder system, using the cylinders individual angular velocities $\Omega_{{1-3}}$ as control parameters (with the convention that $\Omega_k< 0$ for clockwise rotation). This is a versatile experiment well suited to challenge the single-step approach, as the requirement to span large ranges of control parameters emulating a variety of \red{steady and unsteady} actuation (e.g., base bleed, suction) under turbulent conditions makes it especially challenging to rely on the adjoint method (as further discussed in section~\ref{sec:disc}), not to mention DNS.

\subsubsection{Steady actuation}\label{sec:olc:sub:pinball_steady} 

\begin{figure}[!t]
\setlength{\unitlength}{1cm}
\begin{picture}(20,11.75)
\put(0.1,5.5){\includegraphics[trim=175 92.5 155 270pt,clip,height=5.85cm]{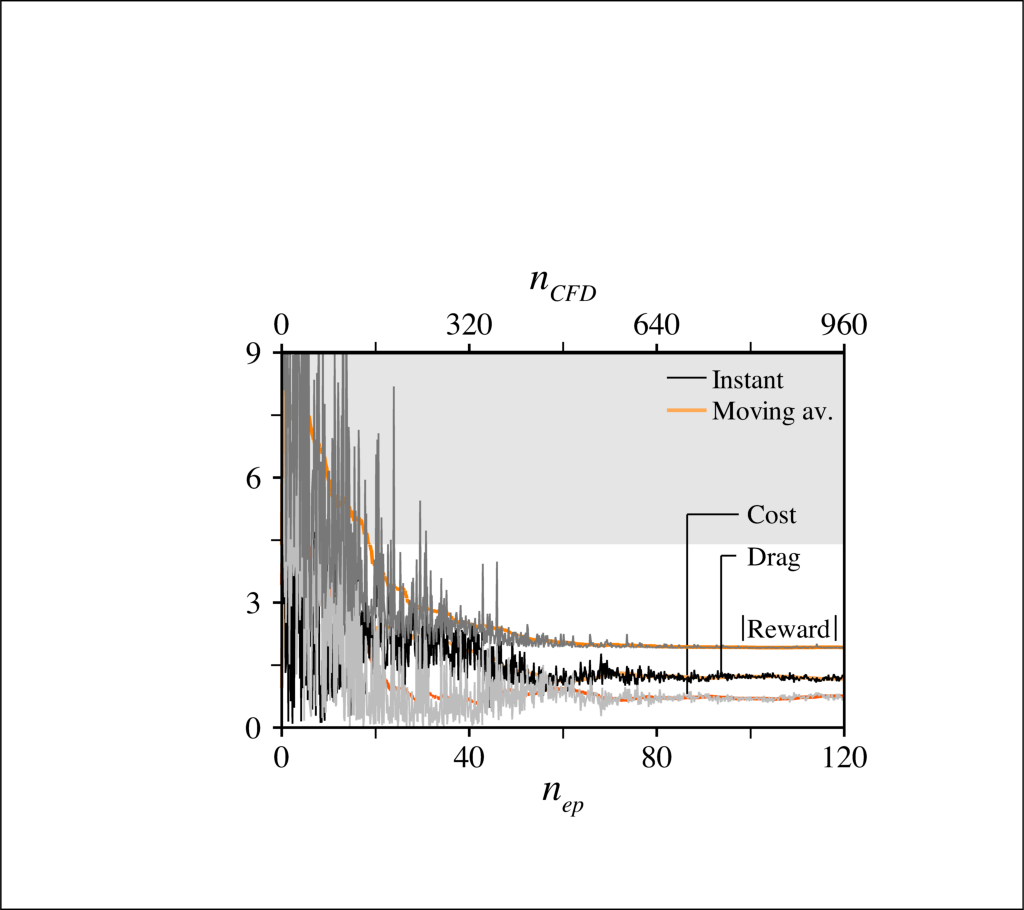}}
\put(7.65,5.5){\includegraphics[trim=175 92.5 155 270pt,clip,height=5.85cm]{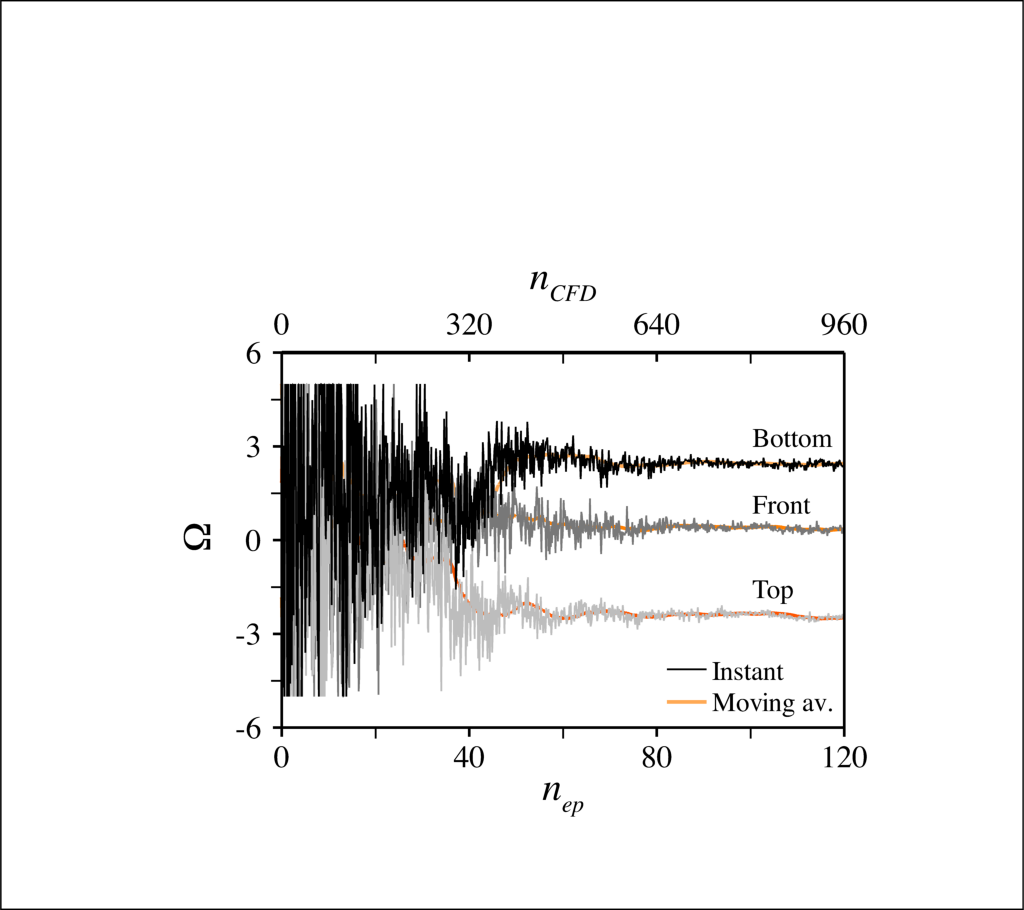}}
\put(0.1,0){\includegraphics[trim=175 92.5 155 340pt,clip,height=5.1cm]{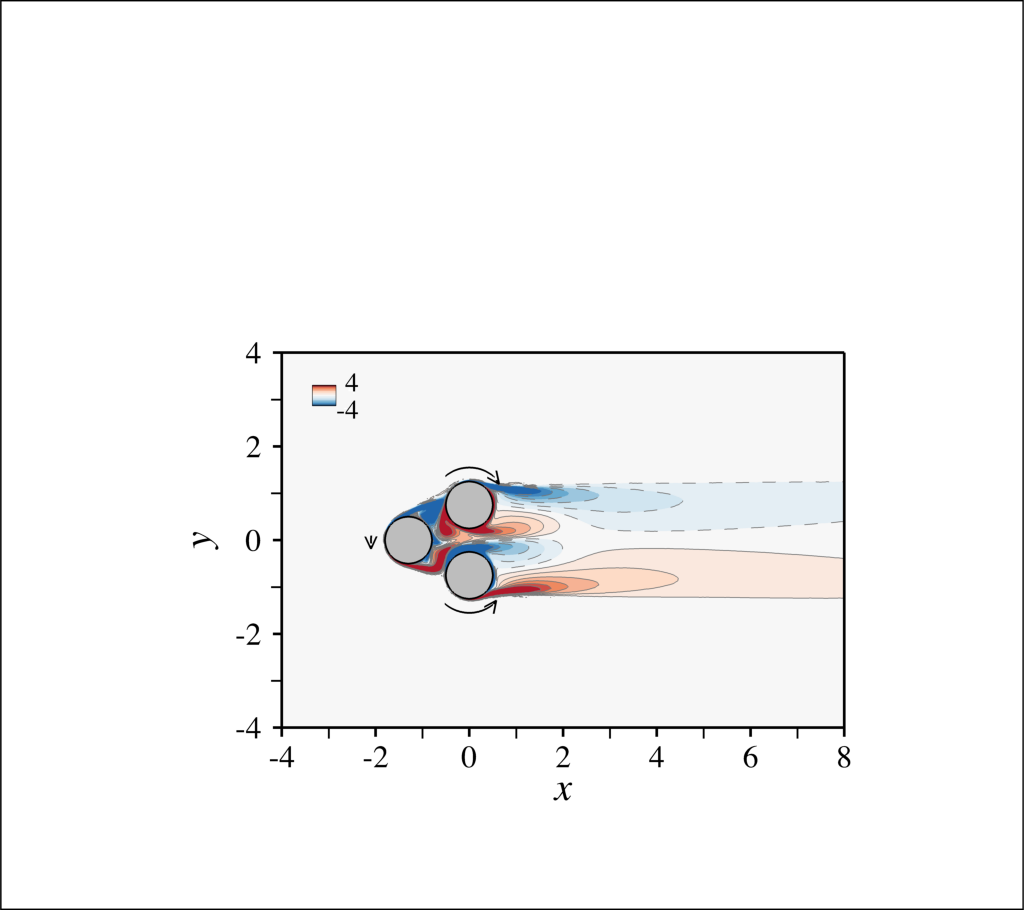}}
\put(7.65,0){\includegraphics[trim=175 92.5 155 340pt,clip,height=5.1cm]{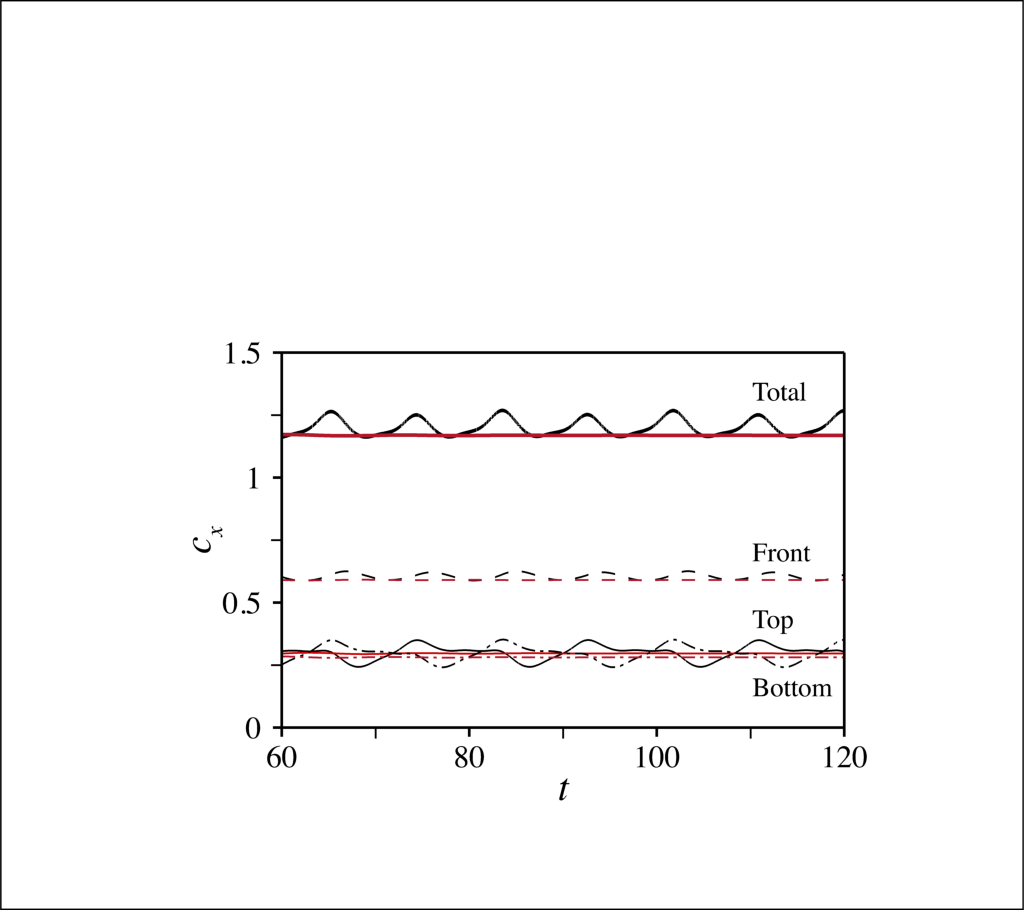}}
\put(4.6,9.6){\includegraphics[trim=440 275 360 445pt,clip,height=0.8cm]{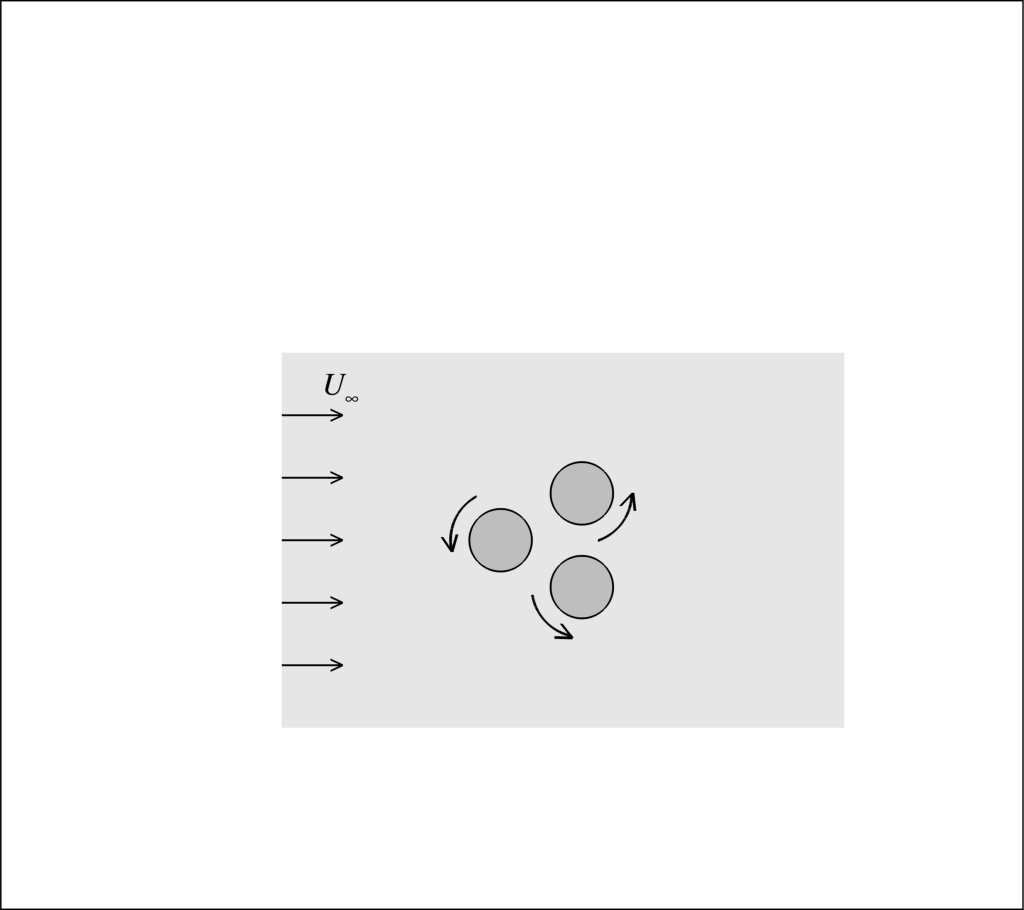}}
\put(12.15,6.575){\includegraphics[trim=440 275 360 445pt,clip,height=0.8cm]{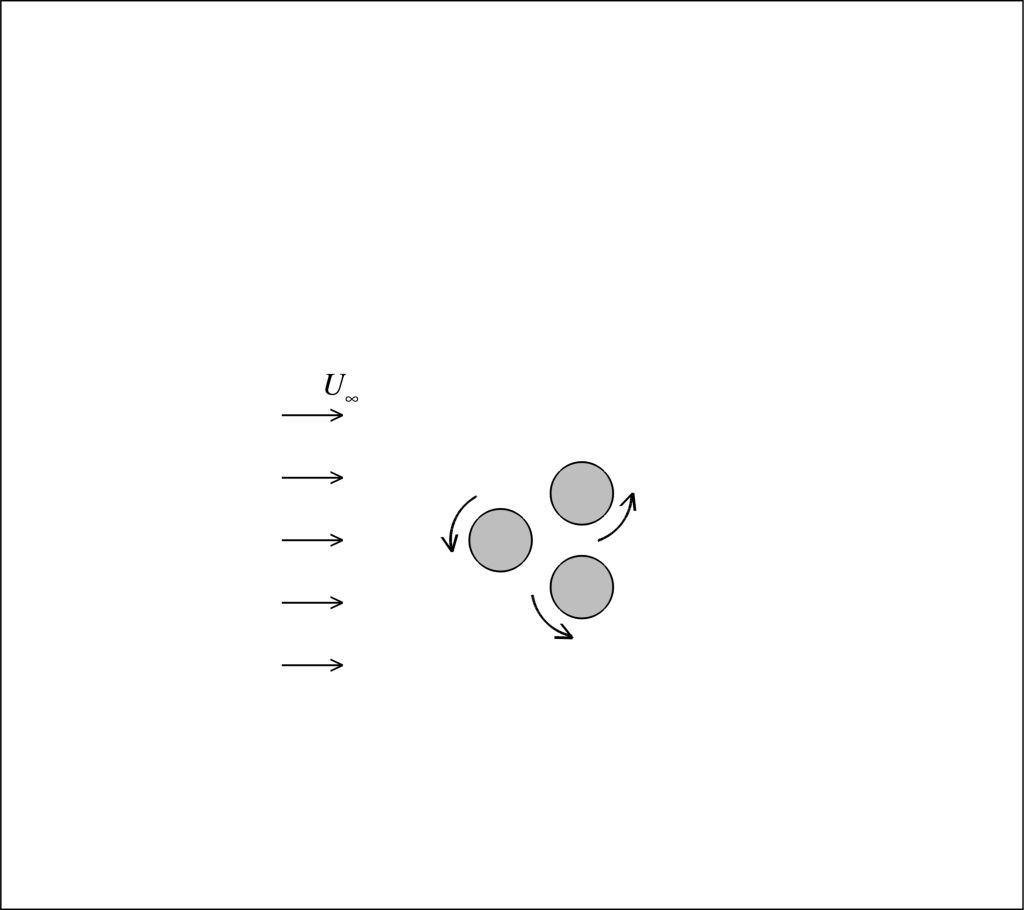}}
\put(0.1,11.75){(a)}
\put(7.65,11.75){(b)}
\put(0.1,5.5){(c)}
\put(7.65,5.5){(d)}
\end{picture}
\caption{Open-loop control of a fluid pinball at $\rey=2200$ - (a) Evolution per episode for the instant (black line) and moving average (over episodes, light orange line) values of the mean drag (over time), together with related cost (light grey/dark orange) and reward (dark grey/orange) information computed under steady actuation for $\beta=0.025$. The uncontrolled drag is at the bottom of the grey shaded area. (b) Same as (a) for the
angular velocities of the front (blak/light orange), top (dark grey/orange) and bottom cylinders (light grey/dark orange). (c) Iso-contours of the vorticity field computed under the optimal velocities $({\Omega_1}^\star,{\Omega_2}^\star,{\Omega_3}^\star)=(0.34,-2.49,2.44)$. (d) Time history of drag computed under the sub-optimal velocities $({\Omega_1},{\Omega_2},{\Omega_3})=(0,-2.47,2.47)$ (black lines), whose cost is identical to that of the optimal (red lines). The thick lines denote the drag of the three-cylinder system. The thin lines pertain to the front (dashed lines), top (solid lines) and bottom cylinders (dash-dotted lines).}
\label{fig:pinball_ppo}
\end{figure} 

\red{First, constant angular velocities are applied to each cylinder to alter} the vorticity flux fed to the wake, as evidenced in figure~\ref{fig:pinball}(b-d) showing instantaneous vorticity fields computed under several control configurations.
\red{Drag is optimized by minimizing the} compound reward function
\bal
r=-\overline{D}-\beta \sum_{k=1}^3|\Omega_k|^3\,,\label{eq:pinball_reward}
\eal
where the leftmost term is the power of the drag force and is thus associated to performance, the rightmost term estimates the power to be supplied to the rotating cylinders and is thus associated to cost, and  $\beta$ is a weighting coefficient set empirically to $\beta=0.025$ (a value found to be large enough for cost considerations to impact the optimization procedure, but not so large as to dominate the reward signal, in which case actuating is meaningless). 
For each PPO-1 learning episode, the network outputs three values $\xi_{{1-3}}$ in $[-1; 1]^3$ mapped into
\bal
\Omega_k=\xi_k\Omega_{\text{max}}\,,
\eal
for the non-dimensional angular velocities to vary in $[-\Omega_{\text{max}};\Omega_{\text{max}}]$ 
with $\Omega_{\text{max}}=5$.
The reward defined in~\myrefeq{eq:pinball_reward} is computed using the simulation parameters documented in table~\ref{tab:pinball}, 
after which the network is updated for 32 epochs using 8 environments and 2 steps mini-batches. 
Note, rotation is actually ramped up over a time-span $[t_{\Omega_i};t_{\Omega_f}]$ to smooth out the transient, using effective rates
\bal
\tilde{\Omega}_k(t)=\frac{\min(\max(t,t_{\Omega_i}),t_{\Omega_f})-t_{\Omega_i}}{t_{\Omega_f}-t_{\Omega_i}}\Omega_k\,,
\eal
forced to zero on $[0,t_{\Omega_i}]$, to $\Omega_k$ on $[t_{\Omega_f};t_f]$, and linearly increasing in between. 

\red{For this case, 120 episodes have been run, which represents 960 simulations, 
each of which lasts $\sim3$h$20$ on 12 cores, hence $\sim 3200$h of total CPU cost (equivalently, $\sim400$h of resolution time).
The moving average reward reaches a plateau after about 80 episodes in figure~\ref{fig:pinball_ppo}(a), where the relevance of the weighing coefficient value $\beta=0.025$ shows through the fact that the performance and cost components of the reward are of the same order of magnitude.
The optimal value of drag $\overline{c_x}^{\,\star}=1.17 \pm 0.01$ computed as the average over the 5 latest episodes represents a tremendous reduction by almost $60\%$ with respect to the uncontrolled value $2.91$.
The associated angular velocities whose evolution is depicted in figure~\ref{fig:pinball_ppo}(b) correspond to a boat tail-like arrangement, i.e., the top cylinder rotates clockwise (${\Omega_2}^\star=-2.49 \pm 0.01$), the bottom cylinder rotates counter-clockwise and almost symmetrically (${\Omega_3}^\star=2.44 \pm 0.01$), and the front cylinder rotates more slowly and also counter-clockwise (${\Omega_1}^\star=0.34 \pm 0.01$). The net rotation is thus in the same direction as the front cylinder, and we show in figure~\ref{fig:pinball_ppo}(c) that the tilting of the shear layers to the centerline alleviates the secondary flow from the gap between the two downstream cylinders, which is found to eventually suppress vortex shedding.
Interestingly, an experimentally implemented machine learning approach using genetic algorithms
yields similar optimal arrangements in~\cite{Raibaudo2020}. 
For two different values of the weighing parameter, the authors therein report optimal angular velocities $(0.68,-2.26,2.56)$ and (1.40,-1.70,2.04) and optimal drag reductions by $78\%$ and  $49\%$, respectively, but it is uneasy to push further the comparison because the latter study uses a different reward function in which drag is approximated from a small, discrete number of sensors distributed in the wake.}

\red{For the purpose of reducing drag, the above asymmetrical boat tailing actuation turns to be more efficient than its pure, symmetrical counterpart emulated by $({\Omega_1},{\Omega_2},{\Omega_3})=(0,-|\Omega|,|\Omega|)$.\footnote{At least if $\beta$ is large enough for cost to matter in the optimization procedure, otherwise the algorithm has been found to converge to the
symmetrical boat tailing configuration $(0,-\Omega_{\text{max}},\Omega_{\text{max}})$, and the reverse flow is completely suppressed (not shown here).}
This is illustrated in figure~\ref{fig:pinball_ppo}(d) comparing the optimal drag to its symmetrical value computed with $|\Omega|=2.47$ (to maintain the same cost efficiency, the associated drag reduction being by $\sim 58\%$). Pure boat tailing is insufficient to inhibit vortex shedding,
as the symmetrical drag of all three individual cylinders is seen to exhibit small but finite-amplitude oscillations. Moreover, the drag of the downstream cylinders turns to be roughly identical on average. This suggests that the edge of asymmetrical over symmetrical boat tailing lies in its ability to reduce the drag of the front cylinder, an effect similar to that of suppressing vortex development and reducing drag by creating circulation around a single rotating bluff body~\cite{Kang1999}. Asymmetrical boat tailing is also more efficient than base bleed, another method widely used to reduce drag by blowing fluid directly into the wake, and that can be emulated by $({\Omega_1},{\Omega_2},{\Omega_3})=(0,|\Omega|,-|\Omega|)$ for the reverse rotation of the downstream cylinders to conversely enhance the gap flow in between them (not shown here).}

\subsubsection{Periodic actuation}\label{sec:olc:sub:pinball_periodic}  

\begin{figure}[!t]
\setlength{\unitlength}{1cm}
\begin{picture}(20,5.5)
\put(0.1,0){\includegraphics[trim=175 92.5 155 340pt,clip,height=5.1cm]{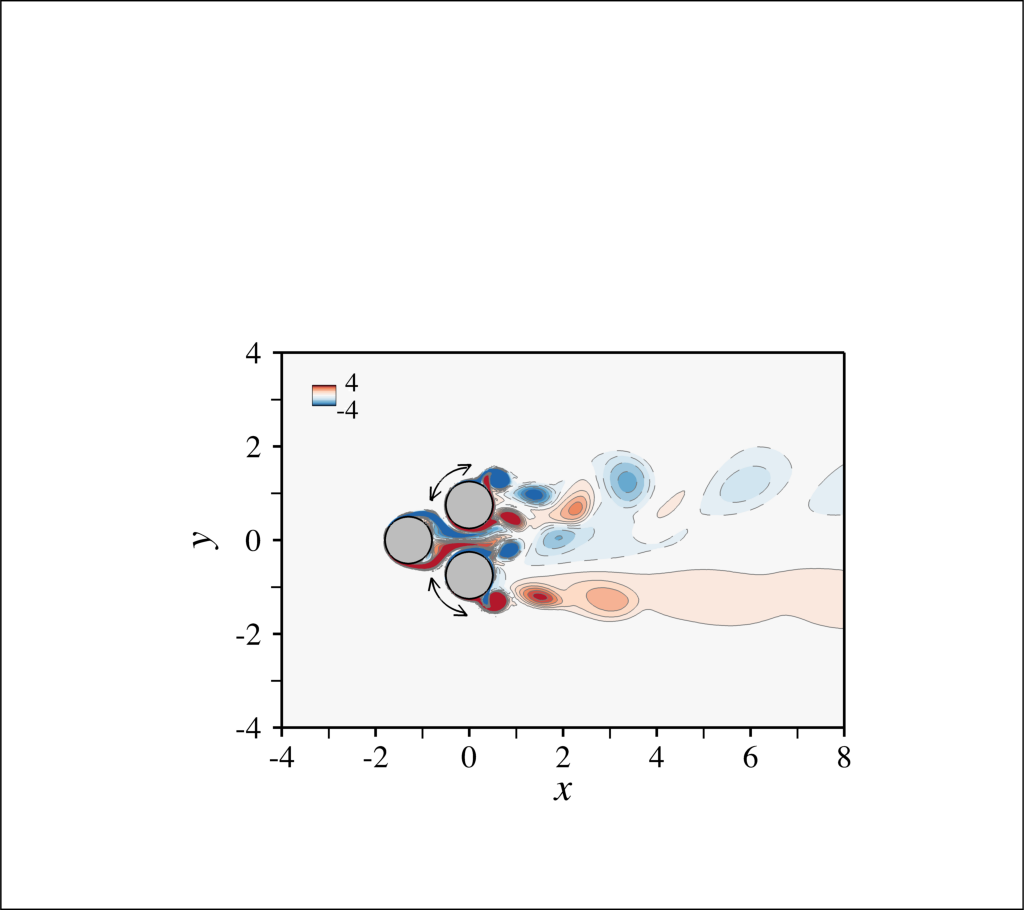}}
\put(7.65,0){\includegraphics[trim=175 92.5 155 340pt,clip,height=5.1cm]{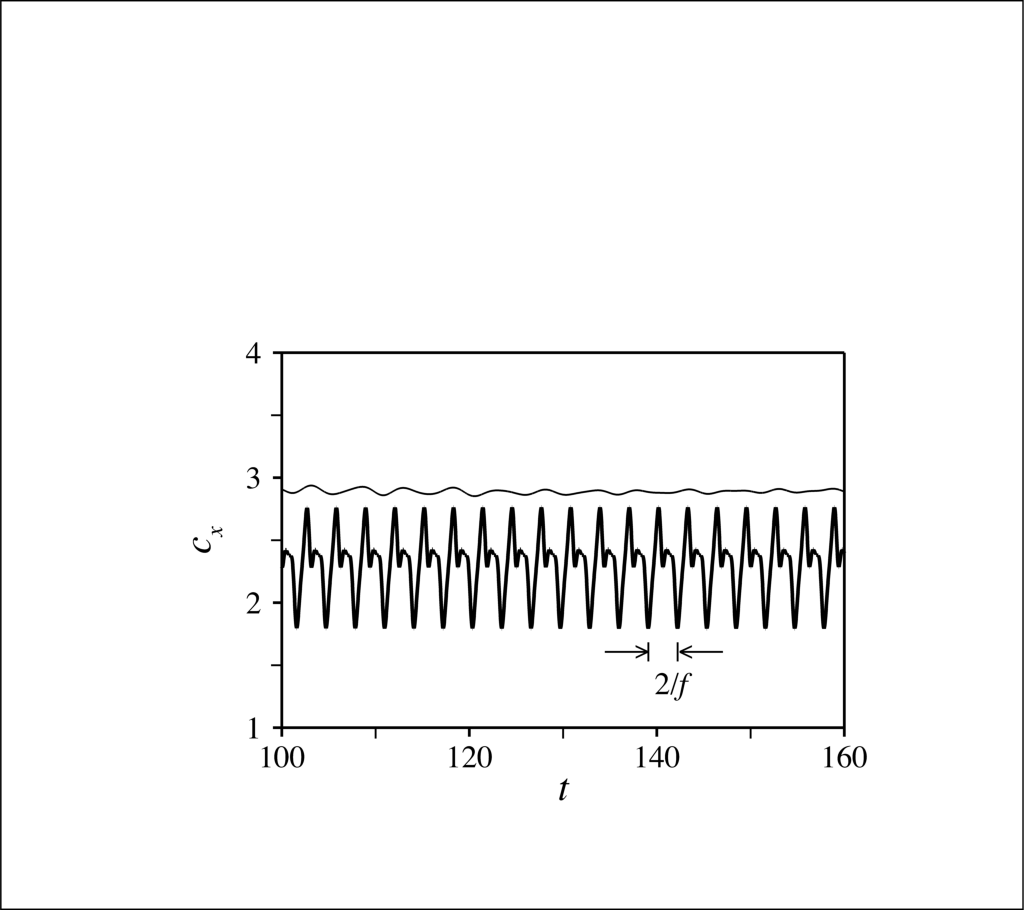}}
\put(0.1,5.5){(a)}
\put(7.65,5.5){(b)}
\end{picture}
\caption{Open-loop control of a fluid pinball at $\rey=2200$ - (a) Iso-contours of the vorticity field and (b)
time history of drag computed under periodic actuation~\myrefeq{eq:periodic} with angular velocity $\Omega=2.47$ and frequency $f=4f_0$. The thick and fine lines in (b) denote the controlled and uncontrolled values, respectively.}
\label{fig:pinball_ppoperiodic}
\end{figure}

Periodic actuation at frequency $f$ has also been considered using a simplified configuration 
\bal
\Omega_1=0\,,\qquad \Omega_2=-\Omega_3=\Omega\sin(2\pi f t)\,,\label{eq:periodic}
\eal
whose front cylinder is fixed, and whose downstream cylinders are periodically and symmetrically driven with maximum angular velocity  $\Omega$.
Such a control oscillates between symmetrical boat tailing (found to be nearly-optimal under steady actuation) and base-bleed, and we assess the extent to which an additional degree of freedom (the oscillation frequency) creates room to improve the performance. The optimization relies on the compound reward
\bal
r=-\overline{D}-2\beta|\Omega|^3\,,\label{eq:pinball_reward_freq}
\eal
computed using the same weighing parameter $\beta=0.025$ as before. 
For each PPO-1 learning episode, the network outputs two values $\xi_{{1, 2}}$ in $[-1; 1]^2$ mapped into
\bal
\Omega=\frac{1+\xi_1}{2}\Omega_{\text{max}}\,,\qquad
\frac{f}{f_0}=\frac{1-\xi_2}{2}\lambda_{\text{min}}+\frac{1+\xi_2}{2}\lambda_{\text{max}}\,,
\eal
where $f_0=0.16$ is the dominant frequency of vortex shedding computed in the absence of control. The angular velocity therefore varies in $[0;\Omega_{\text{max}}]$ 
with $\Omega_{\text{max}}=5$ (the case $\Omega<0$ is covered by periodicity) and the frequency ratio varies in 
$[\lambda_{\text{min}};\lambda_{\text{max}}]$ with $\lambda_{\text{min}}=0.5$ and $\lambda_{\text{max}}=4$. This is a compromise between size of the parameter space and cost control, as investigating
smaller frequencies would require to increase the averaging time-span, and resolving accurately larger frequencies would require to decrease the time-step.
We shall not go into the details of the obtained results, because the frequency ratio ends up oscillating randomly in $[\lambda_{\text{min}};\lambda_{\text{max}}]$, while the angular velocity converges to ${\Omega}^\star=0$. It is definitively possible to reduce drag under the considered periodic actuation, as we show for instance in figure~\ref{fig:pinball_ppoperiodic} that a velocity $\Omega=2.47$ (identical to that \red{used previously to compare asymmetrical and} symmetrical boat tailing) and a frequency ratio $\lambda=4$ reduce drag by $20\%$, but the cost of doing so is too large, as the associated reward actually increases by $5\%$
\red{(note the period doubling bifurcation phenomenon in figure~\ref{fig:pinball_ppoperiodic}(b): drag is found to exhibit sub-harmonic oscillations at half the forcing frequency, which is a classical dynamical responses of harmonically forced nonlinear oscillators). These are only preliminary results intended to compare the efficiency of steady and periodic strategies using identical reward functions. We therefore defer to future work the computation of non-trivial periodic optimal distributions, for which it may be necessary to modify the reward function and/or to reduce the cost (by adequately decreasing the weighing parameter).} 

\begin{table}[!t]
\begin{center}
\begin{tabular}{cccccccccccc}
\toprule
\multicolumn{1}{r}{\multirow{3}{*}{\makecell[r]{\raisebox{-\totalheight}{\includegraphics[trim=400 260 380 460pt,clip,height=1cm]{fig14_thumbb.png}}}}}
& \multicolumn{1}{r}{\makecell[r]{$r$}} & \multicolumn{1}{r}{\makecell[r]{$\overline{c_x}$}}  & \multicolumn{1}{r}{\makecell[r]{$\Omega_1$}} & \multicolumn{1}{r}{\makecell[r]{$\Omega_2$}} & \multicolumn{1}{r}{\makecell[r]{$\Omega_3$}}  
& \multicolumn{1}{r}{\multirow{3}{*}{\makecell[r]{\raisebox{-\totalheight}{\includegraphics[trim=400 260 380 460pt,clip,height=1cm]{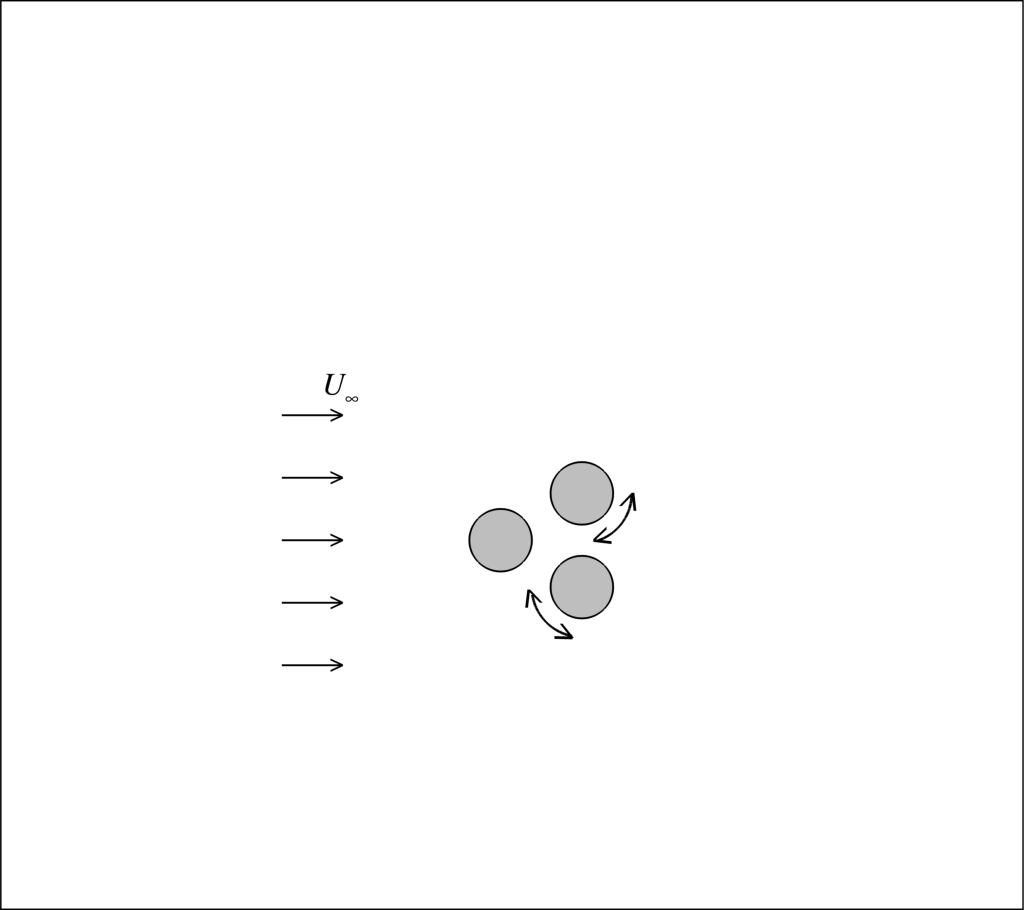}}}}} & \multicolumn{1}{r}{\makecell[r]{$r$}} & \multicolumn{1}{r}{\makecell[r]{$\overline{c_x}$}} & \multicolumn{1}{p{0.6cm}}{\makecell[r]{$\Omega\,$}} & \multicolumn{1}{r}{\makecell[r]{$f$}}\\
\cmidrule(lr){2-2}\cmidrule(lr){3-3}\cmidrule(lr){4-4}\cmidrule(lr){5-5}\cmidrule(lr){6-6}\cmidrule(lr){8-8}\cmidrule(lr){9-9}\cmidrule(lr){10-10}\cmidrule(lr){11-11}
& \multicolumn{1}{r}{$-1.93$} & \multicolumn{1}{r}{$1.17$} & \multicolumn{1}{r}{$0.34$} & \multicolumn{1}{r}{$-2.47$} & \multicolumn{1}{r}{$2.44$} & & $-2.91$ & 2.91 & \multicolumn{1}{p{0.6cm}}{\makecell[r]{$0\,$}} & N/D & \multicolumn{1}{l}{\qquad Optimal}\\
\\
\cmidrule(lr){1-12}
\multicolumn{12}{r}{CFD} \\
\cmidrule(lr){1-12}
\multicolumn{6}{r}{2200} & & \multicolumn{4}{r}{\guillemotright} & \multicolumn{1}{l}{\qquad Reynolds number} \\
\multicolumn{6}{r}{Steady} & & \multicolumn{4}{r}{Periodic} & \multicolumn{1}{l}{\qquad Actuation} \\
\multicolumn{6}{r}{0.05} & &  \multicolumn{4}{r}{0.025} & \multicolumn{1}{l}{\qquad Time-step} \\
\multicolumn{6}{r}{$[5;10]$} & & \multicolumn{4}{r}{\guillemotright} & \multicolumn{1}{l}{\qquad Rotation ramp-up time span} \\
\multicolumn{6}{r}{$[300;400]$} & & \multicolumn{4}{r}{\guillemotright} & \multicolumn{1}{l}{\qquad Averaging time span} \\
\multicolumn{6}{r}{$[-6;20]\times[-6;6]$} & & \multicolumn{4}{r}{\guillemotright} & \multicolumn{1}{l}{\qquad Mesh dimensions} \\
\multicolumn{6}{r}{$110000$} & & \multicolumn{4}{r}{\guillemotright} & \multicolumn{1}{l}{\qquad Nb. mesh elements} \\
\multicolumn{6}{r}{$0.001$} & & \multicolumn{4}{r}{\guillemotright} & \multicolumn{1}{l}{\qquad Interface $\perp$ mesh size} \\
\multicolumn{6}{r}{$12$} & & \multicolumn{4}{r}{\guillemotright} & \multicolumn{1}{l}{\qquad Nb. Cores} \\
\cmidrule(lr){1-12}
\multicolumn{12}{r}{PPO-1} \\
\cmidrule(lr){1-12}
\multicolumn{6}{r}{120} &&  \multicolumn{4}{r}{40} & \multicolumn{1}{l}{\qquad Nb. episodes} \\
\multicolumn{6}{r}{8} & & \multicolumn{4}{r}{\guillemotright} & \multicolumn{1}{l}{\qquad Nb. environments} \\
\multicolumn{6}{r}{32} & & \multicolumn{4}{r}{\guillemotright} & \multicolumn{1}{l}{\qquad Nb. epochs} \\
\multicolumn{6}{r}{2} & & \multicolumn{4}{r}{\guillemotright} & \multicolumn{1}{l}{\qquad Size of mini-baches} \\
\multicolumn{6}{r}{3200h} & & \multicolumn{4}{r}{2100h} & \multicolumn{1}{l}{\qquad CPU time} \\
\multicolumn{6}{r}{400h} & & \multicolumn{4}{r}{260h} & \multicolumn{1}{l}{\qquad Resolution time} \\
\bottomrule
\end{tabular}
\caption{\label{tab:pinball} Simulation parameters and convergence data for open-loop control of a fluid pinball at $\rey=2200$. The interface mesh size yields $\sim 20$ grid points in the boundary-layer of the non-rotating front, top and bottom cylinder, just prior to separation, and the averaging time-span represents $\sim 15-20$ shedding cycles, depending on the angular velocities. For the periodic case, the time step yields $\sim 60$ data points over the smallest actuation period.}
\end{center}
\end{table}

\section{Discussion} \label{sec:disc}

This section is intended to provide insight into the efficiency of the single-step PPO algorithm compared to that of other well-established methods. We skip voluntarily DNS, as systematical optimization procedures are useless if a problem is simple enough that a small number of numerical simulations suffices to find the optimal. This is true of the optimization cases documented in section~\ref{sec:opt}, although the results remain valuable to assess accuracy and highlight the limit of applying conservative policy updates to optimize sharp reward functions (that are common occurrence in low to moderate-Reynolds-number-fluid mechanical systems sustaining linear instabilities).

\subsection{Adjoint methods} \label{sec:disc:sub:adj}

\red{We begin with the adjoint method used in section~\ref{sec:olc:sub:control}
for systematic validation purposes. 
As explained in appendix~\ref{sec:appadj}, this is an approach intended to compute the drag of a control-induced disturbance modeled after the linearized governing equations forced by small-amplitude momentum source $\delta\ff$ and wall velocity $\delta\uu_w$, without ever computing the disturbance itself. The main assumptions and limitations at various levels of sophistication are reviewed in the appendix, so the line of though is to describe only the 
specifics of the control problems considered herein.}
The general picture is that the baseline adjoint method is accurate and fairly efficient in terms of CPU cost, but demanding in terms of storage and increasingly difficult to apply rigorously when turbulence sets in (this is discussed in appendix~\ref{sec:appadj:sub:base}). \red{On the other hand, the frozen Reynolds stresses approximation has marginal CPU and storage costs, it carries over to any turbulence modeling under the
so-called frozen viscosity assumption, but accuracy must be assessed on a case-by-case basis  (see appendix~\ref{sec:appadj:sub:mean}).}

\subsubsection{Open-loop control by a small control cylinder} \label{sec:disc:sub:adjcontrol}

Open-loop control by a small control cylinder is a favorable case in the sense that only the center position of the control cylinder (not its shape, nor its size) is optimized, hence 
the adjoint problem needs be solved only once. Nonetheless, it comes with a substantial modeling component, as the source term $\delta\ff$ used in the adjoint calculations must adequately represent the effect of a true control cylinder. We use here the pointwise reacting force proposed in~\cite{meli14}, equal and opposite to the force felt by a control cylinder of same diameter in a uniform flow at the local, mean velocity. The latter is carefully crafted to reference data, but there are inherent approximations associated with overlooking the lift component of the force induced by the local velocity gradient (since the control cylinder, albeit small, has finite size) and inertia (for the model force at each time instant to be the force that would act if the upstream flow at the same instant was a steady one). This can hurt accuracy and undermine the results in flow regions where the control cylinder drag is close to balancing the decrease in the drag of the main cylinder, all the more so in turbulent regimes where additional simplifications are needed to allow implementing the adjoint method itself (e.g., frozen eddy viscosity and/or Reynolds stresses).

In terms of pure performance, the baseline adjoint method is beyond compare for the laminar, steady case at $\rey=40$, because it merely requires solving a couple of steady solutions (one nonlinear, one linear), and PPO-1 would need converge in less than two episodes to approach that cost. 
\red{Regarding the laminar, time-dependent case at $\rey=100$, the results reported herein rely on a naive implementation of the adjoint method: all time steps of the uncontrolled solution are written to disk, the adjoint equations are solved over the same time interval and with the same time step, and
meaningful time averages of the adjoint-based integrands are computing after discarding the early and late time steps (corresponding to transients of the uncontrolled and adjoint solutions). In practice, this takes 45 Gb of storage. The cost of tackling similarly a three-dimensional (3-D) case with 40 points distributed in the span-wise direction would thus be about $2$ Tb (as estimated by simple cross-multiplication), which is close to intractable without sophisticated integration, interpolation and/or checkpointing schemes. Meanwhile, the storage cost of PPO-1 is barely a few hundred Mb overall, and is expected to jump to a few ten Gb in 3-D without any additional development.}
As for CPU cost, the adjoint method amounts to roughly 7-8 episodes, which is about thrice as less as the number of episodes needed to achieve convergence with PPO-1 (this is an estimation for two numerical simulations oversized by the repeated IO calls, although an exact comparison is difficult because our DRL and adjoint results have been obtained using a different finite element codes on different hardware resources).
Finally, for the turbulent cases at $\rey=3900$ and $\rey=22000$, the cost of the adjoint method is again marginal, as we relied on the frozen Reynolds stresses formulation for which it suffices to compute a  nonlinear uncontrolled mean flow and a linear steady adjoint solution. PPO-1 would need to converge in one single episode to match the cost, but we believe the case at $\rey=3900$ to provide clear evidence that the simplifying assumptions can make it intricate to \red{compare both qualitatively and quantitatively}.

\subsubsection{Open-loop control of a fluidic pinball} \label{sec:disc:sub:adjpinball}

The adjoint modeling of the fluidic pinball is straightforward, since the wall velocity $\delta\uu_w$ is simply the cylinder linear velocity. 
\red{The challenge for this case rather lies in the large value of the optimal angular velocities (found to induce velocities close to the ambient velocity in the vicinity of the downstream cylinders), that suffice to invalidate the linearity assumption inherent to the adjoint method.}
On paper, this problem can still be tackled with a nonlinear steepest descent algorithm recursively solving an adjoint problem and modifying the control parameters in the direction of the negative gradient. 
While it usually takes about ten iterations for fluid mechanical systems to converge (provided relevant 
update strategy and descent step are used), we did not attempt to do so, as it would magnify the limitations of the adjoint method underlined in the appendix. Namely, the storage cost would increase (even a simple conjugate gradient algorithm would require availability of multiple time histories of adjoint solutions) and convergence could be weakened or even sapped if the simplifications made in turbulent regimes yield inaccurate gradient evaluations.

\subsection{Evolution strategies} \label{sec:disc:sub:other}

\red{Evolution strategies (ES) are another popular family of division of population-based algorithms performing black-box optimization in continuous search spaces without computing directly the gradient of the target function. ES imitate principles of organic evolution processes as rules for optimum seeking procedures, using repeated interplay of variation (via recombination and mutation) and selection in populations of candidate solutions. 
They rely on a stochastic description of the variables to optimize, i.e., they consider probability density functions instead of deterministic variables. At each generation (or iteration) new candidate solutions are sampled isotropically by variation of the current parental individuals according to a multivariate normal distribution. After applying recombination and mutation transformations (respectively  amounting to selecting a new mean for the distribution, and to adding a random perturbation with zero mean), the individuals with the highest cost function are then selected to become the parents in the next generation. Improved variants include the covariance matrix adaptation evolution strategy (CMA-ES), that 
also updates its  full covariance matrix to accelerate convergence toward the optimum (which amounts to learning a second-order model of the underlying objective function).}

\red{As has been said for introductory purposes, it lies out of the scope of this paper to provide exhaustive performance comparison data against state-of-the art evolution algorithms. 
The efforts for developing single-step PPO remain at an early stage, so we do not expect the method to be able to compete right away.
Nonetheless, we do not expect it to be utterly outmatched either, as genetic algorithms\footnote{Another class of evolutionary algorithms with slightly different implementation details. Namely, 
most parameters in genetic algorithms (GA) are exogenous, i.e., set by the practitioner, while ES features endogenous parameters associated with individuals, that evolve together with them. Also, only the fittest individuals are selected to become parents in GA, while parents are selected randomly in ES and the fittest offsprings are selected and inserted in the next generation.} have been shown capable to learn optimal open- and closed-loop control strategies within a few hundreds to a few thousands test runs (see~\cite{Deng2018} and the references therein), and it takes a few hundred (resp. less than one thousand) simulations for single-step PPO to learn the
optimal open-loop strategy for control by a small cylinder (resp. for control of the fluidic pinball). 
In present form, the method can be thought as an evolutionary-like algorithm with simpler heuristics (i.e., without an evolutionary update strategy, since the optimal model parameters are learnt via gradient ascent). Its performance should thus be comparable to that of standard ES methods with isotropic covariance matrix, meaning that further characterization and fine-tuning, as well as pre-trained deep learning models (as is done in transfer learning) are likely required to outperform more advanced methods. }


\section{Conclusion}

Open-loop control of laminar and turbulent flow past bluff bodies is achieved here training a
fully connected network with a novel single-step PPO deep reinforcement algorithm, in which it gets only one attempt per learning episode at finding the optimal. The numerical reward fed to the network is computed with a finite elements CFD environment solving stabilized weak forms of the governing equations (Navier--Stokes, otherwise uRANS with negative Spalart--Allmaras as turbulence model) with a combination of variational multiscale approach, immersed volume method and anisotropic mesh adaptation.

Convergence and accuracy are assessed from two optimization cases (maximizing the mean lift of a NACA 0012 airfoil and the fluctuating lift of two side-by-side circular cylinders, both in laminar regimes). Those are simple enough to allow comparison to in-house DNS data, yet they stress that the occurrence of instability yields sharp reward functions for which the conservative policy updates specific to PPO can trap the optimization process into local optima. The method is also applied to two open-loop control problems whose parameter spaces are large enough to dismiss DNS. Single-step PPO is found to successfully reduce the drag of laminar and turbulent cylinder flows by mapping the best positions for placement of a small control cylinder in good agreement with reference data obtained by the adjoint method. The achieved reduction ranges from $~2\%$ using a circular geometry of the main cylinder at $\rey=40$, up to $30\%$ using a square geometry at $\rey=22000$. Second, the method proves fruitful to reduce the drag of the fluidic pinball, an arrangement of three identical, rotating circular cylinders immersed in a turbulent stream. An optimal reduction by almost $60\%$ (consistent with that recently obtained using genetic algorithms) is reported using a boat tailing actuation made up of a slowly rotating front cylinder and two downstream cylinders rotating in opposite directions so as to reduce the gap flow in between them. 
For both cases, convergence is reached after a few ten episodes, which represents a few hundreds CFD runs. \red{Exhaustive computational efficiency data are reported with the hope to foster future comparisons, but it is worth emphasizing that we did not seek to optimize said efficiency, neither by optimizing the hyper parameters, nor by using pre-trained deep learning models.}

Fluid dynamicists have just begun to gauge the relevance of deep reinforcement learning techniques to assist the design of optimal flow control strategies. This research weighs in on this issue and shows that the proposed single-step PPO holds a high potential as a reliable, go-to black-box optimizer for complex CFD problems.
The one advantages here are scope and applicability, as the storage cost of an episode is simply that of a CFD run (times the number of environments), and there is no prerequisite beyond the ability to compute accurate numerical solutions (which behoves the CFD solver, not the RL algorithm). Consequently, we would not anticipate \red{any additional numerical developments before} tackling a 3-D turbulent flow with the same CFD environment, even with a more sophisticated turbulence modeling (since the built-in small-scale component of the VMS solution also acts as an implicit LES).
Despite these achievements, further development, characterization and fine-tuning are needed to 
consolidate the acquired knowledge, whether it be via an 
improved balance between exploration and exploitation to deal with steep global maxima (for instance using Trust Region-Guided PPO, as it effectively encourages the policy to explore more on the potential valuable actions, no matter whether they were preferred by the previous policies or not), via non-normal probability density functions to deal with multiple global maxima, or via coupling with a surrogate model trained on-the-fly. 



\appendix

\section{A quick survey of adjoint-based optimization}\label{sec:appadj}

We briefly review here the various adjoint frameworks used in section~\ref{sec:olc:sub:control} 
for systematic validation purposes of the PPO-1 results. The starting point is a so-called uncontrolled solution $(\uu,p)$ to the non-linear equations of motion (Navier--Stokes, unless specified otherwise) forced by a momentum source $\ff$ and a velocity $\uu_w$ distributed over all solid surfaces $\Gamma_w$ in the computational domain (although it is possible to restrict to a subset).

\subsection{Baseline adjoint method} \label{sec:appadj:sub:base}

The adjoint method computes the change in drag induced by small variations ($\delta\ff, \delta\uu_w$) of these control parameters as
\bal
\delta \overline{c_x}=\int_\Omega \overline{\uu^\dag\cdot\delta\ff}\,\textrm{d}s+\int_{\Gamma_w} \overline{(\bsigma^\dag(-p^\dag,\uu^\dag)\cdot\nn)\cdot\delta\uu_w}\,\textrm{d}l\,,\label{eq:deltacdadj}
\eal
where $\nn$ is the unit outward normal to $\Gamma_w$ annd we note $\bsigma^\dag(-p^\dag,\uu^\dag)=p^\dag\op{I}+\frac{1}{\rey}\nabla\uu^\dag$. Finally, $(\uu^\dag,p^\dag)$ are adjoint velocity and pressure fields solution to
\bal
\nabla\cdot \uu^\dag=0\,,\qquad
-(\partial_t\uu^\dag+\nabla\uu^\dag \cdot \uu) +\nabla \uu^T \cdot \uu^\dag
+\nabla\cdot\bsigma^\dag(-p^\dag,\uu^\dag)=\00\,,\label{eq:nsadj2}
\eal
forced at $\Gamma_w$ by a velocity equal to twice the ambient velocity (the factor of 2 stems from the definition of dynamic pressure), as obtained multiplying $\uu^\dag$ and $p^\dag$ onto the linearized momentum and continuity equations, using the divergence theorem to integrate by parts over the computational domain, and integrating in time over the span of the simulation. In essence, this
amounts to computing the drag of the control induced disturbance modeled after the forced, linearized Navier--Stokes equations, without ever computing the disturbance itself.

A typical implementation consists of two sequential numerical simulations (for the uncontrolled and adjoint solutions, respectively) plus a series of vector dot products, to give the drag variation at each grid point. This is simple on paper, but the method has some limitations :
\smallskip

\paragraph*{- \!\!\!} the adjoint equations are problem-specific and must be derived and implemented manually on a case-by-case basis.
\smallskip

\paragraph*{- \!\!\!} the cost is marginal in steady flow regimes, because the time-independence of the uncontrolled solution makes the adjoint problem purely linear. Otherwise, the entire time history of uncontrolled solutions must be available at every adjoint time step because of the reversal of space-time directionality; see the minus sign ahead of the material derivative term in eqs.~\myrefeq{eq:nsadj2}. \red{This is very demanding in terms of storage (the repeated IO also increases the computational burden compared to a classical CFD run with identical simulation parameters) but these issues can be mitigated 
using checkpointing~\cite{Griewank2000} and high-order time-integration and interpolation schemes~\cite{Tsitouras2011}.}
\smallskip

\paragraph*{- \!\!\!} not all cost functions are admissible due to the need for consistent adjoint boundary conditions, \red{although this can be overcome with augmented Lagrangian methods based on auxiliary boundary equations~\cite{Arian1999}.}
\smallskip

\paragraph*{- \!\!\!} applicability to high-fidelity turbulence modeling is uncertain because the noise-induced sensitivity to initial conditions (the ``butterfly effect'') is expected to yield exponentially diverging solutions if the length of the adjoint simulation exceeds the predictability time scale. Possible solutions include averaging over a large number of ensemble calculations~\cite{Lea00} (which increases significantly the computational cost and decreases the attractiveness of the method) or invoking sophisticated shadowing and space-split techniques sampling on selected flow trajectories~\cite{wang13jcp,Chandramoorthy2019} (which comes at the cost of ease of implementation).
Moreover, \red{the literature somehow oddly
reports} several cases of turbulent adjoint solution blowing up in 2-D~\cite{Barth2010,Nazarov2012} and 3-D~\cite{wang13}, but also several instances in 3-D where no blow-up is observed~\cite{Hoffman2005,Hoffman2006,Jansson2011}. 
\smallskip

\paragraph*{- \!\!\!} applicability to RANS simulations is conversely generally acknowledged. However, discarding the linearization and adjointization of even the simplest turbulence models (using the so-called frozen eddy-viscosity approximation) to avoid massive debugging and validation efforts has somehow become standard lore, even though completeness and exactness are required to ensure numerical accuracy and avoid diverging adjoint solutions due to error propagation and amplification.

\subsection{Frozen Reynolds stresses approximation}\label{sec:appadj:sub:mean}

A simple adjoint formalism has been proposed in~\cite{meli14} to provide insight into the reliability of adjoint-based predictions in practical situations where no complete history of time and space-accurate solutions is available. The approach is closely related to existing studies considering the mean flow an admissible solution for linear stability analysis, as it simply dismisses the way the control-induced modification to the fluctuating uncontrolled solution feeds back onto the mean (hence the frozen Reynolds stress moniker to echo the above frozen eddy viscosity). In doing so, \myrefeq{eq:deltacdadj} can be shown to reduce to 
\bal
\delta \overline{\overline{c_x}}=\int_\Omega \overline{\overline{\uu^\dag}}\cdot\overline{\delta\ff}\,\textrm{d}s+
\int_{\Gamma_s} (\bsigma^\dag(\overline{\overline{p^\dag}},\overline{\overline{\uu^\dag}})\cdot\nn)\cdot\overline{\delta\uu_w}\,\textrm{d}l\,,
\eal
where the double overline denotes approximations to the true time-averaged quantities, and 
the adjoint velocity and pressure fields are solution to 
\bal
\nabla\cdot \overline{\overline{\uu^\dag}}=0\,,\qquad
-\nabla\overline{\overline{\uu^\dag}} \cdot \overline{\uu} +\nabla \overline{\uu}^T \cdot \overline{\overline{\uu^\dag}}
+\nabla\cdot\bsigma(-\overline{\overline{p^\dag}},\overline{\overline{\uu^\dag}})=\00\,,\label{eq:nsadj2mean}
\eal
forced at $\Gamma_w$ by the same velocity equal to twice the ambient velocity.
The strength of the approach lies in the fact that once the mean uncontrolled solution is known, computing the approximated adjoint solution merely requires solving a single linear problem.
Accuracy must be assessed on a case-by-case basis, but the computational and storage costs of doing so are marginal, and the approach carries over to any turbulence modeling method under the frozen viscosity assumption.

\section*{Acknowledgements} 
This work is supported by the Carnot M.I.N.E.S. Institute through the M.I.N.D.S. project.

\bibliography{biblist}

\end{document}